\documentclass{article}      

\usepackage{amsmath,amssymb,url}

\newtheorem{theorem}{Theorem}[section]

\newtheorem{prop}[theorem]{Proposition}

\newtheorem{lem}[theorem]{Lemma}
\newtheorem{rem}[theorem]{Remark}
\numberwithin{equation}{section}

\newcommand{\qed}{\Box}

\begin{document}

\title{Second order term 
of cover time for planar simple random walk}

\author{Yoshihiro Abe \footnote{Department of Mathematics and Informatics,
Chiba University, 1-33 Yayoi-cho Inage-ku Chiba-shi Chiba 263-8522 JAPAN;
              yosihiro@math.s.chiba-u.ac.jp}}

\date{}

\maketitle

\begin{abstract}
We consider the cover time for a simple random walk on the two-dimensional discrete
torus of side length $n$. Dembo, Peres, Rosen, and Zeitouni 
[Ann. Math. 160:433-464, 2004]
identified the leading term in
the asymptotics for the cover time as $n$ goes to infinity. 
In this paper, we study the exact second
order term. This is a discrete analogue of the work on the cover time for planar Brownian
motion by Belius and Kistler 
[Probab. Theory Relat Fields. 167:461-552, 2017].
\end{abstract}

\section{Introduction} \label{intro}
The present paper focuses on
the cover time for the simple random walk (SRW) on the two-dimensional discrete torus
$\mathbb{Z}_n^2 := (\mathbb{Z} / n \mathbb{Z})^2$
with the origin $o$.
Let $S = (S_i, i \geq 0, P_x, x \in \mathbb{Z}_n^2)$
be the discrete-time SRW on $\mathbb{Z}_n^2$.
For each $A \subset \mathbb{Z}_n^2$, 
let $H_A$ be the hitting time of $A$ defined by
\begin{equation} \label{eq:hitting_time}
H_A := \min \{i \geq 0 : S_i \in A \}.
\end{equation}
For each $x \in \mathbb{Z}_n^2$, we will write $H_x$ when $A = \{x \}$.
The cover time for the SRW on $\mathbb{Z}_n^2$
is given by
\begin{equation} \label{eq:cover_time}
\tau_{\text{cov}}^n := \max_{x \in \mathbb{Z}_n^2} H_x.
\end{equation}
Dembo, Peres, Rosen, and Zeitouni \cite{DPRZ1} 
obtained the exact leading term of $\tau_{\text{cov}}^n$;
they showed that $\tau_{\text{cov}}^n/n^2 (\log n)^2$
converges to $4/\pi$ in probability as $n \to \infty$.
Ding \cite{Di} improved the result by proving that
$\sqrt{\tau_{\text{cov}}^n/n^2} - \frac{2}{\sqrt{\pi}} \log n$
is of order $\log \log n$ with probability tending to $1$ as $n \to \infty$.
Related to the precise estimate of the cover time,
geometric properties of the level set
(so-called {\it late points})
of the form
$$\left\{x \in \mathbb{Z}_n^2 : H_x > \alpha \frac{4}{\pi} n^2 (\log n)^2 \right\},
~~\alpha \in (0, 1)$$
were investigated in \cite{DPRZ2}, \cite{Ok}.
Recently, 
large deviations of the cover time
and the local structure around an unvisited point have been studied 
via random interlacement techniques \cite{CGPV}, \cite{CPV}, \cite{CP}, \cite{Ro}.
For a general introduction to the cover time,
we refer the reader to \cite[Chapter 11]{MCMT}.
See also \cite{DLP} \cite{Z} for
close connections with the discrete Gaussian free field.

The aim of this paper is to obtain the exact second order term
of the cover time of $\mathbb{Z}_n^2$
which was conjectured by Belius and Kistler \cite{BeKi}.
Our main result is as follows:

\begin{theorem} \label{thm:main}
There exists $c \in (0, 1)$ such that
\begin{equation} \label{eq:main-statement}
2 \log n - \log \log n - (\log \log n)^c \leq \frac{\tau_{\text{cov}}^n}{\frac{2}{\pi} n^2 \log n}
\leq 2 \log n - \log \log n + (\log \log n)^c
\end{equation}
holds with $P_o$-probability tending to $1$ as $n \to \infty$.
\end{theorem}

\begin{rem} \label{rem:bk-result}
Belius and Kistler \cite{BeKi} have already
established a similar estimate for the $\varepsilon$-cover time 
for the standard
Brownian motion on $(\mathbb{R}/\mathbb{Z})^2$.
In the proof of Theorem \ref{thm:main},
we basically follow techniques
developed in \cite{BeKi},
but it is not straightforward to apply them.
The technical difficulties in our discrete setting
come from lack of rotationally invariance of SRW;
it only has an approximately rotationally invariance, 
so we have to control the errors very carefully.
\end{rem}

\begin{rem}
It is natural to ask whether 
$\frac{\tau_{\text{cov}}^n}{\frac{2}{\pi} n^2 \log n} - (2 \log n - \log \log n)$
is tight or not.
Recently, Belius, Rosen, and Zeitouni \cite{BRZ2}
proved that the cover time for the Brownian motion on the two dimensional sphere 
is tight.
\end{rem}

Let us describe the outline of the paper. 
In Section \ref{sec:2}, we restate Theorem \ref{thm:main}
in terms of the number of excursions around each vertex,
which is called traversal process.
We will recall basic estimates on transition probabilities among circles
(Lemma \ref{lem:circle-to-circle-probab}).
Applying this, we will approximate
the law of the traversal process by that of a (critical) Galton-Watson process
in Lemma \ref{lem:transfer-lemma}
which provides explicit bounds on the approximation errors.
As was revealed by Belius and Kistler \cite{BeKi},
so-called ''barrier estimates" for the Galton-Watson process
play an important role in studying the traversal process.
Recently, Belius, Rosen, and Zeitouni \cite{BRZ1} gave quite general barrier estimates
for the Galton-Watson process. 
We will recall their result in Lemma \ref{lem:GW-process-barrier-estimate}.
In Section \ref{sec:upper-bound}, we obtain the upper bound of the cover time.
In Lemmas \ref{lem:control-above} and \ref{lem:control-below},
we will prove that the traversal process crosses neither the curve
$i \mapsto a_n^+ (i)$ nor the curve
$i \mapsto a_n^- (i)$
(see (\ref{eq:control-above-upper-curve}) and
(\ref{eq:control-below-lower-curve}) for the definitions).
Then, we impose the restriction so that
the traversal process stays in the tube $[a_n^- (\cdot), a_n^+ (\cdot)]$.
The upper envelope $a_n^+(\cdot)$ enables us to apply
the transfer lemma (Lemma \ref{lem:transfer-lemma})
and to make the approximation errors of order $1$. 
The lower envelope $a_n^- (\cdot)$ enables us to apply
the barrier estimate (Lemma \ref{lem:GW-process-barrier-estimate}).
In Section \ref{sec:lower-bound}, we obtain the lower bound of the cover time.
For each $x \in \mathbb{Z}_n^2$,
we will set the event $A_n (x)$
which says that $x$ is an unvisited point
and the traversal process corresponding to $x$ stays in the tube
$[b_n^- (\cdot), b_n^+ (\cdot)]$
(see (\ref{eq:event-lower-bound})--(\ref{eq:lower-bound-upper-control-curves}) for the definitions
of $A_n(x)$, $b_n^{\pm} (\cdot)$)
and apply the second moment method
to $Z_n := \sum_{x \in \mathbb{Z}_n^2} 1_{A_n (x)}$.
The most difficult part of the proof is to obtain the upper bound of $E_o [(Z_n)^2]$.
As is often the case with the study of the two dimensional cover time
(e.g. \cite{DPRZ1}, \cite{DPRZ2}, \cite{BeKi}),
in order to estimate $P_o [A_n (x) \cap A_n (y)]$,
we will consider branching structure of the two traversal processes
corresponding to $x$ and $y$.
We will take a decreasing sequence of radii $(r_k)_{k \geq 1}$
and set
the level $k = \ell (x, y)$ at which
the balls $B(x, r_k)$ and $B(y, r_k)$ of radius $r_k$ 
centered at $x$ and $y$ become disjoint.
See (\ref{eq:radii}) and (\ref{eq:branching-level})
for the precise definitions of $r_k$ and $\ell (x, y)$.
Intuitively speaking, for $k < \ell (x, y)$,
the two balls $B(x, r_k)$ and $B(y, r_k)$ overlap
and the traversal processes in the balls are heavily correlated,
while, for $k \ge \ell (x, y)$,
the two balls $B(x, r_k)$ and $B(y, r_k)$ are disjoint
and the traversal processes in the balls should be almost independent.
Thus,  
as in Section \ref{sec:upper-bound},
we can apply the transfer lemma (Lemma \ref{lem:transfer-lemma}) and the barrier estimate
(Lemma \ref{lem:GW-process-barrier-estimate})
to each of the traversal processes
in $B(x, r_k)$ and $B(y, r_k)$ for $k \ge \ell (x, y)$.
We deal with six cases
(see (\ref{eq:window}), (\ref{eq:length}), (\ref{eq:d-function})
for the definitions of $w_n$, $L_n$, $d_n (s)$):
\begin{equation*}
\begin{aligned}
&(1) \ell (x, y) = 0,\,\,\,  (2) 1 \leq \ell (x, y) \leq w_n + 2,\,\,\, 
(3) w_n + 3 \leq \ell (x, y) \leq d_n (\xi),\\
&(4) d_n (\xi) < \ell (x, y) 
\leq \lceil (1 - \varepsilon) L_n \rceil, \,\,\,
(5) \lceil (1 - \varepsilon) L_n \rceil < \ell (x, y) \leq L_n - w_n - 1,\\
&(6) \ell(x, y) \geq L_n - w_n
\end{aligned}
\end{equation*}
for some $0 < \varepsilon < 1$ and $\xi > 0$.
We will show that
the term corresponding to $x, y$ in the case (1) forms the main part
of the estimate on $E_o [(Z_n)^2]$ (Lemma \ref{lem:2nd-moment-case1})
and that other terms corresponding to $x, y$ in cases (2)--(6)
are negligible (Lemmas \ref{lem:2nd-moment-case2}, \ref{lem:2nd-moment-case3},
\ref{lem:2nd-moment-case4}, \ref{lem:2nd-moment-case5}, (\ref{eq:main-lower-bound-pf-step1})).
Main ingredients in the estimate of $P_o [A_n (x) \cap A_n (y)]$  are decoupling inequalities
which break dependence of excursions around $x$ and $y$
(Lemma \ref{lem:decomposition-trajectories} 
and (\ref{eq:lem-2nd-moment-case4-pf-3})).
In the cases (1)--(3), we will use Lemma \ref{lem:decomposition-trajectories} which claims that
trajectories of SRW left in disjoint regions
are almost independent.
Thanks to Lemma \ref{lem:decomposition-trajectories},
in the case (1), $A_n (x)$ and $A_n (y)$ are almost independent 
(see (\ref{eq:lem-2nd-moment-case1-1-almost-sure-independence}))
and, in the cases (2)--(3), 
the barrier conditions in $A_n(x)$ and $A_n (y)$ for $i \ge \ell (x, y)$ 
are almost independent (see e.g. Lemma \ref{lem:lem-2nd-moment-case3-pf-tech-1}).
In the case (4), we use
the equation (\ref{eq:lem-2nd-moment-case4-pf-3}) which
states that excursions inside a ball are independent of
excursions outside the ball except the initial and terminal points of the excursions.
The Harnack inequality (\ref{eq:uniform-estimate-harmonic-measure})
removes the dependence of the endpoints.
Thanks to this,
the barrier conditions in $A_n(x)$ and $A_n (y)$ for $i \ge \ell (x, y)$ 
are almost independent (Lemma \ref{lem:2nd-moment-case4-decoupling}).
In the case (5), thanks to the strong Markov property,
$A_n (x)$ and the barrier condition in $A_n (y)$ for $i \ge \ell (x, y)$ 
are almost independent (see (\ref{eq:2nd-moment-case5-pf-17})).
In the case (6),
since the number of such pairs
is small, we will just bound $P_o[A_n (x) \cap A_n(y)]$ by
$P[A_n (x)]$ (see (\ref{eq:main-lower-bound-pf-step1})).
In Section \ref{sec:appendix}, we study concentration estimates on excursion lengths.
In Sections \ref{sec:pf-transfer-lem}, \ref{sec:barrier-estimates},
we give proofs for the transfer lemma (Lemma \ref{lem:transfer-lemma})
and the barrier estimate (Lemma \ref{lem:GW-process-barrier-estimate}).

Throughout the paper, we will write $c, c^{\prime}, \dotsc$ 
to denote positive constants.
Values of $c, c^{\prime}, \dotsc$  will change from line to line.
We use $c_1, c_2, \dotsc$ to denote constants whose values are fixed within each argument.
Given sequences $(c_n)_{n \ge 1}$ and $(c_n^{\prime})_{n \ge 1}$, we write 
$c_n^{\prime} = O(c_n)$
if there exists a universal constant $C$ such that $|c_n^{\prime}/c_n| \leq C$ for all $n \geq 1$.
We will write $c_n \asymp c_n^{\prime}$ when $c_n = O (c_n^{\prime})$
and $c_n^{\prime} = O (c_n)$.
We let $o(1)$ denote a term with $o (1) \to 0$ as $n \to \infty$.
For any set $A$, $|A|$ denotes the cardinality of $A$.
For $a, b \in \mathbb{R}$, we set
$a \vee b := \max \{a, b \}$ and $a \wedge b := \min \{a, b \}$.
Let $d(\cdot, \cdot)$ be the $\ell^2$-distance in $\mathbb{Z}_n^2$.
Set $B(x, r) := \{y \in \mathbb{Z}_n^2 : d(x, y) < r \}.$
For each $A \subset \mathbb{Z}_n^2$ with $|A| \geq 2$, 
we define its boundary by $\partial A := \{y \in \mathbb{Z}_n^2 : 
y \in \mathbb{Z}_n^2 \backslash A, ~d(x, y) = 1 ~\text{for some}~x \in A \}$.
When $A = \{x \}$, we set $\partial A := \{x \}$.

\section{Excursion counts} \label{sec:2}
In this section, we reformulate Theorem \ref{thm:main}
in terms of the number of excursions between concentric circles.
Recall that $(S_i)_{i \geq 0}$ is a SRW
on $\mathbb{Z}_n^2$.
Let us give a sequence of random times as follows:
For each $x \in \mathbb{Z}_n^2$ and $0 < r < R < \frac{n}{2}$, set
\begin{align} \label{eq:stopping-times}
&R_1 (x, R, r) := H_{\partial B (x, r)}, \notag \\
&D_k (x, R, r) := \min \{i > R_k (x, R, r) : S_i \in \partial B (x, R) \}, \notag \\
&R_{k+1} (x, R, r) := \min \{i > D_k (x, R, r) : S_i \in \partial B (x, r) \},~~~k \geq 1.
\end{align}
We will call the path 
$(S_{R_k (x, R, r)}, \dotsc, S_{D_k (x, R, r)})$ 
($(S_{D_k (x, R, r)}, \dotsc, S_{R_{k + 1} (x, R, r)})$, respectively)
excursion  from $\partial B (x, r)$ to $\partial B (x, R)$
(from $\partial B (x, R)$ to $\partial B (x, r)$, respectively).
Fix sufficiently small $\delta \in (0, \frac{1}{2})$.
Pick $\gamma \in (0, 1)$ close enough to $1$.
Then, take $\alpha, \beta \in (0, 1)$ close enough to zero so that the following holds:
\begin{equation} \label{eq:assumption-parameter}
2 \gamma - 2 \beta - \alpha > 1,~
(1 - 2 \delta) \beta > \alpha,~ \alpha + \beta > \frac{1}{2} + \delta - \alpha \delta.
\end{equation} 
We will consider concentric circles of radii
\begin{equation} \label{eq:radii}
r_k = r_k (n) := (\ell_n)^{L_n - k},~~~k = 0, \cdots, L_n,
\end{equation}
where
\begin{equation} \label{eq:width}
\ell_n = \ell_n (\alpha) := \exp \{(\log \log n)^{\alpha} \},
\end{equation}
\begin{equation} \label{eq:window}
w_n = w_n (\beta) := (\log \log n)^{\beta},
\end{equation}
\begin{equation} \label{eq:length}
L_n = L_n (\alpha, \beta, c_{\star}):= \left \lfloor \frac{\log n}{\log (\ell_n)} - c_{\star} w_n 
\right \rfloor.
\end{equation}
Here $c_{\star}$ is a sufficiently large positive constant.
Set
\begin{align}
&m_n^{+} := \left \lfloor \left(1 - \frac{\log \log n}{2 \log n} + \frac{s_n}{\log n} \right) \frac{2 (\log n)^2}{\log (\ell_n)}
\right \rfloor, \label{eq:number-crossing1} \\
&m_n^{-} := \left \lceil \left(1 - \frac{\log \log n}{2 \log n} - \frac{s_n}{\log n} \right) \frac{2 (\log n)^2}{\log (\ell_n)}
\right \rceil, \label{eq:number-crossing2}
\end{align}
where
\begin{equation} \label{eq:fluctuation}
s_n = s_n (\gamma):= (\log \log n)^{\gamma}.
\end{equation} 
The following indicates that $m_n^{+}$
and $m_n^{-}$ are approximations of the number of excursions
from $\partial B (x, r_1)$ to $\partial B (x, r_0)$
up to the cover time.
\begin{prop} \label{prop:traversal-process}
The events that
\begin{equation} \label{eq:upper-main-prop}
H_x \leq D_{m_n^{+}} (x, r_0, r_1),~~\forall x \in \mathbb{Z}_n^2,
\end{equation}
\begin{equation} \label{eq:lower-main-prop}
H_x \geq D_{m_n^{-}} (x, r_0, r_1),~~\exists x \in \mathbb{Z}_n^2
\end{equation}
hold with $P_o$-probability tending to $1$ as $n \to \infty$.
\end{prop}
We will give the proof of (\ref{eq:upper-main-prop})
(respectively, (\ref{eq:lower-main-prop}))
in Section \ref{sec:upper-bound}
(respectively, in Section \ref{sec:lower-bound}).

We have uniform controls of $D_{m_n^{+}} (x, r_0, r_1)$ and $D_{m_n^{-}} (x, r_0, r_1)$
in $x \in \mathbb{Z}_n^2$:
\begin{prop} \label{prop:time-number}
The events that
\begin{equation} \label{eq:upper-time}
D_{m_n^{+}} (x, r_0, r_1) \leq \frac{4}{\pi} n^2 (\log n)^2 \left(1 - \frac{\log \log n}{2 \log n} + 
\frac{2 s_n}{\log n} \right),~ \forall x \in \mathbb{Z}_n^2,
\end{equation}
\begin{equation} \label{eq:lower-time}
D_{m_n^{-}} (x, r_0, r_1) \geq \frac{4}{\pi} n^2 (\log n)^2 \left(1 - \frac{\log \log n}{2 \log n} -
\frac{2 s_n}{\log n} \right),~ \forall x \in \mathbb{Z}_n^2
\end{equation}
hold with $P_o$-probability tending to $1$ as $n \to \infty$.
\end{prop}
We will provide the proof of Proposition \ref{prop:time-number}
in Section \ref{sec:appendix}.

We will use the following basic estimates on the SRW on $\mathbb{Z}^2$.
See, for example, \cite[Proposition 1.6.7, Exercise 1.6.8]{La} for the proof.
\begin{lem} \label{lem:circle-to-circle-probab}
(i) There exists $c_1 > 0$ such that for all $0 < r < R < \frac{n}{2}$ and
$x, y \in \mathbb{Z}_n^2$ with $r < d(x, y) < R,$
$$\frac{\log (\frac{R}{d(x, y)}) - \frac{c_1}{r}}{\log (\frac{R}{r})} 
\le P_y [H_{\partial B (x, r)} < H_{\partial B (x, R)}]
\le \frac{\log (\frac{R}{d(x, y)}) + \frac{c_1}{r}}{\log (\frac{R}{r})}.$$
(ii) There exists $c_2 > 0$ such that
for all $0 < R < \frac{n}{2}$
and $x, y \in \mathbb{Z}_n^2$ with $0 < d (x, y) < R$,
$$\frac{\log (\frac{R}{d (x, y)}) - \frac{c_2}{d (x, y)} - \frac{c_2}{\log R}}{\log R}
\leq P_y [H_x < H_{\partial B (x, R)}]
\leq \frac{\log (\frac{R}{d (x, y)}) + \frac{c_2}{d (x, y)} + \frac{c_2}{\log R}}{\log R}.$$
\end{lem}
It is convenient to give specific notation to the probabilities
in Lemma \ref{lem:circle-to-circle-probab}.
Given a decreasing sequence of radii $\overline{R} = (R_i)_{0 \leq i \leq L}$
with $R_L = 1$,
for each $0 \leq i_1 < i_2 < i_3 < L$, set
\begin{align} \label{eq:prob-circle-to-circle-1}
&p_{i_2, i_3}^{i_1, \pm} (\overline{R})
:= \frac{\log \left(\frac{R_{i_1}}{R_{i_2}} \right) \pm 
\frac{c_{\ref{lem:circle-to-circle-probab}.1}}{R_{i_3}}}
{\log \left(\frac{R_{i_1}}{R_{i_3}} \right)} \notag \\
&p_{i_2, i_1}^{i_3, \pm} (\overline{R}) := 1 - p_{i_2, i_3}^{i_1, \mp} (\overline{R}),
\end{align}
where $c_{\ref{lem:circle-to-circle-probab}.1}$
is the constant $c_1$ in Lemma \ref{lem:circle-to-circle-probab}(i).
When $i_3 = L$, set
\begin{align} \label{eq:prob-circle-to-circle-2}
&p_{i_2, L}^{i_1, \pm} (\overline{R})
:= \frac{\log \left(\frac{R_{i_1}}{R_{i_2}} \right) \pm c_{\ref{lem:circle-to-circle-probab}.2}
\left(\frac{1}{R_{i_2}} + \frac{1}{\log (R_{i_1})} \right)}
{\log \left(R_{i_1} \right)} \notag \\
&p_{i_2, i_1}^{L, \pm} (\overline{R}) := 1 - p_{i_2, L}^{i_1, \mp} (\overline{R}),
\end{align}
where $c_{\ref{lem:circle-to-circle-probab}.2}$ is the constant
$c_2$ in Lemma \ref{lem:circle-to-circle-probab}(ii).
We will compare these probabilities
with corresponding transition probabilities on 
$\{0, 1, \cdots, L \}$.
More precisely, we will approximate 
$p_{i_2, i_3}^{i_1, \pm} (\overline{R})$
by $\frac{i_2 - i_1}{i_3 - i_1}$,
which is the probability of the event that
SRW on $\{0, \dotsc, L \}$ starting at $i_2$ reaches $i_3$
before visiting $i_1$. 
In the comparison, we will focus on ``errors" which are defined as follows:
for each $0 \leq i_1 < i_2 < i_3 \leq L$,
\begin{equation} \label{eq:error}
\Delta_{i_2, i_3}^{i_1, \pm} =
\Delta_{i_2, i_3}^{i_1, \pm} (\overline{R})
:= \frac{p_{i_2, i_3}^{i_1, \pm} (\overline{R})}{\frac{i_2 - i_1}{i_3 - i_1}},
~~\Delta_{i_2, i_1}^{i_3, \pm} = \Delta_{i_2, i_1}^{i_3, \pm} (\overline{R})
:= \frac{p_{i_2, i_1}^{i_3, \pm} (\overline{R})}{\frac{i_3 - i_2}{i_3 - i_1}}.
\end{equation}
As is often the case with the study of the cover time in two dimensions,
we will analyze traversal processes rather than hitting times.
Recall the definition of ``excursion" from below (\ref{eq:stopping-times}).
Fix $x \in \mathbb{Z}_n^2$ and $m \in \mathbb{N}$.
For each $0 \leq i < j \leq L$,
we define the traversal process by
\begin{align} \label{eq:traversal-general}
T_{j \to i}^{x, m} (\overline{R})
:= ~&\text{the number of traversals from
$\partial B (x, R_j)$ to $\partial B (x, R_i)$} \notag \\
&\text{by $m$ excursions from 
$\partial B (x, R_1)$ to $\partial B (x, R_0)$}.
\end{align} 
We will write
\begin{equation} \label{eq:traversal-short}
T_i^{x, m} (\overline{R}) := T_{i + 1 \to i}^{x, m} (\overline{R}).
\end{equation}
Similarly, for each $0 \leq k \leq i \leq L$, we set
\begin{align} \label{eq:traversal-process-intermediate}
T_i^{k, x, m} (\overline{R})
:= ~&\text{the number of traversals from
$\partial B (x, R_{i+1})$ to $\partial B (x, R_i)$} \notag \\
&\text{by $m$ excursions from 
$\partial B (x, R_{k+1})$ to $\partial B (x, R_k)$}.
\end{align} 
Let $P_m^{\text{GW}}$ be the law of the Galton-Watson process $(T_i)_{i \geq 0}$
with $T_0 = m$ and geometric offspring distribution of parameter $1/2$.
Then, the following lemma gives explicit bounds on errors in the approximation of 
the law of the traversal process by that of the Galton-Watson process.
Since each traversal makes an error of the form (\ref{eq:error}),
the approximation errors are given by the products of such errors.
\begin{lem} \label{lem:transfer-lemma}
(Transfer lemma)
Fix any decreasing sequence of radii $\overline{R} = (R_i)_{i = 0}^L$
with $R_L = 1$, $x, y \in \mathbb{Z}_n^2$,
and $m, L \in \mathbb{N}$. \\
(i) For any $k, \widetilde{k} \in \{0, \dotsc, L - 1 \}$
with $k \leq L - \widetilde{k} - 1$,
$m_i \geq 0$, $k \leq i \leq L - \widetilde{k} - 1$
with
$P_m^{\text{GW}} [\{T_i = m_i, ~k \leq \forall i \leq L - \widetilde{k} - 1\}
\cap \{T_{L - 1} = 0 \}] \neq 0$,
\begin{align} \label{eq:transfer-lemma-statement-1}
&\frac{
P_y \left[\left\{T_i^{x, m} (\overline{R}) = m_i, ~k \leq \forall i \leq L - \widetilde{k} - 1 \right \}
\cap \left\{T_{L - 1}^{x, m} (\overline{R}) = 0 \right \} \right]}
{P_m^{\text{GW}} \left[\left\{T_i = m_i, ~k \leq \forall i \leq L - \widetilde{k} - 1 \right \}
\cap \left\{T_{L - 1} = 0 \right \} \right]} \notag \\
&\in \left[\Delta_1^- \cdot \Delta_{\star}^-,~\Delta_1^+ \cdot \Delta_{\star}^+ \right],
\end{align}
where $\Delta_1^+$ and $\Delta_{\star}^+$ are defined by
\begin{align} \label{eq:transfer-lemma-main-error-1}
\Delta_1^{+} &:=
\left(\Delta_{1, 0}^{k + 1, +} \vee \Delta_{1, k + 1}^{0, +}  \right)^m
\left(\Delta_{k, 0}^{k + 1, +}  \vee \Delta_{k, k + 1}^{0, +}  \right)^{m_k}
\notag \\
&\times \prod_{i = k + 1}^{L - \widetilde{k} - 1}
\left(\Delta_{i, i - 1}^{i + 1, +} \vee \Delta_{i, i + 1}^{i - 1, +} 
\right)^{m_{i - 1} + m_i},
\end{align}
and 
$\Delta_{\star}^{+} := \left(\Delta_{L - \widetilde{k}, L - \widetilde{k} - 1}^{L, +} 
\right)^{m_{L - \widetilde{k} - 1}}$.
$\Delta_1^-$ and $\Delta_{\star}^-$
are defined by replacing ``$+$" and ``$\vee$'' 
in the definitions of $\Delta_1^+$ and $\Delta_{\star}^+$
with ``$-$'' and ``$\wedge$''.
Moreover
(\ref{eq:transfer-lemma-statement-1}),
with
$\{T_{L - 1}^{x, m} (\overline{R}) = 0 \}$
and $\{T_{L - 1} = 0 \}$ removed and
with $\Delta_{\star}^{\pm}$ replaced by $1$,
holds true. \\\\
(ii) For any $\ell \in \mathbb{N}$ and $1 \leq i \leq L - 1$,
\begin{equation} \label{eq:transfer-lemma-statement-2}
P_y \left[T_i^{x, m} (\overline{R}) \geq \ell,~T_{L - 1}^{x, m} (\overline{R}) = 0 \right]
\leq \Delta_2^+ \cdot
P_m^{\text{GW}} \left[T_i \geq \ell,~T_{L - 1} = 0 \right],
\end{equation}
where
\begin{equation} \label{eq:transfer-lemma-error-2}
\Delta_2^+
:= \left(\Delta_{1, 0}^{i + 1, +} \vee \Delta_{1, i + 1}^{0, +} \right)^m
\left(\Delta_{i, 0}^{i + 1, +} \vee \Delta_{i, i + 1}^{0, +} \right)^{\ell}
\left(\Delta_{i + 1, i}^{L, +} \right)^{\ell}
\Delta_{i + 1, 0}^{L, +}
\left(\Delta_{1, 0}^{L, +} \right)^m.
\end{equation} \\\\
(iii) For any $k, \widetilde{k} \in \{1, \dotsc, L - 1 \}$
with $k < L - \widetilde{k} - 1$,
$m_i \geq 0$, $k + 1 \leq i \leq L - \widetilde{k} - 1$,
$y \in \partial B (x, R_k)$,
\begin{equation} 
\label{eq:transfer-lemma-statement-3}
\begin{aligned}
&P_y \Biggl[
\bigcap_{i=k+1}^{L - \widetilde{k} - 1}
\left\{T_i^{k, x, m} (\overline{R}) = m_i \right\} \cap
\left\{T_{L - 1}^{k, x, m} (\overline{R}) = 0 \right\} \\
&\,\,\,\,\,\,\,\,\,\,\,\,\,\cap
\left\{R_m (x, R_k, R_{k + 1}) < H_{\partial B (x, R_0)} < R_{m + 1} (x, R_k, R_{k + 1}) \right\}
\Biggr] \\
&\leq \Delta_3^+ P_m^{\text{GW}}
\left[\bigcap_{i=1}^{L - k - \widetilde{k} - 1} \{T_i = m_{k+i}\} \cap \{T_{L - k - 1} = 0\} \right]
\cdot P [G = m],
\end{aligned}
\end{equation}
where $G$ is a geometric random variable with success probability
$\frac{1}{k + 1}$ and
$\Delta_3^+$ is defined by
\begin{equation} \label{eq:transfer-lemma-error3}
\prod_{i = k + 1}^{L - \widetilde{k} - 1}
\left(\Delta_{i, i - 1}^{i + 1, +} \vee \Delta_{i, i + 1}^{i - 1, +} \right)^{m_{i - 1} + m_i}
\left(\Delta_{L - \widetilde{k}, L - \widetilde{k} - 1}^{L, +} \right)^{m_{L - \widetilde{k} - 1}}
\left(\Delta_{k, k + 1}^{0, +} \right)^m \cdot \Delta_{k, 0}^{k + 1, +}.
\end{equation}
\end{lem}

We will give the proof of Lemma \ref{lem:transfer-lemma} in Section \ref{sec:pf-transfer-lem}.

In the proof of Proposition \ref{prop:traversal-process}, we heavily use so-called barrier estimates
for the Galton-Watson process.
For $a, b \in \mathbb{R}$ and $L \in \mathbb{N}$, set a linear barrier
\begin{equation*}
f_{a, b} (i; L) := a + (b - a) \frac{i}{L},~~~0 \leq i \leq L.
\end{equation*}
We will consider small perturbations of the linear barrier of the form
$f_{a, b} (i; L) \pm C (i_L)^c$ for some $c, C > 0$, where
$$i_L := i \wedge (L - i),~~~0 \leq i \leq L.$$
The following is a slightly modified version of barrier estimates for 
the Galton-Watson process
by Belius, Rosen, and Zeitouni \cite{BRZ1}:
\begin{lem} \label{lem:GW-process-barrier-estimate}
(Barrier estimate)
(i) For any $\delta, C \in (0, \infty)$, $\varepsilon \in (0, \frac{1}{2})$,
$\eta > 1$,
$a, b, x, y$ with
$\sqrt{2} \leq x, y \leq \eta L$, $\frac{x^2}{2} \in \mathbb{N}$,
$0 \leq a \leq x$, $0 \leq b \leq y$, $b \leq a$,
\begin{align} \label{eq:GW-process-barrier-upper-bound}
&P_{\frac{x^2}{2}}^{\text{GW}} \left[
\bigcap_{i=1}^{L-1} \left\{f_{a, b} (i; L) - C i_L^{\frac{1}{2} - \varepsilon} \leq \sqrt{2 T_i} \right\}
\cap \left\{\sqrt{2 T_L} \in [y, y + \delta] \right\}
\right] \notag \\
&\leq c_1 e^{c_2 \eta^2}
\sqrt{\frac{x}{y}} \frac{1}{\sqrt{L}}
e^{- \frac{(y - x)^2}{2L}} 
\frac{(\eta + x - a)(\eta + y - b)}{L},
\end{align}
where $c_1$ and $c_2$ are positive constants which
depend only on $\delta$, $C$, $\varepsilon$. \\\\
(ii) There exists $r_0 > 0$ such that
for any $C > 0$, $\eta > 1$, $\varepsilon \in (0, \frac{1}{2})$, $\mu \in (0, 1)$,
$4 r^{\frac{1}{2} + 2 \varepsilon} \leq \mu L \leq a \leq x \leq \eta L$,
$\frac{x^2}{2} \in \mathbb{N}$,
$C r^{\frac{1}{2} - \varepsilon} > \eta$, $L > 2r > r_0$,
\begin{align} \label{eq:GW-process-barrier-lower-bound}
&P_{\frac{x^2}{2}}^{\text{GW}} \Biggl[
\bigcap_{i=r}^{L-1-r}
\left\{f_{a, b} (i; L) + C i_L^{\frac{1}{2} - \varepsilon} \leq \sqrt{2 T_i}
\leq f_{x, 0} (i; L) + \widetilde{C} i_L^{\frac{1}{2} + \varepsilon} \right\} \notag \\
&\,\,\,\,\,\,\,\,\,\,\,\,\,\,\,\,\,\,\,\,\,\,\,\,\,\,\,\,\,\,\,\,\,\,\,\,\,\,\,\,\,\,
\,\,\,\,\,\,\,\,\,\,\,\,\,\,\,\,\,\,\,\,\,\,\,\,\,\,\,\,\,\,\,\,\,\,\,\,\,\,\,\,\,\,
\,\,\,\,\,\,\,\,\,\,\,\,\,\,\,\,\,\,\,\,\,\,\,\,\,\,\,\,\,\,\,\,\,\,\,\,\,\,\,\,\,\,
\cap \left\{\sqrt{2 T_{L - 1}} = 0 \right\}
\Biggr] \notag \\
&\geq c_3 \frac{r}{L - 2r} \left(1 - \frac{1}{L} \right)^{\frac{x^2}{2}},
\end{align}
where $c_3$ is a positive constant which depends only on $r_0$, $C$, $\varepsilon$, $\mu$.
\end{lem}
We will give the proof of Lemma \ref{lem:GW-process-barrier-estimate}
in Section \ref{sec:barrier-estimates}.
\begin{rem} \label{rem:barrier-estimate}
Unfortunately, we cannot directly apply
\cite[Theorem 1.1]{BRZ1} for the following reason:
Let $x, L, \eta$ be the constants in \cite[Theorem 1.1]{BRZ1}.
In our setting, we will typically take 
$x \asymp \frac{\log n}{\sqrt{\log (\ell_n)}}$
and $L \asymp \frac{\log n}{\log (\ell_n)}$.
Thus, we need to take $\eta \asymp \sqrt{\log (\ell_n)}$.
Since this $\eta$ depends heavily on $n$,
we need explicit information on
how the constant $c$ in \cite[Theorem 1.1]{BRZ1} depends on $\eta$.
\end{rem}

\section{Upper bound of cover time} \label{sec:upper-bound}
In this section, we prove (\ref{eq:upper-main-prop}).
We will use the same notation as in Section \ref{sec:2}.
To simplify notation, we will write
\begin{equation} \label{eq:abbreviation-stopping-time}
D_m^{x, i} := D_m (x, r_i, r_{i+1}),~R_m^{x, i} := R_m (x, r_i, r_{i+1}),~x \in \mathbb{Z}_n^2,~i \geq 0,
~m \in \mathbb{N}.
\end{equation}
For $m \in \mathbb{N}$, $x \in \mathbb{Z}_n^2$, $1 \leq i \leq L_n - 1$, set
\begin{equation} \label{eq:modified-ver-traversal-process}
\widetilde{T}_i^{x, m} := 
\max \{\ell \geq 1 : R_{\ell}^{x, i} \circ \theta_{H_{\partial B (x, r_1)}} < D_m^{x, 0} \},
\end{equation}
where $\theta_s, s \geq 0$ are shift operators.
We will work on the traversal processes
$(\widetilde{T}_k^{x, m_n^+})_{k \in \{0, 1, \cdots, L_n -1 \}}$, $x \in \mathbb{Z}_n^2$.
We first prove that the traversal process $\widetilde{T}^{x, m_n^+}$
stays in the tube $[a_n^- (\cdot), a_n^+ (\cdot)]$ with high probability
in Lemmas \ref{lem:control-above} and \ref{lem:control-below} 
(see (\ref{eq:control-above-upper-curve}) and
(\ref{eq:control-below-lower-curve}) for the definitions of $a_n^+ (\cdot)$ and $a_n^- (\cdot)$).
Using this together with the transfer lemma (Lemma \ref{lem:transfer-lemma})
and the barrier estimate (Lemma \ref{lem:GW-process-barrier-estimate}), 
we will prove (\ref{eq:upper-main-prop}).

We need control of the traversal process from above
to deal with the approximation errors which appear
in applying the transfer lemma.
\begin{lem} \label{lem:control-above}
There exists $\kappa_+ > 0$ such that
\begin{equation} \label{eq:control-above}
\lim_{n \to \infty}
P_o \left[
\bigcup_{x \in \mathbb{Z}_n^2} \bigcup_{i=1}^{L_n - 1 - \left \lceil \frac{2 \log \log n}{\log (\ell_n)} \right \rceil}
\left\{\widetilde{T}_i^{x, m_n^+}
\geq a_n^{+} (i),~
H_x > D_{m_n^+}^{x, 0} \right\}
\right] = 0,
\end{equation}
where for each $i \in \{1, \cdots, L_n - 1 \}$
\begin{equation} \label{eq:control-above-upper-curve}
a_n^{+} (i) := \left \lceil 
\left \{\sqrt{m_n^{+}} \left(1 - \frac{i}{L_n} \right) + \kappa_+ \sqrt{\frac{(i+1)(L_n - i)}{L_n + 1}} \sqrt{\log \log n}
\right \}^2 \right \rceil.
\end{equation}
\end{lem}
{\it Proof.}
To simplify notation, we will write
\begin{equation} \label{eq:d-function}
d_n (s) := \left\lceil \frac{s\log \log n}{\log (\ell_n)} \right\rceil,~~s > 0.
\end{equation}
Fix $x \in \mathbb{Z}_n^2$
and $i \in \{1, \cdots, L_n - d_n (2) - 1 \}$.
By the transfer lemma (Lemma \ref{lem:transfer-lemma}(ii))
with $\overline{R} = (r_i)_{i = 0}^{L_n}$,
we have
\begin{equation} \label{eq:lem-control-above-4-0}
P_o \left[\widetilde{T}_{i}^{x, m_n^+} \geq a_n^{+} (i),
H_x > D_{m_n^+}^{x, 0} \right]
\leq \Delta_2^+ P_{m_n^+}^{\text{GW}} \left[T_i \geq a_n^{+} (i), ~T_{L_n - 1} = 0 \right],
\end{equation}
where $\Delta_2^+$ is defined by (\ref{eq:transfer-lemma-error-2})
with $m = m_n^+$, $\ell = a_n^+ (i)$, $L = L_n$.
Recalling the definitions of $\Delta_{i_2, i_3}^{i_1, +}$, $\Delta_{i_2, i_1}^{i_3, +}$
from (\ref{eq:error}),
by a simple calculation,
we have
\begin{equation*}
\Delta_{1, 0}^{i+1, +} \vee \Delta_{1, i+1}^{0, +} 
= 1 + O \left(\frac{1}{(\log n)^2 \log (\ell_n)} \right),
\end{equation*}
\begin{equation*}
\Delta_{i+1, i}^{L_n, +}, ~\Delta_{i, 0}^{i+1, +} \vee \Delta_{i, i+1}^{0, +}
= 1 + O \left(\frac{1}{(L_n - i)^2 (\log (\ell_n))^2} \right),
\end{equation*}
\begin{equation*}
\Delta_{i+1, 0}^{L_n, +} = 1 + O\left(\frac{1}{(\log \log n)^2} \right),~
\Delta_{1, 0}^{L_n, +} = 1 + O \left(\frac{1}{(\log n)^2} \right).
\end{equation*}
Since $m_n^+ \asymp \frac{(\log n)^2}{\log (\ell_n)}$
and $a_n^+ (i) \leq c (L_n - i)^2 \log (\ell_n) + c (L_n - i) \log \log n$,
we have $\Delta_2^+ = 1 + o(1)$.
By the Markov property, the probability 
on the right-hand side of (\ref{eq:lem-control-above-4-0})
is equal to
\begin{equation} \label{eq:lem-control-above-5}
\sum_{m = a_n^+ (i)}^{\infty}
P_{m_n^+}^{\text{GW}} [T_i = m] P_m^{\text{GW}} [T_{L_n - i - 1} = 0].
\end{equation}
By \cite[(5.3)]{BRZ1},
the first probability of the $m$-th term in (\ref{eq:lem-control-above-5})
is bounded from above by
$e^{- (\sqrt{m} - \sqrt{m_n^{+}}~)^2/(i+1)}$.
By a simple calculation, 
the second probability of the $m$-th term in (\ref{eq:lem-control-above-5})
is equal to $(1- \frac{1}{L_n - i})^m$
and this is bounded from above by $e^{- \frac{m}{L_n - i}}$.
Thus, the right of (\ref{eq:lem-control-above-5})
is bounded from above by 
\begin{align} \label{eq:lem-control-above-6}
&\sum_{m = a_n^{+} (i)}^{\infty} e^{- \frac{m_n^+}{L_n + 1}}
\exp \left \{- \frac{L_n + 1}{(i+1)(L_n - i)} \left(\sqrt{m} - \frac{L_n - i}{L_n + 1} \sqrt{m_n^+} \right)^2 \right \} \notag \\
&\leq n^{- 2} (\log n) \sum_{m = a_n^+ (i)}^{4 m_n^+} e^{- (\kappa_+)^2 \log \log n}
+ n^{-2} (\log n) \sum_{m = 4 m_n^+ + 1}^{\infty} e^{- \frac{m}{4 L_n}} \notag \\
&\leq c_1 n^{-2} (\log n)^{- (\kappa_+)^2 + 3} + c_1 n^{- 4} (\log n)^3,
\end{align}
where we used the definition of $a_n^+ (i)$
and the inequality $\sqrt{m} - \frac{L_n - i}{L_n + 1} \sqrt{m_n^+} \geq \frac{\sqrt{m}}{2}$
for each $m > 4 m_n^+$
in the first inequality.
By this and the union bound,
taking $\kappa_+$ large enough,
we have the desired result.
$\qed$ \\\\
We need some control of the traversal process from below
to apply the barrier estimate.
\begin{lem} \label{lem:control-below}
There exists $\kappa_- > 0$ such that 
\begin{equation} \label{eq:control-below}
\lim_{n \to \infty} P_o \left[
\bigcup_{x \in \mathbb{Z}_n^2} \bigcup_{i=1}^{L_n-1}
\left\{\sqrt{\widetilde{T}_i^{x, m_n^+}} < \sqrt{m_n^{+}} \left(1- \frac{i}{L_n} \right)
- \kappa_- \frac{\log \log n}{\sqrt{\log (\ell_n)}} \right\}
\right] = 0.
\end{equation}
\end{lem} 
{\it Proof.}
Fix $i \in \{1, \cdots, L_n - 1 \}$.
Set $r_i^{+} := (1 + \frac{\sqrt{2}}{(\log n)^2}) r_i$, $r_i^{-} := (1 - \frac{\sqrt{2}}{(\log n)^2}) r_i$. 
Set
\begin{equation} \label{eq:kth-level}
F_i := \left \{ \left(k \cdot \widetilde{r}_i, 
\ell \cdot \widetilde{r}_i \right) : 
 k, \ell \in \left \{0, \cdots ,\left\lfloor \frac{n}{\widetilde{r}_i} \right\rfloor \right \} \right \},
\end{equation}
where $\widetilde{r}_i := \lfloor \frac{r_i}{\ell_n (\log n)^2} \rfloor$.
Fix any $x \in \mathbb{Z}_n^2$.
There exists $y \in F_i$ such that
$x \in \left(y + [0, \widetilde{r}_i]^2 \right) \cap \mathbb{Z}^2$ mod $n \mathbb{Z}^2$.
One can easily check the following:
$$B (y, r_{i+1}^{-}) \subset B (x, r_{i+1}) \subset B (x, r_i) \subset B (y, r_i^+),$$
$$B (x, r_1) \subset B (y, r_1^{+}) \subset B (y, r_0^{-}) \subset B (x, r_0).$$
Thus, we have
\begin{equation} \label{eq:comparison-traversal-processes}
\widehat{T}_i^{x, m_n^+} \leq \widetilde{T}_i^{x, m_n^+},
\end{equation}
where for $z \in \mathbb{Z}_n^2$, $m \in \mathbb{N}$, we set 
\begin{equation} \label{eq:modified-modified-ver-traversal-process}
\widehat{T}_i^{z, m} := 
\max \{\ell \geq 0 : R_{\ell} (z, r_i^+, r_{i + 1}^-) \circ \theta_{H_{\partial B (z, r_1^{+})}} 
< D_m (z, r_0^{-}, r_1^{+}) \}.
\end{equation} 
Thus, for sufficiently large $\kappa_-$, 
the probability in (\ref{eq:control-below}) is bounded from above by
\begin{equation} \label{eq:lem-control-below-1}
\sum_{i = 1}^{L_n - d_n (\sqrt{\kappa_-}) - 1}
\sum_{y \in F_i}
P_o \left[\widehat{T}_i^{y, m_n^+} < a_n^{-} (i) \right],
\end{equation}
where
\begin{equation} \label{eq:control-below-lower-curve}
a_n^{-} (i) := \left \lfloor \left \{\left(\sqrt{m_n^{+}} \left(1- \frac{i}{L_n} \right)
- \kappa_{-} \frac{\log \log n}{\sqrt{\log (\ell_n)}} \right) \vee 1 \right \}^2 \right \rfloor
\end{equation}
and we have used the fact that
for sufficiently large $\kappa_-$,
$\sqrt{m_n^{+}} (1 - \frac{i}{L_n}) - \kappa_- \frac{\log \log n}{\sqrt{\log (\ell_n)}} < 0$
for all $i \geq L_n - d_n (\sqrt{\kappa_-})$.
Fix $1 \leq i \leq  L_n - d_n (\sqrt{\kappa_{-}}) - 1$
and $y \in F_i$.
By the transfer lemma
(Lemma \ref{lem:transfer-lemma}(i))
with any sequence of radii $(R_k)_{k = 0}^{L_n}$
with $R_0 = r_0^-$, $R_1 = r_1^+$, $R_i = r_i^+$, $R_{i+1} = r_{i+1}^-$,
the $(i, y)$-th term in (\ref{eq:lem-control-below-1})
is bounded from above by
\begin{equation} \label{eq:lem-control-below-2}
(1 + o(1)) P_{m_n^+}^{\text{GW}} [T_i < a_n^- (i)].
\end{equation}
By \cite[(5.3)]{BRZ1},
the probability in (\ref{eq:lem-control-below-2}) is bounded from above by
\begin{equation} \label{eq:lem-control-below-3}
e^{- \frac{\left(\sqrt{m_n^+} - \sqrt{a_n^- (i)} \right)^2}{i + 1}}.
\end{equation}
By the definitions of $m_n^+$, $a_n^- (i)$,
and the condition
$i < L_n - d_n (\sqrt{\kappa_{-}})$,
(\ref{eq:lem-control-below-3})
is bounded from above by
$(\ell_n)^{- 2 (i+1)} (\log n)^{- 2 \sqrt{\kappa_-} + 1}$.
One can check that
$|F_i| \leq c_1 (\ell_n)^{2(i+1)} (\ell_n)^{2 c_{\star} w_n + 2} (\log n)^4$.
Thus, the sum in (\ref{eq:lem-control-below-1}) is bounded from above by
$c_2 e^{3 c_{\star} (\log \log n)^{\alpha + \beta} } (\log n)^{- 2 \sqrt{\kappa_-} + 6}$.
This goes to $0$ as $n \to \infty$ for sufficiently large $\kappa_- > 0$
since $\alpha + \beta < 1$ by the assumption (\ref{eq:assumption-parameter}).
$\qed$ \\\\

{\it Proof of (\ref{eq:upper-main-prop}).}
Recall the constants $\kappa_-$ and $\kappa_+$ from Lemmas \ref{lem:control-below} and
\ref{lem:control-above}.
Recall $a_n^{-} (k)$ and $a_n^{+} (k)$ from (\ref{eq:control-below-lower-curve})
and (\ref{eq:control-above-upper-curve}).
Recall the notation $d_n (\cdot)$ from (\ref{eq:d-function}).
To simplify notation, we set $L_n^{\star} := L_n - d_n (2 \kappa_{-}) - 1$. 
We have
\begin{align} \label{eq:pf-upper-main-prop-1}
&P_o \left[\bigcup_{x \in \mathbb{Z}_n^2} \left\{H_x > D_{m_n^{+}}^{x, 0} \right\} \right] \notag \\
\leq &P_o \left[
\bigcup_{x \in \mathbb{Z}_n^2} \bigcap_{i=1}^{L_n^{\star}}
\left\{
a_n^{-} (i) \leq \widetilde{T}_i^{x, m_n^+} \leq a_n^{+} (i)
\right\} \cap \left\{H_x > D_{m_n^{+}}^{x, 0} \right\}
\right] \notag \\
+ &P_o \left[
\bigcup_{x \in \mathbb{Z}_n^2}
\bigcup_{i=1}^{L_n - 1}
\left\{\sqrt{\widetilde{T}_i^{x, m_n^+}} < \sqrt{m_n^{+}} \left(1- \frac{i}{L_n} \right)
- \kappa_- \frac{\log \log n}{\sqrt{\log (\ell_n)}} \right\}
\right] \notag \\
+ &P_o \left[
\bigcup_{x \in \mathbb{Z}_n^2}
\bigcup_{i=1}^{L_n - d_n (2) - 1} 
\left\{\widetilde{T}_i^{x, m_n^+}
> a_n^{+} (i),~
H_x > D_{m_n^+}^{x, 0} \right\}
\right],
\end{align}
where we have used the fact that $\sqrt{m_n^{+}} \left(1- \frac{i}{L_n} \right)
- \kappa_- \frac{\log \log n}{\sqrt{\log (\ell_n)}} > 1$
for all $i \in \{1, \cdots, 
L_n^{\star} \}.$
The first term on the right of (\ref{eq:pf-upper-main-prop-1})
is the probability of the event that the traversal process $\widetilde{T}^{x, m_n^+}$
stays in the tube $[a_n^-, a_n^+]$ and equals zero at $L_n-1$ for some $x$.
The second term (resp. the third term) on the right of (\ref{eq:pf-upper-main-prop-1})
is the probability of the event that the traversal process
$\widetilde{T}^{x, m_n^+}$ hits the lower curve $a_n^-$ (resp. the upper curve $a_n^+$)
for some $x$.
By Lemmas \ref{lem:control-below} and \ref{lem:control-above}, the second and third terms
of (\ref{eq:pf-upper-main-prop-1}) go to $0$ as $n \to \infty$.
Thus, we only need to deal with the first term 
on the right-hand side of (\ref{eq:pf-upper-main-prop-1}).

Fix $x \in \mathbb{Z}_n^2$.
By the transfer lemma (Lemma \ref{lem:transfer-lemma}), we have
\begin{align} \label{eq:pf-upper-main-prop-2}
&P_o \left[
H_x > D_{m_n^{+}}^{x, 0},~
a_n^{-} (i) \leq \widetilde{T}_i^{x, m_n^+} \leq a_n^{+} (i),~
1 \leq \forall i 
\leq L_n^{\star}
\right]
 \notag \\
&\leq (1 + o(1))
P_{m_n^+}^{\text{GW}} \left[
a_n^{-} (i) \leq T_i \leq a_n^{+} (i),~
1 \leq \forall i \leq L_n^{\star},~
T_{L_n - 1} = 0
\right],
\end{align}
where $o(1) \to 0$ as $n \to \infty$
uniformly in $x$.
By conditioning on $T_{L_n^{\star}}$,
the probability on
the right of (\ref{eq:pf-upper-main-prop-2})
is bounded from above by
\begin{equation} 
\label{eq:pf-upper-main-prop-3}
\begin{aligned}
\sum_{m = a_n^- (L_n^{\star})}^{a_n^+ (L_n^{\star})}
&P_{m_n^+}^{\text{GW}}
\Biggl[
\bigcap_{i=1}^{L_n^{\star}-1} \left\{
\sqrt{T_i} \geq \sqrt{m_n^+} \left(1 - \frac{i}{L_n} \right)
- \frac{\kappa_{-} \log \log n}{\sqrt{\log (\ell_n)}} - 1 \right\} \\
&\,\,\,\,\,\,\,\,\,\,\,\,\,\,\,\,\,\,\,\,\,\,\,\,\,\,\,\,\,\,\,\,\,\,\,\,\,\,\,\,\,\,\,\,\,\,
\,\,\,\,\,\,\,\,\,\,\,\,\,\,\,\,\,\,\,\,\,\,\,\,\,\,\,\,\,\,\,\,\,\,\,\,\,\,\,\,\,\,\,\,\,\,
\cap \left\{T_{L_n^{\star}} = m \right \} \Biggr] \\
&\times P_m^{\text{GW}} [T_{d_n (2 \kappa_{-})} = 0].
\end{aligned}
\end{equation}
To estimate the first probability of the $m$-th term in (\ref{eq:pf-upper-main-prop-3}),
we will apply the barrier estimate (Lemma \ref{lem:GW-process-barrier-estimate}(i))
with $x = \sqrt{2 m_n^+}$, $y = \sqrt{2m}$, 
$a = \sqrt{2} (\sqrt{m_n^+} - \kappa_- \frac{\log \log n}{\sqrt{\log (\ell_n)}} - 1)$,
$b = \sqrt{2} \{\sqrt{m_n^+} (1 - \frac{L_n^{\star}}{L_n})
- \kappa_- \frac{\log \log n}{\sqrt{\log (\ell_n)}} - 1\}$
(we can take $\eta = c \sqrt{\log (\ell_n)}$ for some positive constant $c$.)
Then, the probability is bounded from above by
\begin{equation} \label{eq:pf-upper-main-prop-4}
c_1 \frac{(\ell_n)^{c_2}}{L_n} e^{- \frac{(\sqrt{m_n^+} - \sqrt{m})^2}{L_n^{\star}}}.
\end{equation}
By a simple calculation,
the second probability of the $m$-th term in (\ref{eq:pf-upper-main-prop-3})
is equal to $(1 - \frac{1}{d_n (2 \kappa_-) + 1})^m$ and this is bounded from above by
$e^{- \frac{m}{d_n (2 \kappa_-) + 1}}$.
The product of this and the exponential factor in (\ref{eq:pf-upper-main-prop-4})
is bounded from above by
$e^{- \frac{m_n^+}{L_n}}$.
By this, (\ref{eq:pf-upper-main-prop-2}) is bounded from above by
\begin{equation} \label{eq:pf-upper-main-prop-5}
c_3 (\ell_n)^{c_4} n^{-2} e^{- 2 s_n}.
\end{equation}

Therefore,
the first term on the right-hand side
of (\ref{eq:pf-upper-main-prop-1}) is bounded from above by
$c_5 (\ell_n)^{c_4} e^{- 2s_n}$
and this goes to $0$ as $n \to \infty$.
Therefore, the right of (\ref{eq:pf-upper-main-prop-1}) goes to $0$ as $n \to \infty$. $\qed$
\\\\

{\it Proof of the upper bound of Theorem \ref{thm:main} via (\ref{eq:upper-time}).}
(\ref{eq:upper-main-prop}) and (\ref{eq:upper-time})
immediately yield the upper bound. $\qed$

\section{Lower bound of cover time} \label{sec:lower-bound}
In this section, we prove (\ref{eq:lower-main-prop}).
Recall some notation from (\ref{eq:stopping-times})-(\ref{eq:fluctuation})
and (\ref{eq:abbreviation-stopping-time}).
For each $x \in \mathbb{Z}_n^2$, $1 \leq i \leq L_n - 1$, and $m \in \mathbb{N}$, set
\begin{equation} \label{eq:traversal-process}
T_i^{x, m} := \max \{k \geq 0 : R_k^{x, i} < D_m^{x, 0} \}.
\end{equation}
As in Section \ref{sec:upper-bound}, we will study
the traversal process
$(T_i^{x, m_n^{-}})_{i \leq L_n-1}$, $x \in \mathbb{Z}_n^2$. 
Set
\begin{equation} \label{eq:lower-bound-lower-bump}
f_n (s) := \min \{s^{1/2 - \delta}, ~(L_n - 1 - s)^{1/2 - \delta} \},
~s \in [0, L_n - 1],
\end{equation}
\begin{equation} \label{eq:lower-bound-upper-bump}
g_n (s) := \min \{s^{1/2 + \delta}, ~(L_n - 1 - s)^{1/2 + \delta} \},
~s \in [0, L_n - 1],
\end{equation}
where $\delta$ is the one in (\ref{eq:assumption-parameter}).
For each $x \in \mathbb{Z}_n^2$, set
\begin{equation} \label{eq:event-lower-bound}
A_n (x) := \left \{
b_n^{-} (i) \leq T_i^{x, m_n^{-}} \leq b_n^{+} (i),~
w_n \leq \forall i \leq L_n - 1 - w_n,~
H_x > D_{m_n^{-}}^{x, 0}
\right \},
\end{equation}
where
\begin{equation} \label{eq:lower-bound-lower-control-curves}
b_n^{-} (s) := \left \lceil \left \{ \left(1 - \frac{s}{L_n} \right) \sqrt{m_n^{-}} + f_n (s) \right \}^2 
\right \rceil,
~s \in [0, L_n - 1],
\end{equation}
\begin{equation} \label{eq:lower-bound-upper-control-curves}
b_n^{+} (s) := \left \lfloor \left \{ \left(1 - \frac{s}{L_n} \right) \sqrt{m_n^{-}} + g_n (s) \right \}^2
\right \rfloor,
~s \in [0, L_n - 1].
\end{equation}
We will apply the second moment method to
\begin{equation} \label{eq:1st-moment}
Z_n := \sum_{x \in \mathbb{Z}_n^2 \backslash B (o, r_0)} 1_{A_n (x)}
\end{equation}
and we need a lower bound of $E_o [Z_n]$ and
an upper bound of $E_o [Z_n^2]$.
To estimate $E_o [Z_n]$, we need a lower bound of $P_o (A_n (x))$, $x \in \mathbb{Z}_n^2$.
\begin{lem} \label{lem:lower-bound-first-moment}
There exist $c_1, c_2 \in (0, \infty)$ and $n_0 \in \mathbb{N}$
such that for all $n \geq n_0$ and $x \in \mathbb{Z}_n^2 \backslash B(o, r_0)$,
\begin{equation} \label{eq:lower-bound-first-moment}
P_o (A_n (x)) \geq c_1 n^{-2} w_n \log (\ell_n) e^{2 s_n} (\ell_n)^{- 2 c_{\star} w_n - c_2}.
\end{equation}
\end{lem}
{\it Proof.}
Fix $x \in \mathbb{Z}_n^2 \backslash B(o, r_0)$.
By the transfer lemma (Lemma \ref{lem:transfer-lemma}(i)),
we have
\begin{equation} \label{eq:lem-lower-bound-first-moment-3}
P_o [A_n (x)] \geq (1 + o (1)) P_{m_n^-}^{\text{GW}} \left[
\bigcap_{i=w_n}^{L_n - 1 - w_n}
\left\{b_n^{-} (i) \leq T_i \leq b_n^{+} (i) \right\}
\cap \left\{T_{L_n - 1} = 0 \right\}
\right].
\end{equation}
By the barrier estimate (Lemma \ref{lem:GW-process-barrier-estimate}(ii)),
the right of (\ref{eq:lem-lower-bound-first-moment-3}) is bounded from below
by
$c_1 \frac{w_n}{L_n} (1 - \frac{1}{L_n})^{m_n^-}$,
which yields the desired result. $\qed$ \\\\
From now, to obtain an upper bound of $E_o[Z_n^2]$,
we estimate $P_{o} [A_n (x) \cap A_n (y)]$, $x, y \in \mathbb{Z}_n^2 \backslash B (o, r_0)$.
When two balls centered at $x$ and $y$ overlap,
the traversal processes in those balls are heavily correlated,
while, when two balls are disjoint,
the traversal processes in the balls should be almost independent.
We will use the term ``decoupling" to describe
such procedure to break dependence of the traversal processes.
Based on this observation, we define the branching level by
\begin{equation} \label{eq:branching-level}
\ell (x, y) := \min \{k \geq 0 : B (x, r_k) \cap B (y, r_k) = \emptyset \}.
\end{equation}
Recall the definition of $d_n (\cdot)$ from (\ref{eq:d-function}).
We will consider six cases:
\begin{equation*}
\begin{aligned}
&(1) \ell (x, y) = 0,\,\,\, (2) 1 \leq \ell (x, y) \leq w_n + 2,\,\,\, 
(3) w_n + 3 \leq \ell (x, y) \leq d_n (\xi),\\
&(4) d_n (\xi) < \ell (x, y) 
\leq \lceil (1 - \varepsilon) L_n \rceil,\,\,\,
(5) \lceil (1 - \varepsilon) L_n \rceil < \ell (x, y) \leq L_n - w_n - 1,\\
&(6) \ell(x, y) \geq L_n - w_n.
\end{aligned}
\end{equation*}
In the case $(4)$, following the argument in \cite[Section 6.2]{BeKi},
we will decouple the traversal processes
by conditioning on excursions outside the balls
and applying the Harnack inequality to remove the dependence of the endpoints of the excursions.
In the case $(5)$, 
following the argument in \cite[Section 6.1]{BeKi},
we use a recursion argument based on the strong Markov property of SRW
to obtain a decoupling inequality.
In the case $(6)$, we do not need decoupling estimates
since the number of the pairs $(x, y)$ is small. 
In cases $(1)$--$(3)$,
to decouple the traversal processes,
we use the following lemma
which states that excursions in disjoint regions
are almost independent.
\begin{lem} \label{lem:decomposition-trajectories}
Let $n \in \mathbb{N}$ and $0 < r < R < R^{\prime} < \widetilde{R} < n/4$.
Let $x, y \in \mathbb{Z}_n^2$
with $B (x, R) \cap B (y, R) = \emptyset$ and $B (y, R) \subset B (x, R^{\prime})$.
We define the space $\mathcal{W}$
of nearest-neighbor paths of finite length on $\mathbb{Z}_n^2$ by
$$\mathcal{W} := \left \{\omega : ~
\begin{minipage}{160pt}
$\exists I_{\omega} \in \mathbb{N},~\forall i \leq I_{\omega},~\omega (i) \in \mathbb{Z}_n^2$, \\
$d (\omega (i), \omega (i - 1)) = 1,~~1 \leq \forall i \leq I_{\omega}$
\end{minipage}
\right \}.$$
For each $z \in \{x, y \}$, 
we define the space $\mathcal{W}_z$ of
paths from $\partial B (z, r)$ to $\partial B (z, R)$ by
$$\mathcal{W}_z := \left \{\omega \in \mathcal{W} : ~
\begin{minipage}{160pt}
$\omega (0) \in \partial B (z, r),~\omega (I_{\omega}) \in \partial B (z, R)$, \\
$\omega (i) \in B (z, R),~~\forall i < I_{\omega}$ 
\end{minipage}
\right \}.$$
Similarly, we define the space $\widetilde{\mathcal{W}}_x$ of
paths from $\partial B (z, \widetilde{R})$ to $\partial B (z, R^{\prime})$ by
$$\widetilde{\mathcal{W}}_x := \left \{\omega \in \mathcal{W} : ~
\begin{minipage}{170pt}
$\omega (0) \in \partial B (x, \widetilde{R}),~\omega (I_{\omega}) \in \partial B (x, R^{\prime})$, \\
$\omega (i) \in \mathbb{Z}_n^2 \backslash \overline{B (x, R^{\prime})},~~\forall i < I_{\omega}$ 
\end{minipage}
\right \},$$
where for $A \subset \mathbb{Z}_n^2$, we set $\overline{A} := A \cup \partial A$.\\\\
For any $k, \ell, m \in \mathbb{N}$, 
$v \in \overline{B (x, R^{\prime})} \backslash \left(B (x, R) \cup B (y, R) \right)$,
$\widetilde{F}_i \subset \widetilde{\mathcal{W}}_x$, $1 \leq i \leq k$, 
$F_i^x \subset \mathcal{W}_x$, $1 \leq i \leq \ell$,
$F_i^y \subset \mathcal{W}_y$, $1 \leq i \leq m$, 
set the event $\widetilde{\mathcal{E}_k}$
for $k$ excursions from $\partial B (x, \widetilde{R})$ to $\partial B (x, R^{\prime})$
by
\begin{equation*}
\widetilde{\mathcal{E}}_k := \bigcap_{i = 1}^k \left\{S_{\cdot \wedge H_{\partial B (x, R^{\prime})}}
\circ \theta_{D_i (x, \widetilde{R}, R^{\prime})} \in \widetilde{F}_i \right \},
\end{equation*}
set the event $\mathcal{E}_{\ell} (x)$
for $\ell$ excursions from $\partial B (x, r)$ to $\partial B (x, R)$
by 
\begin{equation*}
\mathcal{E}_{\ell} (x) := \bigcap_{i = 1}^{\ell} \left\{S_{\cdot \wedge H_{\partial B (x, R)}}
\circ \theta_{R_i (x, R, r)} \in F_i^x \right \},
\end{equation*}
set the event $\mathcal{E}_m (y)$
for $m$ excursions from $\partial B (y, r)$ to $\partial B (y, R)$
by 
\begin{equation*}
\mathcal{E}_m (y) := \bigcap_{i = 1}^m \left\{S_{\cdot \wedge H_{\partial B (y, R)}}
\circ \theta_{R_i (y, R, r)} \in F_i^y \right \}.
\end{equation*}
Then,
\begin{align} \label{eq:decomposition-trajectories}
&P_v
\left[\widetilde{\mathcal{E}}_k \cap \mathcal{E}_{\ell} (x) \cap \mathcal{E}_m (y) \right] \notag \\ 
&\leq \prod_{i=1}^k \max_{z \in \partial B (x, \widetilde{R})}
P_z \left[S_{\cdot \wedge H_{\partial B (x, R^{\prime})}} \in \widetilde{F}_i \right] \notag \\
&~~~\times \prod_{i=1}^{\ell} \max_{z \in \partial B (x, r)}
P_z \left[S_{\cdot \wedge H_{\partial B (x, R)}} \in F_i^x \right] 
\cdot \prod_{i=1}^m \max_{z \in \partial B (y, r)}
P_z \left[S_{\cdot \wedge H_{\partial B (y, R)}} \in F_i^y \right]. 
\end{align}
\end{lem}
\begin{rem}
Note that excursions in the events $\widetilde{\mathcal{E}}_k$, $\mathcal{E}_{\ell} (x)$,
and $\mathcal{E}_m (y)$ are disjoint.
Lemma \ref{lem:decomposition-trajectories} states that
the probability of the event for those disjoint excursions
is bounded by the product of probabilities of events
for each excursion.
\end{rem}
{\it Proof.}
The claim (\ref{eq:decomposition-trajectories}) follows from
the strong Markov property of the SRW and the transfinite induction in $(k, \ell, m)$.
We omit the details. $\qed$ \\\\
We first deal with the case 
$w_n + 3 \leq \ell (x, y) \leq d_n (\xi)$.
\begin{lem} \label{lem:2nd-moment-case3}
For any $\xi > 0$, as $n \to \infty$,
\begin{equation} \label{eq:2nd-moment-case3}
\sum_{k=w_n + 3}^{d_n (\xi)} 
\sum_{\begin{subarray}{c} x, y \in \mathbb{Z}_n^2 \backslash B (o, r_0) \\ \ell (x, y) = k 
\end{subarray}}
P_o [A_n (x) \cap A_n (y)]
= o(1) \left(E_o [Z_n] \right)^2.
\end{equation}
\end{lem}
Recall the notation
$D_m^{z, i}$, $R_m^{z, i}$
from (\ref{eq:abbreviation-stopping-time}).
Fix $w_n + 3 \leq k \leq d_n (\xi)$.
To apply Lemma \ref{lem:decomposition-trajectories},
we prepare some notation.
For each $z \in \mathbb{Z}_n^2$, $k \leq i \leq L_n - 1$, and $m \in \mathbb{N}$, set
\begin{equation} \label{eq:2nd-moment-case3-pf-traversal-process}
T_i^{k, z, m} := \max \left \{\ell \geq 1 : R_{\ell}^{z, i} < D_m^{z, k} \right \}.
\end{equation}
We define the interval $J_m^{k, z}$ by
\begin{equation*}
J_m^{k, z} :=
\begin{cases}
\left(D_m^{z, k},~D_{m + 1}^{z, k} \right) ~&\text{if}~m \neq b_n^+ (k), \\
\left(D_{b_n^+ (k)}^{z, k},~\infty \right) ~&\text{if}~m = b_n^+ (k).
\end{cases}
\end{equation*}
We define
$\widetilde{b}_m^{(k)} (i)$ by
\begin{equation} \label{eq:2nd-moment-3-case3-upper-functions-traversal-process}
\left \lceil \left \{\sqrt{m} \left(1 -  \frac{i + 1 - k}{L_n + 1 - k} \right)
+ \kappa \sqrt{\frac{(i - k + 1)(L_n - i)}{L_n + 1 - k}} \sqrt{\log \log n} \right \}^2 \right \rceil,
\end{equation}
where $\kappa > 2$ is a sufficiently large constant.
The curve $\widetilde{b}_m^{(k)} (\cdot)$ will play a role of an upper barrier
for the traversal process $T^{k, z, m}$
which controls approximation errors when we apply the transfer lemma (Lemma \ref{lem:transfer-lemma}).
Take any $x, y \in \mathbb{Z}_n^2 \backslash B (o, r_0)$
with $\ell (x, y) = k$.
Fix $z \in \{x, y \}$.
On the event $A_n (z)$,
there exists $b_n^- (k) \leq m \leq b_n^+ (k)$ such that
$H_z \in J_m^{k, z}$.
By this and the monotonicity of the traversal process
$T_{\cdot}^{k, z, \ell}$ in $\ell$,
we have
\begin{equation} \label{eq:lem-2nd-moment-case3-pf-00}
A_n (z) \subset C_{n, k} (z) \cap  \bigcup_{m = b_n^- (k)}^{b_n^+ (k)}
\bigcap_{i=k+1}^{L_n - d_n (2) - 1}
\left\{
b_n^{-} (i) \leq T_i^{k, z, m} \right\}
\cap \left\{H_z \in J_m^{k, z} \right\},
\end{equation}
where
\begin{equation} \label{eq:2nd-moment-case3-pf-initial-crossing}
C_{n, k} (z) := \left \{b_n^{-} (k-3) \leq T_{k - 3}^{z, m_n^-} \leq b_n^{+} (k-3) \right \}.
\end{equation}
Furthermore, we decompose the $m$-th event in (\ref{eq:lem-2nd-moment-case3-pf-00})
into the event where the traversal process $T_{\cdot}^{k, z, m}$
stays below the curve $\widetilde{b}_m^{(k)} (\cdot)$
and the event where $T_{\cdot}^{k, z, m}$ crosses the curve $\widetilde{b}_m^{(k)} (\cdot)$.
Then, we have
\begin{equation} \label{eq:lem-2nd-moment-case3-pf-0}
A_n (z) \subset \bigcup_{m = b_n^- (k)}^{b_n^+ (k)}
\left\{\left(C_{n, k} (z) \cap A_{n, m}^{(k)} (z) \right) \cup
\left(C_{n, k} (z) \cap B_{n, m}^{(k)} (z) \right) \right \},
\end{equation}
where for each $b_n^{-} (k) \leq m \leq b_n^{+} (k)$, set
\begin{equation} \label{eq:2nd-moment-case3-pf--modified-event}
A_{n, m}^{(k)} (z) := 
\bigcap_{i=k+1}^{L_n - d_n (2) - 1}
\left \{
b_n^{-} (i) \leq T_i^{k, z, m} \leq \widetilde{b}_m^{(k)} (i) \right\}
\cap \left\{H_z \in J_m^{k, z} \right\},
\end{equation}
\begin{equation} \label{eq:2nd-moment-case3-technical-event}
B_{n, m}^{(k)} (z) := 
\bigcup_{i=k+1}^{L_n - d_n (2) - 1}
\left \{
T_i^{k, z, m} > \widetilde{b}_m^{(k)} (i) \right\}
\cap \left\{H_z \in J_m^{k, z} \right\}.
\end{equation}
By (\ref{eq:lem-2nd-moment-case3-pf-0}),
$P_o [A_n (x) \cap A_n (y)]$ is bounded from above by
\begin{equation} \label{eq:lem-2nd-moment-case3-pf-1}
\sum_{M_x = b_n^{-} (k)}^{b_n^{+} (k)} \sum_{M_y = b_n^{-} (k)}^{b_n^{+} (k)}
(I_1 + I_2 + I_3 + I_4),
\end{equation}
where
\begin{align} \label{eq:lem-2nd-moment-case3-pf-2}
&I_1 := P_o [C_{n, k} (x) \cap A_{n, M_x}^{(k)} (x) \cap A_{n, M_y}^{(k)} (y)], \notag \\
&I_2 := P_o [C_{n, k} (x) \cap A_{n, M_x}^{(k)} (x) \cap B_{n, M_y}^{(k)} (y)], \notag \\
&I_3 := P_o [C_{n, k} (x) \cap B_{n, M_x}^{(k)} (x) \cap A_{n, M_y}^{(k)} (y)], \notag \\
&I_4 := P_o [C_{n, k} (x) \cap B_{n, M_x}^{(k)} (x) \cap B_{n, M_y}^{(k)} (y)].
\end{align}
Fix $M_x, M_y \in \{b_n^{-} (k), \cdots ,b_n^{+} (k) \}$.
Since the estimates of $I_1, I_3, I_4$
are similar to that of $I_2$,
we will mainly focus on $I_2$.
In the following lemma,
we decompose $I_2$ into the product of probabilities of events
corresponding to 
$C_{n, k} (x)$, $A_{n, M_x}^{(k)}(x)$, $B_{n, M_y}^{(k)} (y)$
and transfer the law of the traversal process
to that of the Galton-Watson process:
\begin{lem} \label{lem:lem-2nd-moment-case3-pf-tech-1}
$I_2$ in (\ref{eq:lem-2nd-moment-case3-pf-2}) is bounded from above by
$(1 + o (1)) P_1 P_2 P_3$, where
\begin{align} \label{eq:lem-2nd-moment-case3-pf-12}
&P_1 :=
P_{m_n^-}^{\text{GW}} \left[T_{k - 3} \leq b_n^+ (k - 3)
\right], \notag \\
&P_2 := P_{M_x}^{\text{GW}} \left[
\bigcap_{i=1}^{L_n - d_n (2) - k - 1}
\left\{T_i \geq b_n^- (k + i) \right\}
\cap \left\{T_{L_n - k - 1} = 0 \right\}
\right] \notag \\
&\,\,\,\,\,\,\,\,\,\,\,\,\,\,\,\,\,\,\,\,\,\,\,\,\,\,\,\,\,\,\,\,\,\,\,\,\,\,\,\,\,\,\,\,\,\,\,\,\,\,\,\,\,\,
\,\,\,\,\,\,\,\,\,\,\,\,\,\,\,\,\,\,\,\,\,\,\,\,\,\,\,\,\,\,\,\,\,\,\,\,\,\,\,\,\,\,\,\,\,\,\,\,\,\,\,\,\,\,
\,\,\,\,\,\,\,\,\,\,\,\,\,\,\,\,\,\,\,\,\,\,\,\
\times \left(\frac{1}{L_n - k} \right)^{1_{\{M_x \neq b_n^+ (k) \}}},
\notag \\
&P_3 :=
\sum_{q = k + 1}^{L_n - d_n (2) - 1} 
P_{M_y}^{\text{GW}} \left[T_{q - k} \geq \widetilde{b}_{M_y}^{(k)} (q),~
T_{L_n - k - 1} = 0 \right] 
\left(\frac{1}{L_n - k} \right)^{1_{\{M_y \neq b_n^+ (k) \}}}.
\end{align}
\end{lem}
\begin{rem}
$P_1, P_2, P_3$ in Lemma \ref{lem:lem-2nd-moment-case3-pf-tech-1} 
correspond to probabilities of events
$C_{n, k} (x)$, $A_{n, M_x}^{(k)}(x)$, $B_{n, M_y}^{(k)} (y)$ respectively.
\end{rem}
{\it Proof.}
For each $m \in \mathbb{N}$,
we define the number of traversals between
$\partial B (x, r_1)$ and $\partial B (x, r_0)$
by $m$ excursions 
from $\partial B (x, r_{k - 3})$ to $\partial B (x, r_{k - 2})$ by
$$\widetilde{T}^{k, x, m} := \max \{\ell \geq 1 : D_{\ell}^{x, 0} < R_m^{x, k-3} \}.$$
We first decompose $B_{n, M_y}^{(k)} (y)$ into events
with respect to the first time 
at which the traversal process crosses the curve $\widetilde{b}_{M_y}^{(k)} (\cdot)$:
\begin{equation*}
B_{n, M_y}^{(k)} (y)
= \bigcup_{k < q < L_n - d_n (2)} \bigcup_{m < M_y} B_{q, m},
\end{equation*}
where
\begin{equation*}
B_{q, m} :=
\left\{T_q^{k, y, m} < \widetilde{b}_{M_y}^{(k)} (q) \leq T_q^{k, y, m + 1},~
H_y \in J_{M_y}^{k, y} \right\}.
\end{equation*}
Then, we decompose $C_{n, k} (x)$
into events which are measurable to the traversal processes
before and after $R_1^{x, k - 3}$:
\begin{equation*}
C_{n, k} (x)
= \bigcup_{m^{\prime} = 0}^{m_n^- - 1} \bigcup_{p = b_n^- (k - 3) - 1}^{b_n^+ (k - 3) - 1}
C_{m^{\prime}}^1
\cap C_{m^{\prime}, p}^2,
\end{equation*}
where
$C_{m^{\prime}}^1
:= \left \{
D_{m^{\prime}}^{x, 0} < R_1^{x, k-3} < D_{m^{\prime} + 1}^{x, 0}
\right \}$ and
\begin{equation*}
C_{m^{\prime}, p}^2 :=
 \left \{
\widetilde{T}^{k, x, p} \circ \theta_{D_1^{x, k-3}} < m_n^- - m^{\prime}
\leq \widetilde{T}^{k, x, p + 1} \circ \theta_{D_1^{x, k-3}}
\right \}.
\end{equation*}
By these, we have
\begin{equation} \label{eq:lem-2nd-moment-case3-pf-3}
I_2 \leq \sum_{q, m, m^{\prime}, p} I_2 (q, m, m^{\prime}, p),
\end{equation}
where
\begin{equation} \label{eq:lem-2nd-moment-case3-pf-4}
I_2 (q, m, m^{\prime}, p) := E_o \left[1_{C_{{m}^{\prime}}^1}
P_{S_{H_{\partial B (x, r_{k-2})}}} \left[
C_{m^{\prime}, p}^2 \cap A_{n, M_x}^{(k)} (x) \cap B_{q, m} \right] \right],
\end{equation}
and the sum is taken over
$q = k + 1, \dotsc, L_n - d_n (2) - 1$,
$m = 0, \dotsc, M_y - 1$,
$m^{\prime} = 0, \dotsc, m_n^- - 1$,
$p = b_n^- (k-3) - 1, \dotsc, b_n^+ (k-3) - 1$.
Fix such $q, m, m^{\prime}, p$.
In order to apply Lemma \ref{lem:decomposition-trajectories},
we decompose the events
$C_{m^{\prime}, p}^2$, $A_{n, M_x}^{(k)} (x)$, $B_{q, m}$
into events for excursions.
First, we decompose $C_{m^{\prime}, p}^2$
into events for excursions from $\partial B (x, r_{k - 3})$ 
to $\partial B (x, r_{k - 2})$ as follows:
\begin{equation*}
\bigcup_{(j_{\ell}^c) \in U_c}
\bigcap_{\ell = 1}^p
\left\{\widetilde{T}^{k, x, 1} \circ \theta_{D_{\ell}^{x, k-3}} = j_{\ell}^c \right\}
\cap \left\{\widetilde{T}^{k, x, 1} \circ \theta_{D_{p + 1}^{x, k-3}}
\geq m_n^- - m^{\prime} - \sum_{\ell = 1}^p j_{\ell}^c \right\},
\end{equation*} 
where 
\begin{equation} \label{eq:condition-traversal-c}
U_c := \left\{(j_{\ell}^c)_{\ell = 1}^p \in \{0, 1, \dotsc \}^p : 
\sum_{\ell = 1}^p j_{\ell}^c < m_n^- - m^{\prime} \right\}.
\end{equation}
Then, we decompose $A_{n, M_x}^{(k)} (x)$
into events for excursions from $\partial B (x, r_{k + 1})$ to $\partial B (x, r_k)$ as follows:
\begin{equation*}
\begin{aligned}
\bigcup_{(j_{\ell}^{a, i}) \in U_a} \bigcap_{\ell = 1}^{M_x}
\bigcap_{i=k+1}^{L_n - d_n (2) - 1}
&\left \{
T_i^{k, x, 1} \circ \theta_{R_{\ell}^{x, k}} = j_{\ell}^{a, i} \right\} \\
&\cap \left\{H_x \circ \theta_{R_{\ell}^{x, k}} > D_{\ell}^{x, k} \right\}
\cap \left\{H_x \circ \theta_{R_{M_x + 1}^{x, k}} \in J_x \right\},
\end{aligned}
\end{equation*}
where 
\begin{equation} \label{eq:condition-traversal-a}
\begin{aligned}
U_a := &\bigl\{\{(j_{\ell}^{a, i})_{\ell = 1}^{M_x} \in \{0, 1, \dotsc\}^{M_x}
: k + 1 \leq i \leq L_n - d_n (2) - 1 \} :\\
&\,\,\,\,\,\,\,\,\,b_n^- (i) \leq \sum_{\ell = 1}^{M_x} j_{\ell}^{a, i} \leq \widetilde{b}_{M_x}^{(k)} (i),\,\,
k + 1 \leq \forall i \leq L_n - d_n (2) - 1 \bigr\}
\end{aligned}
\end{equation}
and $J_x$ is the interval $(0, D_{M_x + 1}^{x, k})$ if $M_x \neq b_n^+ (k)$
and is the interval $[0, \infty)$ if $M_x = b_n^+ (k)$.
Finally, we decompose $B_{q, m}$ into events
for excursions from $\partial B (y, r_{k + 1})$ to $\partial B (y, r_k)$ as follows:
\begin{align*}
\bigcup_{(j_{\ell}^b) \in U_b}
&\bigcap_{\ell = 1}^m
\left \{
T_q^{k, y, 1} \circ \theta_{R_{\ell}^{y, k}} = j_{\ell}^b,\,\,
H_y \circ \theta_{R_{\ell}^{y, k}} > D_{\ell}^{y, k}
\right \} \\
&\cap
\left \{
T_q^{k, y, 1} \circ \theta_{R_{m + 1}^{y, k}} \geq \widetilde{b}_{M_y}^{(k)} (q) - 
\displaystyle \sum_{\ell = 1}^m j_{\ell}^b,\,\,
H_y \circ \theta_{R_{m + 1}^{y, k}} > D_{m + 1}^{y, k}
\right \} \\
&\cap \bigcap_{\ell = m + 2}^{M_y} \left\{H_y \circ \theta_{R_{\ell}^{y, k}} > D_{\ell}^{y, k} \right\}
\cap \left\{H_y \circ \theta_{R_{M_y + 1}^{y, k}} \in J_y \right\},
\end{align*}
where
\begin{equation} \label{eq:condition-traversal-b}
U_b := \left\{(j_{\ell}^b)_{\ell = 1}^m \in \{0, 1, \dotsc \}^m :
\sum_{\ell = 1}^m j_{\ell}^b < \widetilde{b}_{M_y}^{(k)} (q) \right\}
\end{equation}
and $J_y$ is obtained from $J_x$ by replacing $x$ with $y$. 
By these decompositions and
Lemma \ref{lem:decomposition-trajectories},
the probability in (\ref{eq:lem-2nd-moment-case3-pf-4})
is bounded from above by
$\sum P_c P_a P_b$,
where the sum is taken over $(j_{\ell}^c) \in U_c$,
$(j_{\ell}^{a, i}) \in U_a$,
$(j_{\ell}^b) \in U_b$
and $P_c$, $P_a$, $P_b$ are defined as follows:
\begin{align} \label{eq:lem-2nd-moment-case3-pf-pc}
P_c := &\left\{\prod_{\ell = 1}^p
\max_{z \in \partial B (x, r_{k-3})} P_z [\widetilde{T}^{k, x, 1} = j_{\ell}^c] \right \} \notag \\
&\times \max_{z \in \partial B (x, r_{k - 3})}
P_z \left[\widetilde{T}^{k, x, 1} \geq m_n^- - m^{\prime} - \sum_{\ell = 1}^p j_{\ell}^c \right],
\end{align}
\begin{align} \label{eq:lem-2nd-moment-case3-pf-pa}
P_a := &\left\{\prod_{\ell = 1}^{M_x}
\max_{z \in \partial B (x, r_{k+1})}
P_z \left[
\bigcap_{i=k+1}^{L_n - d_n (2) - 1}
\left\{T_i^{k, x, 1} = j_{\ell}^{a, i} \right\} \cap \left\{H_x > H_{\partial B (x, r_k)} \right\}
\right] \right\} \notag \\
&\times \max_{z \in \partial B (x, r_{k+1})} P_z [H_x \in \widetilde{J}_x],
\end{align}
where $\widetilde{J}_x$ is the interval $(0, H_{\partial B (x, r_k)})$ if $M_x \neq b_n^+ (k)$
and is the interval $[0, \infty)$ if $M_x = b_n^+ (k)$,
\begin{align} \label{eq:lem-2nd-moment-case3-pf-pb}
P_b := &\prod_{\ell = 1}^m \max_{z \in \partial B (y, r_{k+1})} 
P_z \left[T_q^{k, y, 1} = j_{\ell}^b,~H_y > H_{\partial B (y, r_k)} \right] \notag \\
&\times \max_{z \in \partial B (y, r_{k+1})}
P_z \left[
T_q^{k, y, 1} \geq \widetilde{b}_{M_y}^{(k)} (q) - \sum_{\ell = 1}^m j_{\ell}^b,~
H_y > H_{\partial B (y, r_k)}
\right] \notag \\
&\times \prod_{\ell = m+2}^{M_y} \max_{z \in \partial B (y, r_{k+1})}
P_z [H_y > H_{\partial B (y, r_k)}]
\cdot \max_{z \in \partial B (y, r_{k+1})} P_z [H_y \in \widetilde{J}_y],
\end{align}
where $\widetilde{J}_y$ is obtained from $\widetilde{J}_x$
by replacing $x$ with $y$.
Recall the definitions of $P_1, P_2, P_3$ from (\ref{eq:lem-2nd-moment-case3-pf-12}).
First, applying the transfer lemma (Lemma \ref{lem:transfer-lemma})
to the first line of (\ref{eq:lem-2nd-moment-case3-pf-pa})
and Lemma \ref{lem:circle-to-circle-probab}
to the last lines of (\ref{eq:lem-2nd-moment-case3-pf-pa}),
we bound
$\sum_{(j_{\ell}^{a, i}) \in U_a} P_a$ by $(1 + o(1)) P_2$
(Recall the definition of $U_a$ from (\ref{eq:condition-traversal-a})).
Then, applying the transfer lemma (Lemma \ref{lem:transfer-lemma})
to the first and second lines of (\ref{eq:lem-2nd-moment-case3-pf-pb})
and Lemma \ref{lem:circle-to-circle-probab}
to the third line of (\ref{eq:lem-2nd-moment-case3-pf-pb}),
we bound
$\sum_{q, m, (j_{\ell}^b)} P_b$ by $(1 + o(1)) P_3$,
where the sum is taken over 
$q = k + 1, \dotsc, L_n - d_n (2) - 1$,
$m = 0, \dotsc, M_y - 1$,
$(j_{\ell}^b) \in U_b$
(Recall the definition of $U_b$ from (\ref{eq:condition-traversal-b})).
Finally, using Lemma \ref{lem:circle-to-circle-probab}, we have
\begin{equation} \label{eq:lem-2nd-moment-case3-pf-pc-last}
\sum_{(j_{\ell}^c), p, m^{\prime}} P_o [C_{m^{\prime}}^1] P_c 
\leq (1 + o(1)) P_1,
\end{equation} 
where the sum is taken over $(j_{\ell}^c) \in U_c$,
$p = b_n^- (k-3) - 1, \dotsc, b_n^+ (k-3) - 1$,
$m^{\prime} = 0, \dotsc, m_n^- - 1$
(Recall the definition of $U_c$ from (\ref{eq:condition-traversal-c})).
The proof of (\ref{eq:lem-2nd-moment-case3-pf-pc-last})
is not difficult but not straightforward,
so we give the proof in Section \ref{sec:pf-transfer-lem} for the sake of completeness.
Then, we have the desired result. $\qed$ \\\\
{\it Proof of Lemma \ref{lem:2nd-moment-case3}.}
We first estimate $I_2$ in (\ref{eq:lem-2nd-moment-case3-pf-2}).
(Since $I_3$ is basically $I_2$ with the roles of $M_x$ and $M_y$ changed,
we only treat $I_2$.)
By Lemma \ref{lem:lem-2nd-moment-case3-pf-tech-1},
we should estimate $P_1$, $P_2$, $P_3$ in (\ref{eq:lem-2nd-moment-case3-pf-12}).
First, we estimate $P_1$.
By \cite[(5.3)]{BRZ1}, $P_1$
is bounded from above by
\begin{equation} \label{eq:lem-2nd-moment-case3-pf-12-1}
\begin{aligned}
c_1 e^{- \frac{(\sqrt{m_n^-} - \sqrt{b_n^+ (k-3)})^2}{k-2}}
&\leq c_1 (\ell_n)^{-2k + 8} e^{c g_n (k-3) \sqrt{\log (\ell_n)}} \\
&\leq c_1 (\ell_n)^{-2k + 8} e^{c_2 (\log \log n)^{\frac{1}{2} + \delta - \alpha \delta}}.
\end{aligned}
\end{equation}
To simplify notation, we set
$$v := \sqrt{2 b_n^- (L_n - d_n (2) - 1)}.$$
Next, we estimate $P_2$ in (\ref{eq:lem-2nd-moment-case3-pf-12}).
By conditioning on $T_{L_n - k - d_n (2) - 1}$,
$P_2$
is bounded from above by
\begin{equation} 
\label{eq:lem-2nd-moment-case3-pf-12-2}
\begin{aligned}
\sum_{j = 0}^{\infty}
&P_{M_x}^{\text{GW}} \Biggl[
\bigcap_{i=1}^{L_n - d_n (2) - k - 2}
\left\{\sqrt{2 T_i} \geq \sqrt{2 b_n^- (k + i)} \right\} \\
&\,\,\,\,\,\,\,\,\,\,\,\,\,\,\,\,\,\,\,\,\,\,\,\,\,\,\,\,\,\,
\cap \left\{\sqrt{2 T_{L_n - k - d_n (2) - 1}} \in [v + j, v + j + 1) \right\}
\Biggr] \\
&\times
\sup_{u}
P_u^{\text{GW}} \left[T_{d_n (2)} = 0 \right],
\end{aligned}
\end{equation} 
where the supremum is taken over 
$u \in [\frac{(v + j)^2}{2}, \frac{(v + j + 1)^2}{2}] \cap \mathbb{N}$.
When $j \leq 100 L_n \sqrt{\log (\ell_n)}$,
to the $j$-th term of the sum in (\ref{eq:lem-2nd-moment-case3-pf-12-2}),
we can apply the barrier estimate
(Lemma \ref{lem:GW-process-barrier-estimate})
with
$L = L_n - k - d_n (2) - 1$,
$a = (1 - \frac{k}{L_n}) \sqrt{2 m_n^-}$,
$b = (1 - \frac{L_n - d_n (2) - 1}{L_n}) \sqrt{2 m_n^-}$,
$x = \sqrt{2M_x}$, $y = v + j$, $\eta = c \sqrt{\log (\ell_n)}$ ($c$ is a positive constant).
Then, the $j$-th term of the sum in (\ref{eq:lem-2nd-moment-case3-pf-12-2})
is bounded from above by
\begin{equation} \label{eq:lem-2nd-moment-case3-pf-12-3}
c_3 (\ell_n)^{c_4} \frac{j + 1}{\log n} e^{- \frac{(\sqrt{2 M_x} - v - j)^2}{2 (L_n - k - d_n (2) - 1)}}.
\end{equation}
The supremum in (\ref{eq:lem-2nd-moment-case3-pf-12-2})
is bounded from above by
\begin{equation} \label{eq:lem-2nd-moment-case3-pf-12-4}
\left(1 - \frac{1}{d_n (2) + 1} \right)^{\frac{(v+j)^2}{2}} \leq e^{- \frac{(v + j)^2}{2 (d_n (2) + 1)}}.
\end{equation}
Note that the product of exponential factors in (\ref{eq:lem-2nd-moment-case3-pf-12-3})
and (\ref{eq:lem-2nd-moment-case3-pf-12-4}) is equal to
$e^{- \frac{M_x}{L_n - k}}
e^{- \frac{L_n - k}{2 (L_n - k - d_n (2) - 1)(d_n (2) + 1)}
\left(v + j - \frac{\sqrt{2 M_x} (d_n (2) + 1)}{L_n - k} \right)^2}$.
By this, the product of
the sum over $j \leq 100 L_n \sqrt{\log (\ell_n)}$
and the supremum in (\ref{eq:lem-2nd-moment-case3-pf-12-2})
is bounded from above by
\begin{equation} \label{eq:lem-2nd-moment-case3-pf-12-5}
c_5  \frac{(\ell_n)^{c_6}}{\log n} e^{- \frac{M_x}{L_n - k}}.
\end{equation}
When $j > 100 L_n \sqrt{\log (\ell_n)}$,
bounding the supremum just by $1$
and using \cite[(5.3)]{BRZ1},
the product of the sum over $j > 100 L_n \sqrt{\log (\ell_n)}$
and the supremum in (\ref{eq:lem-2nd-moment-case3-pf-12-2})
is bounded from above by
\begin{equation} \label{eq:lem-2nd-moment-case3-pf-12-6}
P_{M_x}^{\text{GW}} \left[\sqrt{2 T_{L_n - k - d_n (2) - 1}} \geq 100 L_n \sqrt{\log (\ell_n)} \right] 
\leq n^{- 3}.
\end{equation}
Next, we estimate $P_3$ in (\ref{eq:lem-2nd-moment-case3-pf-12}). By conditioning on $T_{q - k}$,
the $q$-th term of the sum in $P_3$ is equal to
\begin{equation} \label{eq:lem-2nd-moment-case3-pf-12-7}
\sum_{m = \widetilde{b}_{M_y}^{(k)} (q)}^{\infty}
P_{M_y}^{\text{GW}} \left[T_{q-k} = m \right] P_m^{\text{GW}} \left[T_{L_n - q - 1} = 0 \right].
\end{equation}
For each $m \geq \widetilde{b}_{M_y}^{(k)} (q)$,
by \cite[(5.3)]{BRZ1}, we have
$P_{M_y}^{\text{GW}} [T_{q-k} = m] \leq e^{- \frac{(\sqrt{m} - \sqrt{M_y})^2}{q-k+1}}$.
We also have
$P_m^{\text{GW}} [T_{L_n - q - 1} = 0] = (1 - \frac{1}{L_n - q})^m \leq e^{- \frac{m}{L_n - q}}$.
By these, (\ref{eq:lem-2nd-moment-case3-pf-12-7}) is bounded from above by
a constant times
\begin{equation} \label{eq:lem-2nd-moment-case3-pf-12-8}
e^{- \frac{M_y}{L_n - k + 1}}
\sum_{m = \widetilde{b}_{M_y}^{(k)} (q)}^{\infty}
\exp \left \{- \frac{L_n - k + 1}{(q - k + 1)(L_n - q)}
\left(\sqrt{m} - \frac{L_n - q}{L_n - k + 1} \sqrt{M_y} \right)^2 \right \}.
\end{equation}
We decompose the sum in (\ref{eq:lem-2nd-moment-case3-pf-12-8})
into the sums over $\widetilde{b}_{M_y}^{(k)} (q) \leq m \leq 4 M_y$
and over $m > 4 M_y$.
Recall the definition of $\widetilde{b}_{M_y}^{(k)} (q)$ from 
(\ref{eq:2nd-moment-3-case3-upper-functions-traversal-process}).
In the sum over $\widetilde{b}_{M_y}^{(k)} (q) \leq m \leq 4 M_y$,
we use the bound 
$$\sqrt{m} - \frac{L_n - q}{L_n - k + 1} \sqrt{M_y}
\geq \kappa \sqrt{\frac{(q - k + 1)(L_n - q)}{L_n + 1 - k}} \sqrt{\log \log n}.$$
In the sum over $m > 4 M_y$,
we use the bound
$$\frac{L_n - q}{L_n - k + 1} \sqrt{M_y} \leq \frac{\sqrt{m}}{2}.$$
By these,
(\ref{eq:lem-2nd-moment-case3-pf-12-8})
is bounded from above by
\begin{equation} \label{eq:lem-2nd-moment-case3-pf-12-9}
c_7 e^{- \frac{M_y}{L_n - k + 1}}
\frac{(\log n)^{- \kappa^2 + 2}}{\log (\ell_n)}.
\end{equation}

By a direct calculation, we have
\begin{align} \label{eq:lem-2nd-moment-case3-pf-12-10}
&\sum_{M_y = b_n^- (k)}^{b_n^+ (k)} e^{- \frac{M_y}{L_n - k + 1}}
\left(\frac{1}{L_n - k} \right)^{1_{\{M_y \neq b_n^+ (k) \} }} 
\leq c_8 e^{- \frac{b_n^- (k)}{L_n - k + 1}} \notag \\
&\leq c_9 (\ell_n)^{- 2 (L_n - k) - 4 c_{\star} w_n + 2}
(\log n) e^{2 s_n - c_{10} f_n (k) \sqrt{\log (\ell_n)}}.
\end{align}

Note that the number of pairs
$x, y \in \mathbb{Z}_n^2 \backslash B (o, r_0)$
with $\ell (x, y) = k$
is bounded from above by
$n^2 (2 r_{k - 1})^2$.
By this, 
(\ref{eq:lem-2nd-moment-case3-pf-12-1}),
(\ref{eq:lem-2nd-moment-case3-pf-12-5}),
(\ref{eq:lem-2nd-moment-case3-pf-12-6}),
(\ref{eq:lem-2nd-moment-case3-pf-12-9}),
(\ref{eq:lem-2nd-moment-case3-pf-12-10}),
and Lemmas \ref{lem:lem-2nd-moment-case3-pf-tech-1} and \ref{lem:lower-bound-first-moment},
we have
\begin{align} \label{eq:lem-2nd-moment-case3-pf-12-11}
&\sum_{k = w_n + 3}^{d_n (\xi)}
\sum_{\begin{subarray}{c} x,y \in \mathbb{Z}_n^2 \backslash B(o, r_0), \\ 
\ell (x, y) = k \end{subarray}}
\sum_{b_n^- (k) \leq M_x, M_y \leq b_n^+ (k)} I_2 \notag \\
&\leq c_{11}
(\ell_n)^{- 2 c_{\star} w_n + c_{12}}
e^{c_2 (\log \log n)^{\frac{1}{2} + \delta - \alpha \delta}}
(\log n)^{- \kappa^2 + 4} \left\{E_o [Z_n] \right \}^2.
\end{align}
Since $\alpha + \beta > 1/2 + \delta - \alpha \delta$
by the assumption (\ref{eq:assumption-parameter}),
the right of (\ref{eq:lem-2nd-moment-case3-pf-12-11}) is
$o (1) E_o [Z_n]^2$.

Recall the definition of $I_3$
from (\ref{eq:lem-2nd-moment-case3-pf-2}).
As in Lemma \ref{lem:lem-2nd-moment-case3-pf-tech-1},
$I_3$ is bounded by $(1+o(1)) P_1 P_2^{\prime} P_3^{\prime}$,
where $P_1$ is the one in (\ref{eq:lem-2nd-moment-case3-pf-12})
and $P_2^{\prime}$, $P_3^{\prime}$
are $P_2$, $P_3$ in (\ref{eq:lem-2nd-moment-case3-pf-12})
with $x$ and $y$ swapped.
Thus,
$\sum_{k = w_n + 3}^{d_n (\xi)}
\sum_{\begin{subarray}{c} x,y \in \mathbb{Z}_n^2 \backslash B(o, r_0), \\ 
\ell (x, y) = k \end{subarray}}
\sum_{b_n^- (k) \leq M_x, M_y \leq b_n^+ (k)} I_3$
is $o (1) E_o [Z_n]^2$.

Recall the definition of $I_1$ from (\ref{eq:lem-2nd-moment-case3-pf-2}).
Since the event in $I_1$ is just the one in $I_2$ with $B_{n, M_y}^{(k)} (y)$ replaced by
$A_{n, M_y}^{(k)} (y)$, 
one can repeat
the proof of Lemma \ref{lem:lem-2nd-moment-case3-pf-tech-1}
(e.g. replace the equations around (\ref{eq:condition-traversal-b})
and (\ref{eq:lem-2nd-moment-case3-pf-pb})
by those around (\ref{eq:condition-traversal-a}) and (\ref{eq:lem-2nd-moment-case3-pf-pa})
(with $x$ swapped for $y$), repectively)
and bound
$I_1$ from above by
$(1+o(1)) P_1 P_2 P_2^{\prime}$,
where $P_1$ and $P_2$ are 
the ones in (\ref{eq:lem-2nd-moment-case3-pf-12})
and $P_2^{\prime}$ is $P_2$ with $x$ replaced by $y$. 
By this together with
(\ref{eq:lem-2nd-moment-case3-pf-12-1}), (\ref{eq:lem-2nd-moment-case3-pf-12-5}),
(\ref{eq:lem-2nd-moment-case3-pf-12-6}),
(\ref{eq:lem-2nd-moment-case3-pf-12-10}),
and Lemma \ref{lem:lower-bound-first-moment},
we have
\begin{align} \label{eq:lem-2nd-moment-case3-pf-30}
&\sum_{k = w_n + 3}^{d_n (\xi)}
\sum_{\begin{subarray}{c} x, y \in \mathbb{Z}_n^2 \backslash B (o, r_0); \\
\ell (x, y) = k \end{subarray}}
\sum_{M_x, M_y = b_n^{-} (k)}^{b_n^{+} (k)}
I_1 \notag \\
&\leq c_{13} (\ell_n)^{- 2 c_{\star} w_n + c_{14}}
e^{c_2 (\log \log n)^{\frac{1}{2} + \delta - \alpha \delta}}
\{E_o [Z_n] \}^2.
\end{align}
Since $\alpha + \beta > 1/2 + \delta - \alpha \delta$
by the assumption (\ref{eq:assumption-parameter}),
the right of (\ref{eq:lem-2nd-moment-case3-pf-30}) is
$o (1) \{E_o [Z_n] \}^2$.

Recall the definition of $I_4$ from (\ref{eq:lem-2nd-moment-case3-pf-2}).
Since the event in $I_4$ is just the one in $I_2$
with $A_{n, M_x}^{(k)} (x)$ replaced by $B_{n, M_x}^{(k)} (x)$,
one can repeat
the proof of Lemma \ref{lem:lem-2nd-moment-case3-pf-tech-1}
(e.g. replace equations around (\ref{eq:condition-traversal-a}) and (\ref{eq:lem-2nd-moment-case3-pf-pa})
by those around (\ref{eq:condition-traversal-b}) and (\ref{eq:lem-2nd-moment-case3-pf-pb})
(with $y$ swapped for $x$), respectively) and bound
$I_4$ from above by
$(1+o(1)) P_1 P_3^{\prime} P_3$,
where $P_1$ and $P_3$ are the ones in (\ref{eq:lem-2nd-moment-case3-pf-12})
and $P_3^{\prime}$ is $P_3$ with $y$ replaced by $x$.
By this together with 
(\ref{eq:lem-2nd-moment-case3-pf-12-1}), (\ref{eq:lem-2nd-moment-case3-pf-12-9}),
(\ref{eq:lem-2nd-moment-case3-pf-12-10}), and Lemma \ref{lem:lower-bound-first-moment},
we have
\begin{align} \label{eq:2nd-moment-case3-pf-step32}
&\sum_{k = w_n + 3}^{d_n (\xi)}
\sum_{\begin{subarray}{c} x, y \in \mathbb{Z}_n^2 \backslash B (o, r_0); \\
\ell (x, y) = k \end{subarray}}
\sum_{M_x, M_y = b_n^{-} (k)}^{b_n^{+} (k)}
I_4 \notag \\
&\leq c_{15} (\ell_n)^{- 2 c_{\star} w_n + c_{16}} 
e^{c_2 (\log \log n)^{\frac{1}{2} + \delta - \alpha \delta}}
(\log n)^{- 2 \kappa^2 + 8}
\{E_o [Z_n] \}^2.
\end{align}
Since $\alpha + \beta > 1/2 + \delta - \alpha \delta$
by the assumption (\ref{eq:assumption-parameter})
and $\kappa$ is sufficiently large,
the right  of (\ref{eq:2nd-moment-case3-pf-step32}) is $o (1) E_0 [Z_n]^2$.
Therefore, we have (\ref{eq:2nd-moment-case3}). $\qed$ \\\\

Next, we deal with the case $1 \leq \ell (x, y) \leq w_n + 2$.
\begin{lem} \label{lem:2nd-moment-case2}
\begin{equation} \label{eq:2nd-moment-case2}
\sum_{\begin{subarray}{c} x, y \in \mathbb{Z}_n^2 \backslash B (o, r_0), \\ 
1 \leq \ell (x, y) \leq w_n + 2 \end{subarray}}
P_o [A_n (x) \cap A_n (y)]
= o(1) \left(E_o [Z_n] \right)^2.
\end{equation}
\end{lem}
{\it Proof.}
The proof is almost the same as that of Lemma \ref{lem:2nd-moment-case3}
with $k = w_n + 3$
except that we do not consider the event $C_{n, k} (x)$,
so the proof is simpler.
We omit the detail.~~~$\qed$ \\\\

Recall the definition of $d_n (\cdot)$ from (\ref{eq:d-function}).
Next, we deal with the case 
$d_n (\xi)
< \ell (x, y) 
\leq \lceil (1 - \varepsilon) L_n \rceil$.
\begin{lem} \label{lem:2nd-moment-case4}
Fix $\varepsilon \in (0, 1/2)$.
There exists $\xi_0 > 0$
such that for any $\xi \geq \xi_0$,
\begin{equation} \label{eq:2nd-moment-case4}
\sum_{k = d_n (\xi) + 1}
^{\lceil (1 - \varepsilon) L_n \rceil} 
\sum_{\begin{subarray}{c} x, y \in \mathbb{Z}_n^2 \backslash B (o, r_0), \\ \ell (x, y) = k 
\end{subarray}}
P_o [A_n (x) \cap A_n (y)]
= o(1) \left(E_o [Z_n] \right)^2.
\end{equation}
\end{lem}
In the proof of Lemma \ref{lem:2nd-moment-case4},
we cannot apply the argument in the proof of Lemma \ref{lem:2nd-moment-case3}
because the factor $e^{c g_n (k-3) \sqrt{\log (\ell_n)}}$ in (\ref{eq:lem-2nd-moment-case3-pf-12-1})
can be too large.
Instead, we will follow the argument in \cite[Section 6.2]{BeKi}. 
We first prove the following decoupling lemma:
\begin{lem} \label{lem:2nd-moment-case4-decoupling}
There exist $c_1, c_2 \in (0, \infty)$ and $\xi_0 > 0$ such that
for any $\xi \geq \xi_0$,
$d_n (\xi) < k \leq \lfloor (1 - \varepsilon) L_n \rfloor$,
and $x, y \in \mathbb{Z}_n^2 \backslash B (o, r_0)$
with $d (x, y) = k$,
\begin{equation} \label{eq:2nd-moment-case4-decoupling}
P_o \left[A_n (x) \cap A_n (y) \right]
\leq c_1 (\ell_n)^{c_2} 
g_n \left(d_{n, k} (\eta) \right)^2 (\log n)^{-2} P_o \left[B_{b_n^- (k)} \right],
\end{equation}
where  $d_{n, k} (\eta) := k + d_n (\eta) + 1$
and
for each $m \in \mathbb{N}$, $B_m$ is defined by
\begin{equation} \label{eq:lem-2nd-moment-case4-pf-34}
B_m :=
\left\{
H_x > D_{m_n^-}^{x, 0},~
H_y > D_m^{y, k}
\right \}.
\end{equation}
\end{lem}
\begin{rem}
Roughly speaking,
the equation in (\ref{eq:2nd-moment-case4-decoupling})
says that the barrier events in $A_n (x)$ and $A_n (y)$
for $i \ge d_{n, k} (\eta)$
are almost independent 
conditioned on the event that $x$ and $y$ are avoided.
Thanks to the strong Markov property and the fact that
$B(x, r_k)$ and $B(y, r_k)$ are disjoint,
excursions
inside $B(x, r_k)$ and $B(y, r_k)$
are almost independent
conditioned on excursions outside the balls
except the dependence on endpoints of the excursions.
In order to get rid of the dependence on endpoints,
we will apply the Harnack inequality. 
\end{rem}
{\it Proof of Lemma \ref{lem:2nd-moment-case4-decoupling}.}
Take any $d_n (\xi) < k \leq \lceil (1 - \varepsilon) L_n \rceil$
and any $x, y \in \mathbb{Z}_n^2 \backslash B (o, r_0)$ with $\ell (x, y) = k$.
(Note that the constants appearing in the proof do not depend on
$x, y$.)
Fix $\eta > 2$ large enough.
To simplify notation, for $z \in \{x, y \}$ and $\ell \geq 1$, set
\begin{equation} \label{eq:lem-2nd-moment-case4-pf-notation-0}
\widetilde{R}_{\ell}^{z, k} := R_{\ell} (z, r_k, r_{d_{n, k} (\eta)}),~
\widetilde{D}_{\ell}^{z, k} := D_{\ell} (z, r_k, r_{d_{n, k} (\eta)}).
\end{equation}
Recall the definitions of
$R_{\ell}^{z, i}$, $D_{\ell}^{z, i}$
from (\ref{eq:abbreviation-stopping-time}).
For each $m \geq 1$ and $k \leq i \leq L_n - 1$,
we define $M_n^z$ and $\widehat{T}_i^{k, z, m}$ by
\begin{equation} \label{eq:lem-2nd-moment-case4-pf-notation-1}
M_n^z := \max \left \{j \geq 1 : 
\widetilde{D}_j^{z, k} < D_{m_n^{-}}^{z, 0} \right \},
\end{equation}
\begin{equation} \label{eq:lem-2nd-moment-case4-pf-notation-2}
\widehat{T}_i^{k, z, m}:=  \max \left \{\ell \geq 1 : R_{\ell}^{z, i}
< \widetilde{D}_m^{z, k} \right \}.
\end{equation}
Then,
$P_o [A_n (x) \cap A_n (y)]$ is bounded from above by
\begin{equation} \label{eq:lem-2nd-moment-case4-pf-1}
\sum_{m_x, m_y = 1}^{b_n^{+} (k + d_n (\eta))}
P_o \left[A_n (x) \cap A_n (y) \cap \bigcap_{z \in \{x, y \}} \{M_n^z = m_z \} \right].
\end{equation}
For each $z \in \{x, y \}$,
let
$\mathcal{G}_z$ be the $\sigma$-algebra generated by
$S_{\cdot \wedge \widetilde{R}_1^{z, k}}$ and
excursions
$$S_{\cdot \wedge H_{\partial B (z, d_{n, k} (\eta))}}
\circ \theta_{\widetilde{D}_{\ell}^{z, k}},~\ell \geq 1$$
from $\partial B (z, r_k)$ to $\partial B (z, r_{d_{n, k} (\eta)})$.
Fix $m_x, m_y \in \{1, \cdots, b_n^+ (k + d_n (\eta)) \}$.
The term in (\ref{eq:lem-2nd-moment-case4-pf-1}) 
corresponding to $m_x, m_y$ is bounded from above by
\begin{equation} \label{eq:lem-2nd-moment-case4-pf-2}
E_o \left[
1_B
P_o \left[\displaystyle \bigcup_{q = b_n^- (d_{n, k} (\eta))}^{b_n^+ (d_{n, k} (\eta))} 
\widehat{A}_{n, q} (y)
\Biggm | \mathcal{G}_y \right]
\right],
\end{equation}
where
\begin{equation*}
\begin{aligned}
B := &\bigcap_{i=d_{n, k} (\eta)}^{L_n - w_n - 1} \left \{
b_n^{-} (i) \leq \widehat{T}_i^{k, x, m_x} \leq b_n^{+} (i) \right\} \\
&\,\,\,\,\,\cap \left\{H_x > \widetilde{D}_{m_x}^{x, k},~
M_n^x = m_x, M_n^y = m_y,~
T_k^{y, m_n^{-}} \geq b_n^{-} (k) \right\},
\end{aligned}
\end{equation*}
\begin{equation*}
\begin{aligned}
\widehat{A}_{n, q} (y)
:= &\left\{
\widehat{T}_{d_{n, k} (\eta)}^{k, y, m_y} = q \right\}
\cap \bigcap_{i=d_{n, k} (\eta) + 1}^{L_n - w_n - 1}
\left\{b_n^{-} (i) \leq \widehat{T}_i^{k, y, m_y} \leq b_n^{+} (i) \right\} \\
&\cap \left\{H_y > \widetilde{D}_{m_y}^{y, k} \right\}.
\end{aligned}
\end{equation*}
By the strong Markov property, we have the following:
For any $M \in \mathbb{N}$,
\begin{align} \label{eq:lem-2nd-moment-case4-pf-3}
&P_o \left[
\bigcap_{\ell = 1}^M
\left\{S_{\cdot \wedge H_{\partial B (y, r_k)}} \circ
\theta_{\widetilde{R}_{\ell}^{y, k}} 
\in~\cdot~ \right \} \biggm| \mathcal{G}_y \right] \notag \\
&= \prod_{\ell = 1}^M
P_{S_{\widetilde{R}_{\ell}^{y, k}}}
\left[S_{\cdot \wedge H_{\partial B (y, r_k)}}^{\prime} \in~\cdot~
\biggm | S_{H_{\partial B (y, r_k)}}^{\prime}
= S_{\widetilde{D}_{\ell}^{y, k}} \right],
\end{align}
where 
$S^{\prime}$ denotes a SRW
independent from $S$.

Fix any $q \in \{b_n^- (d_{n, k} (\eta)), \dotsc, b_n^+(d_{n, k} (\eta)) \}$.
We can decompose the event $\widehat{A}_{n, q} (y)$
into events for excursions
from $\partial B (y, r_{d_{n, k} (\eta)})$
to $\partial B (y, r_k)$ as follows:
\begin{equation} 
\label{eq:lem-2nd-moment-case4-pf-decomposition-event}
\begin{aligned}
\widehat{A}_{n, q} (y) =
\bigcup \bigcap_{\ell = 1}^{m_y}
&\left\{
\widehat{T}_{d_{n, k} (\eta)}^{k, y, 1} \circ \theta_{\widetilde{R}_{\ell}^{y, k}} = q_{\ell} \right\} 
\cap \bigcap_{i=d_{n, k} (\eta) + 1}^{L_n - w_n - 1}
\left\{\widehat{T}_i^{k, y, 1} \circ \theta_{\widetilde{R}_{\ell}^{y, k}} = j_{\ell}^{y, i} \right\} \\
&\cap \left\{H_y \circ \theta_{\widetilde{R}_{\ell}^{y, k}} > \widetilde{D}_{\ell}^{y, k} \right\},
\end{aligned}
\end{equation}
where the union is taken over 
all sequences of nonnegative integers
$(q_{\ell})_{\ell = 1}^{m_y}$ and $(j_{\ell}^{y, i})_{\ell = 1}^{m_y}$,
$d_{n, k} (\eta) + 1 \leq i \leq L_n - w_n - 1$ with
\begin{equation} \label{eq:lem-2nd-moment-case4-pf-specific-number-of-traversals}
\sum_{\ell = 1}^{m_y} q_{\ell} = q,~~
b_n^- (i) \leq \sum_{\ell = 1}^{m_y} j_{\ell}^{y, i} \leq b_n^+ (i),~
d_{n, k} (\eta) + 1 \leq i \leq L_n - w_n - 1.
\end{equation}
By (\ref{eq:lem-2nd-moment-case4-pf-3}),
the probability of the big intersection in 
(\ref{eq:lem-2nd-moment-case4-pf-decomposition-event})
conditioned on $\mathcal{G}_y$
is equal to
\begin{equation} 
\label{eq:eq:lem-2nd-moment-case4-pf-3-1}
\begin{aligned}
&\prod_{\ell = 1}^{m_y} 
P_{S_{\widetilde{R}_{\ell}^{y, k}}} \Biggl[
\left\{\widehat{T}_{d_{n, k} (\eta)}^{k, y, 1} = q_{\ell} \right\}
\cap \displaystyle\bigcap_{i=d_{n, k} (\eta) + 1}^{L_n - w_n - 1}
\left\{\widehat{T}_i^{k, y, 1} = j_{\ell}^{y, i} \right\} \\
&\,\,\,\,\,\,\,\,\,\,\,\,\,\,\,\,\,\,\,\,\,\,\,\,\,\,\,\,\,\,\,\,\,\,\,\,
\,\,\,\,\,\,\,\,\,\,\,\,\,\,\,\,
\cap \left\{
H_y > H_{\partial B (y, r_k)},\,
S_{H_{\partial B (y, r_k)}}^{\prime} = S_{\widetilde{D}_{\ell}^{y, k}}
\right\}
\Biggr] \\
&\times \prod_{\ell = 1}^{m_y} P_{S_{\widetilde{R}_{\ell}^{y, k}}}
\left[S_{H_{\partial B (y, r_k)}}^{\prime} = S_{\widetilde{D}_{\ell}^{y, k}} \right]^{-1}.
\end{aligned}
\end{equation}
Fix any 
$1 \leq \ell \leq m_y$.
Recall the notation
$T_i^{d_{n, k} (\eta), y, q_{\ell}}$ from (\ref{eq:2nd-moment-case3-pf-traversal-process}).
By the strong Markov property at $D_{q_{\ell}}^{y, d_{n, k} (\eta)}$,
the $\ell$-th probability in the first line
of (\ref{eq:eq:lem-2nd-moment-case4-pf-3-1})
is equal to
\begin{equation} 
\label{eq:lem-2nd-moment-case4-pf-added01}
\begin{aligned}
E_{S_{\widetilde{R}_{\ell}^{y, k}}} &\Biggl[
1_{\left\{
R_{q_{\ell}}^{y, d_{n, k} (\eta)} < H_{\partial B (y, r_k)} \right\}}
1_{\bigcap_{i=d_{n, k} (\eta) + 1}^{L_n - w_n - 1}
\left\{T_i^{d_{n, k} (\eta), y, q_{\ell}}
= j_{\ell}^{y, i} \right\}} 1_{\left\{H_y > D_{q_{\ell}}^{y, d_{n, k} (\eta)} \right\}} \\
&\,\,\,\,\,\,\times P_{S_{D_{q_{\ell}}^{y, d_{n, k} (\eta)}}^{\prime}}
\left[
S_{H_{\partial B (y, r_k)}}^{\prime \prime} = S_{\widetilde{D}_{\ell}^{y, k}},\, 
H_{\partial B (y, r_k)} < H_{\partial B (y, d_{n, k} (\eta) + 1)}
\right]
\Biggr],
\end{aligned}
\end{equation}
where $S^{\prime}$, $S^{\prime \prime}$ denote SRWs
independent from $S$ and from each other.

We will estimate the probability in (\ref{eq:lem-2nd-moment-case4-pf-added01}).
To do so, we fix
$z \in \partial B (y, r_{d_{n, k} (\eta)})$
and $w \in \partial B (y, r_k)$.
By the strong Markov property, we have
\begin{align} \label{eq:lem-2nd-moment-case4-pf-3-3}
&P_z \left[S_{H_{\partial B (y, r_k)}} = w,~H_{\partial B (y, r_{d_{n, k} (\eta) + 1})} 
< H_{\partial B (y, r_k)} \right] \notag \\
&= E_z \left[
1_{\{H_{\partial B (y, r_{d_{n, k} (\eta) + 1})} < H_{\partial B (y, r_k)} \}}
P_{S_{D_1^{y, d_{n, k} (\eta)}}} \left[S_{H_{\partial B (y, r_k)}}^{\prime} = w \right] \right].
\end{align}
We use the Harnack inequality from \cite[Lemma 2.1]{DPRZ2}:
Uniformly in $\delta^{\prime} < 1/2$, 
$u, u^{\prime} \in B (0, \delta^{\prime} N)$,
$u^{\prime \prime} \in \partial B (0, N)$,
\begin{equation} \label{eq:uniform-estimate-harmonic-measure}
P_u \left[S_{H_{\partial B (0, N)}} = u^{\prime \prime} \right]
= (1 + O(\delta^{\prime}) + O(1/N))
P_{u^{\prime}} \left[S_{H_{\partial B (0, N)}} = u^{\prime \prime} \right].
\end{equation}
By (\ref{eq:uniform-estimate-harmonic-measure}),
the $P_{S_{D_1^{y, d_{n, k} (\eta)}}}$-probability in (\ref{eq:lem-2nd-moment-case4-pf-3-3})
is equal to
\begin{equation} \label{eq:lem-2nd-moment-case4-pf-3-4}
\left(1 + O \left(\frac{r_{d_{n, k} (\eta)}}{r_k} \right) \right)
P_z \left[S_{H_{\partial B (y, r_k)}} = w \right].
\end{equation}
By (\ref{eq:lem-2nd-moment-case4-pf-3-3})
and Lemma \ref{lem:circle-to-circle-probab},
the $P_{S_{D_{q_{\ell}}^{y, d_{n, k} (\eta)}}^{\prime}}$-probability 
in (\ref{eq:lem-2nd-moment-case4-pf-added01})
is bounded from above by
$(1 + \frac{(\log \log n)^{c_1}}{(\log n)^{\eta}})$ times
\begin{equation} \label{eq:lem-2nd-moment-case4-pf-3-5}
P_{S_{D_{q_{\ell}}^{y, d_{n, k} (\eta)}}^{\prime}}
 \left[H_{\partial B (y, r_k)} < H_{\partial B (y, r_{d_{n, k} (\eta) + 1})} \right]
P_{S_{D_{q_{\ell}}^{y, d_{n, k} (\eta)}}^{\prime}} 
\left[S_{H_{\partial B (y, r_k)}}^{\prime \prime} = S_{\widetilde{D}_{\ell}^{y, k}} \right].
\end{equation}
By this, (\ref{eq:lem-2nd-moment-case4-pf-added01})
is bounded from above by
$(1 + \frac{(\log \log n)^{c_1}}{(\log n)^{\eta}})$
times
\begin{equation}
\label{eq:lem-2nd-moment-case4-pf-added02}
\begin{aligned}
&P_{S_{\widetilde{R}_{\ell}^{y, k}}} \Biggl[
\left\{R_{q_{\ell}}^{y, d_{n, k} (\eta)} < H_{\partial B (y, r_k)} < R_{q_{\ell} + 1}^{y, d_{n, k} (\eta)} \right\} \\
&\,\,\,\,\,\,\,\,\,\,\,\,\,\,\,\,\,\,\,\,\,\,\,\,\,\,\,\,
\cap \bigcap_{i=d_{n, k} (\eta) + 1}^{L_n - w_n - 1}
\left\{T_i^{d_{n, k} (\eta), y, q_{\ell}}
= j_{\ell}^{y, i} \right\} 
\cap \left\{H_y > D_{q_{\ell}}^{y, d_{n, k} (\eta)} \right\}
\Biggr] \\
&\times 
\max_{z \in \partial B (y, d_{n, k} (\eta))}
P_z
\left[
S_{H_{\partial B (y, r_k)}}^{\prime} = S_{\widetilde{D}_{\ell}^{y, k}}
\right].
\end{aligned}
\end{equation}

By (\ref{eq:lem-2nd-moment-case4-pf-added02}), 
the Harnack inequality (\ref{eq:uniform-estimate-harmonic-measure}),
and the transfer lemma (Lemma \ref{lem:transfer-lemma}(iii))
with $L = L_n - k$, $m = q_{\ell}$, $R_i = r_{k + i}$, $\widetilde{k} = w_n$,
$k = d_n (\eta) + 1$,
the conditional probability in (\ref{eq:lem-2nd-moment-case4-pf-2})
is bounded from above by $(1 + o (1))$ times
\begin{equation} 
\label{eq:lem-2nd-moment-case4-pf-3-6}
\begin{aligned}
&\sum_q P_q^{\text{GW}} \Biggl[
\bigcap_{i=1}^{L_n - d_{n, k} (\eta) - w_n - 1}
\left\{b_n^- (d_{n, k} (\eta) + i) \leq T_i \leq b_n^{+} (d_{n, k} (\eta) + i) \right\} \\
&\,\,\,\,\,\,\,\,\,\,\,\,\,\,\,\,\,\,\,\,\,\,\,\,\,\,\,\,\,\,\,\,\,\,\,\,\,\,\,\,\,\,\,\,\,
\,\,\,\,\,\,\,\,\,\,\,\,\,\,\,\,\,\,\,\,\,\,\,\,\,\,\,\,\,\,\,\,\,\,\,\,\,\,\,\,\,\,\,\,\,
\,\,\,\,\,\,\,\,\,\,\,\,\,\,\,\,\,\,\,\,\,\,\,\,\,\,\,\,\,\,\,\,\,\,\,\,\,\,\,\,\,\,\,\,\,
\cap
\left\{T_{L_n - d_{n, k} (\eta) - 1} = 0 \right\}
\Biggr] \\
&\times \sum_{(q_{\ell})_{\ell = 1}^{m_y}}
\left(\frac{d_n (\eta) + 1}{d_n (\eta) + 2} \right)^q
\left(\frac{1}{d_n (\eta) + 2} \right)^{m_y},
\end{aligned}
\end{equation}
where the sums in the first and the second lines of (\ref{eq:lem-2nd-moment-case4-pf-3-6}) 
are taken over $q \in \{b_n^- (d_{n, k} (\eta)), \dotsc, b_n^+(d_{n, k} (\eta)) \}$ 
and $(q_{\ell})_{\ell = 1}^{m_y}$ satisfying the first condition in
(\ref{eq:lem-2nd-moment-case4-pf-specific-number-of-traversals}).

By a direct calculation, one can show that
the sum in the second line of (\ref{eq:lem-2nd-moment-case4-pf-3-6}) is equal to
\begin{equation} \label{eq:lem-2nd-moment-case4-pf-3-7}
P \left[\sum_{\ell = 1}^{m_y} G_{\ell} = q \right],
\end{equation}
where $G_{\ell}$, $\ell \geq 1$
are independent geometric random variables with
success probability $\frac{1}{d_n (\eta) + 2}$.
By conditioning on $T_{L_n - d_{n, k} (\eta) - w_n - 1}$,
the $q$-th term of the sum in the first line of (\ref{eq:lem-2nd-moment-case4-pf-3-6})
is bounded from above by
\begin{equation} 
\label{eq:lem-2nd-moment-case4-pf-3-8}
\begin{aligned}
\sum_{v = b_n^- (L_n - w_n - 1)}^{b_n^+ (L_n - w_n - 1)}
&P_q^{\text{GW}} \Biggl[
\bigcap_{i=1}^{L_n - d_{n, k} (\eta) - w_n - 2}
\left\{T_i \geq b_n^- (d_{n, k} (\eta) + i) \right\} \\
&\,\,\,\,\,\,\,\,\,\,\,\,\,\,\,\,\,\,\,\,\,\,\,\,\,\,\,\,\,\,\,\,
\,\,\,\,\,\,\,\,\,\,\,\,\,\,\,\,\,\,\,\,\,\,\,\,\,\,\,\,\,\,\,\,
\cap 
\left\{T_{L_n - d_{n, k} (\eta) - w_n - 1} = v \right\}
\Biggr] \\
&\times P_v^{\text{GW}} \left[T_{w_n} = 0 \right].
\end{aligned}
\end{equation}
Applying the barrier estimate (Lemma \ref{lem:GW-process-barrier-estimate})
with $x = \sqrt{2q}$, $y = \sqrt{2v}$,
\footnote{Note that these are notations used in Lemma \ref{lem:GW-process-barrier-estimate}.
Please do not confuse these $x, y$ with the points on $\mathbb{Z}_n^2$.}
$a = \frac{L_n - d_{n, k} (\eta)}{L_n} \sqrt{2 m_n^-}$,
and $b = \frac{w_n + 1}{L_n} \sqrt{2 m_n^-}$ 
(note that we can take $\eta = c \sqrt{\log (\ell_n)}$ for some positive constant $c$),
we can bound
the first probability of the $v$-th term in (\ref{eq:lem-2nd-moment-case4-pf-3-8})
from above by
\begin{equation} \label{eq:lem-2nd-moment-case4-pf-3-9}
(\ell_n)^{c_2} g_n (d_{n, k} (\eta)) \frac{1}{\log n}
e^{- \frac{(\sqrt{q} - \sqrt{v})^2}{L_n - d_{n, k} (\eta) - w_n - 1}}.
\end{equation} 
For each $v \in \mathbb{N}$, we have
\begin{equation} \label{eq:lem-2nd-moment-case4-pf-3-10}
P_v^{\text{GW}} [T_{w_n} = 0] = \left(1 - \frac{1}{w_n + 1} \right)^v \leq e^{- \frac{v}{w_n + 1}}.
\end{equation}
We can bound the product of the exponential factors in 
(\ref{eq:lem-2nd-moment-case4-pf-3-9}) and (\ref{eq:lem-2nd-moment-case4-pf-3-10})
by
$e^{- \frac{q}{L_n - d_{n, k} (\eta)}}$
and this is bounded from above by
$(1 + o(1)) (1 - \frac{1}{L_n - d_{n, k} (\eta)})^q$.
Thus, (\ref{eq:lem-2nd-moment-case4-pf-3-6})
is bounded from above by
\begin{equation} \label{eq:lem-2nd-moment-case4-pf-3-11}
c_3 (\ell_n)^{c_4} g_n(d_{n, k} (\eta)) \frac{1}{\log n}
E \left[\left(1 - \frac{1}{L_n - d_{n, k} (\eta)} \right)^{\sum_{\ell = 1}^{m_y} G_{\ell}} \right].
\end{equation}
Using the geometric distribution of $G_{\ell}$, one can compute the expectation in 
(\ref{eq:lem-2nd-moment-case4-pf-3-11}) 
and show that (\ref{eq:lem-2nd-moment-case4-pf-3-11}) 
is equal to
\begin{equation}
\label{eq:lem-2nd-moment-case4-pf-3-11-equal}
c_3 (\ell_n)^{c_4} g_n(d_{n, k} (\eta)) \frac{1}{\log n}
\left(1 - \frac{d_n (\eta) + 1}{L_n - k} \right)^{m_y}.
\end{equation}

Next task is to prove that 
the last factor $(1 - \frac{d_n (\eta) + 1}{L_n - k})^{m_y}$ 
in (\ref{eq:lem-2nd-moment-case4-pf-3-11-equal})
is bounded from above by
\begin{equation} \label{eq:lem-2nd-moment-case4-pf-24}
(1 + o(1))P_o \left[H_y > \widetilde{D}_{m_y}^{y, k}|\mathcal{G}_y \right].
\end{equation}
Since the event in (\ref{eq:lem-2nd-moment-case4-pf-24})
can be written by
$$\bigcap_{\ell = 1}^{m_y} \{H_y \circ \theta_{\widetilde{R}_{\ell}^{y, k}}
> H_{\partial B (y, r_k)} \circ \theta_{\widetilde{R}_{\ell}^{y, k}} \},$$
by (\ref{eq:lem-2nd-moment-case4-pf-3}),
the probability in (\ref{eq:lem-2nd-moment-case4-pf-24})
is equal to
\begin{equation} \label{eq:lem-2nd-moment-case4-pf-25}
\prod_{\ell = 1}^{m_y}
\frac{P_{S_{\widetilde{R}_{\ell}^{y, k}}}
\left[H_y > H_{\partial B (y, r_k)},
~S_{H_{\partial B (y, r_k)}}^{\prime}
= S_{\widetilde{D}_{\ell}^{y, k}}
\right]}
{P_{S_{\widetilde{R}_{\ell}^{y, k}}} \left[S_{H_{\partial B (y, r_k)}}^{\prime}
= S_{\widetilde{D}_{\ell}^{y, k}}
\right]}.
\end{equation}
We will estimate the probabilities in the numerator in (\ref{eq:lem-2nd-moment-case4-pf-25}).
Fix any $z \in \partial B (y, r_{d_{n, k} (\eta)})$
and $z^{\prime} \in \partial B (y, r_k)$.
By the strong Markov property, we have
\begin{align} \label{eq:lem-2nd-moment-case4-pf-26}
&P_z \left[H_y < H_{\partial B (y, r_k)},~
S_{H_{\partial B (y, r_k)}} = z^{\prime} \right] \notag \\
&= P_z \left[H_y < H_{\partial B (y, r_k)} \right]
E_y \left[
P_{S_{H_{\partial B (y, r_{d_{n, k} (\eta)})}}}
\left[S_{H_{\partial B (y, r_k)}} = z^{\prime} \right] \right].
\end{align}
By the Harnack inequality 
(\ref{eq:uniform-estimate-harmonic-measure}),
the right in (\ref{eq:lem-2nd-moment-case4-pf-26})
is equal to
\begin{equation} \label{eq:lem-2nd-moment-case4-pf-27}
(1 + O((\log n)^{- \eta}))
P_z \left[H_y < H_{\partial B (y, r_k)} \right]
P_z \left[S_{H_{\partial B (y, r_k)}} = z^{\prime} \right].
\end{equation} 
Note that we can relate probabilities in (\ref{eq:lem-2nd-moment-case4-pf-25})
and (\ref{eq:lem-2nd-moment-case4-pf-26}) by
\begin{equation}
\label{eq:lem-2nd-moment-case4-pf-27-continue}
\begin{aligned}
&P_z [H_y > H_{\partial B (y, r_k)},~ S_{H_{\partial B (y, r_k)}} = z^{\prime}] \\
&= P_z [S_{H_{\partial B (y, r_k)}} = z^{\prime}]
- P_z [H_y < H_{\partial B (y, r_k)}, S_{H_{\partial B (y, r_k)}} = z^{\prime}].
\end{aligned}
\end{equation}
By (\ref{eq:lem-2nd-moment-case4-pf-27}), (\ref{eq:lem-2nd-moment-case4-pf-27-continue}),
and Lemma \ref{lem:circle-to-circle-probab},
(\ref{eq:lem-2nd-moment-case4-pf-25})
is bounded from below by
\begin{equation} \label{eq:lem-2nd-moment-case4-pf-29}
(1 + o(1))
\prod_{\ell = 1}^{m_y}
P_{S_{R_{\ell}^{y, k}}} 
\left[H_y > H_{\partial B (y, r_k)} \right].
\end{equation}
By (\ref{eq:lem-2nd-moment-case4-pf-29}), Lemma \ref{lem:circle-to-circle-probab},
and the condition $k \leq \lceil (1 - \varepsilon) L_n \rceil$,
we have shown that the last factor in (\ref{eq:lem-2nd-moment-case4-pf-3-11-equal})
is bounded from above by (\ref{eq:lem-2nd-moment-case4-pf-24}).

By (\ref{eq:lem-2nd-moment-case4-pf-3-11})--(\ref{eq:lem-2nd-moment-case4-pf-24}),
(\ref{eq:lem-2nd-moment-case4-pf-2})
is bounded by
$c_5 (\ell_n)^{c_6} g_n (d_{n, k} (\eta)) \frac{1}{\log n}$
times
\begin{equation} 
\label{eq:lem-2nd-moment-case4-pf-31}
\begin{aligned}
&P_o \Biggl[
\bigcap_{i=d_{n, k} (\eta)}^{L_n - w_n - 1}
\left\{b_n^{-} (i) \leq \widehat{T}_i^{k, x, m_x} \leq b_n^+ (i) \right\}  \\
&\,\,\,\,\,\,\,\,\,\,\,\,\,\,\,\,
\cap \bigcap_{z \in \{x, y \}}
\left\{H_z > \widetilde{D}_{m_z}^{z, k},\,
M_n^z = m_z \right\}
\cap \left\{T_k^{y, m_n^-} \geq b_n^- (k) \right\}
\Biggr].
\end{aligned}
\end{equation}
Conditioning on $\mathcal{G}_x$
and repeating the argument above,
(\ref{eq:lem-2nd-moment-case4-pf-31})
is bounded from above by
$c_7 (\ell_n)^{c_8} g_n(d_{n, k} (\eta))^2 \frac{1}{(\log n)^2}$
times
\begin{equation} \label{eq:lem-2nd-moment-case4-pf-32}
P_o \left[
M_n^x = m_x, ~M_n^y = m_y,
T_k^{y, m_n^-} \geq b_n^- (k),\,
H_x > \widetilde{D}_{m_x}^{x, k},~
H_y > \widetilde{D}_{m_y}^{y, k}
\right].
\end{equation}
Summing over $m_x, m_y$,
we obtain the desired result (\ref{eq:2nd-moment-case4-decoupling}).
$\qed$ \\\\
{\it Proof of Lemma \ref{lem:2nd-moment-case4}.}
Recall the definition of the event $B_m$ from (\ref{eq:lem-2nd-moment-case4-pf-34}).
To estimate $P_o [B_{b_n^- (k)}]$,
we will compare $P_o [B_m]$ with $P_o [B_{m - 1}]$. 
Fix $m \in \mathbb{N}$.
Since $R_m^{y, k}$
is contained in one of the intervals
$(D_{\ell}^{x, 0}, D_{\ell + 1}^{x, 0})$, $0 \leq \ell \leq m_n^- - 1$
or in $(D_{m_n^-}^{x, 0}, \infty)$,
$P_o [B_m]$ is equal to
\begin{align} \label{eq:lem-2nd-moment-case4-pf-35}
&\sum_{\ell = 0}^{m_n^- - 1}
P_o \left[
D_{\ell}^{x, 0} < R_m^{y, k} < D_{\ell + 1}^{x, 0},~
H_x > D_{m_n^-}^{x, 0},~
H_y > D_m^{y, k}
\right] \notag \\
& + P_o \left[
D_{m_n^-}^{x, 0} < R_m^{y, k},~
H_x > D_{m_n^-}^{x, 0},~
H_y > D_m^{y, k}
\right].
\end{align}
Fix $\ell \in \{0, \cdots, m_n^- - 1 \}$.
By the strong Markov property at $R_m^{y, k}$,
the $\ell$-th term of the sum in (\ref{eq:lem-2nd-moment-case4-pf-35})
is equal to
\begin{equation} \label{eq:lem-2nd-moment-case4-pf-36}
E_o \left[
\begin{minipage}{180pt}
$1_{\left\{
D_{\ell}^{x, 0} < R_m^{y, k} < D_{\ell + 1}^{x, 0},~
H_x > R_m^{y, k},~H_y > D_{m - 1}^{y, k}
\right \}}$ \\
$\times P_{S_{R_m^{y, k}}} 
\left[H_x > D_{m_n^- - \ell}^{x, 0},~H_y > H_{\partial B (y, r_k)} \right]$
\end{minipage}
\right].
\end{equation}
We will estimate the probability in the expectation in (\ref{eq:lem-2nd-moment-case4-pf-36}).

Fix any $z \in \partial B (y, r_{k + 1})$.
Note that
we have
$B (y, r_k) \subset B (x, r_{k - 2})$.
Thus, 
by the strong Markov property at 
$H_{\partial B (x, r_{k - 2})}$
and at $H_{\partial B (x, r_0)}$,
we have
\begin{align} \label{eq:lem-2nd-moment-case4-pf-37}
&P_z \left[
H_x > D_{m_n^- - \ell}^{x, 0},~H_y > H_{\partial B (y, r_k)} \right]
\notag \\
&= \sum_{u \in \partial B (x, r_{k - 2})} \sum_{v \in \partial B (x, r_0)}
P_z \left[
\begin{minipage}{160pt}
$H_y > H_{\partial B (y, r_k)}$,~$H_x > H_{\partial B (x, r_{k - 2})}$, \\
$S_{H_{\partial B (x, r_{k - 2})}} = u$
\end{minipage}
\right] \notag \\
&\times P_u \left[H_x > H_{\partial B (x, r_0)},~S_{H_{\partial B (x, r_0)}} = v \right]
P_v \left[H_x > D_{m_n^- - \ell - 1}^{x, 0} \right].
\end{align}
We will estimate the second factor on the right-hand side of (\ref{eq:lem-2nd-moment-case4-pf-37}).

Fix $u \in \partial B (x, r_{k - 2})$ and $v \in \partial B (x, r_0)$.
By the strong Markov property at
$H_x$ and $H_{\partial B (x, r_{k - 2})} \circ \theta_{H_x}$,
$P_u [H_x < H_{\partial B (x, r_0)},~S_{H_{\partial B (x, r_0)}} = v]$
is equal to
\begin{equation} \label{eq:lem-2nd-moment-case4-pf-38}
P_u [H_x < H_{\partial B (x, r_0)}]
E_x \left[
P_{S_{H_{\partial B (x, r_{k - 2})}}}
\left[S_{H_{\partial B (x, r_0)}}^{\prime} = v \right] \right].
\end{equation}
By the Harnack inequality
(\ref{eq:uniform-estimate-harmonic-measure}),
(\ref{eq:lem-2nd-moment-case4-pf-38})
is equal to
\begin{equation} \label{eq:lem-2nd-moment-case4-pf-39}
\left\{1 + O \left((\ell_n)^{- k + 2} \right) \right\}
P_u \left[H_x < H_{\partial B (x, r_0)} \right]
P_u \left[S_{H_{\partial B (x, r_0)}} = v \right].
\end{equation}
By (\ref{eq:lem-2nd-moment-case4-pf-39}),
Lemma \ref{lem:circle-to-circle-probab},
and the condition $k \leq \lceil (1 - \varepsilon) L_n \rceil$,
the first probability in the second line of
(\ref{eq:lem-2nd-moment-case4-pf-37})
is equal to
\begin{equation} \label{eq:lem-2nd-moment-case4-pf-40}
\left\{1 - O \left((\ell_n)^{- k + 2} \right) \right \}
P_u \left[H_x > H_{\partial B (x, r_0)} \right]
P_u \left[S_{H_{\partial B (x, r_0)}} = v \right].
\end{equation}
By the Harnack inequality
(\ref{eq:uniform-estimate-harmonic-measure}),
the second probability in (\ref{eq:lem-2nd-moment-case4-pf-40})
is equal to
\begin{equation} \label{eq:lem-2nd-moment-case4-pf-41}
\left\{1 + O \left((\ell_n)^{- k + 2} \right) \right \}
P_{u_{\star}} \left[S_{H_{\partial B (x, r_0)}} = v \right],
\end{equation}
where $u_{\star}$ is a fixed vertex on $\partial B (x, r_{k - 2})$.
By Lemma \ref{lem:circle-to-circle-probab},
the first probability in (\ref{eq:lem-2nd-moment-case4-pf-40})
is equal to
\begin{equation} \label{eq:lem-2nd-moment-case4-pf-42}
\left(1 + O \left((\log n)^{- 2} \right) \right) \frac{L_n - k + 2}{L_n}.
\end{equation}
By (\ref{eq:lem-2nd-moment-case4-pf-40})-(\ref{eq:lem-2nd-moment-case4-pf-42}),
for $\xi$ large enough,
the right of (\ref{eq:lem-2nd-moment-case4-pf-37})
is equal to
\begin{align} \label{eq:lem-2nd-moment-case4-pf-43}
&\left\{1 + O \left((\log n)^{- 2} \right) \right \}
\frac{L_n - k + 2}{L_n}
P_z \left[H_y > H_{\partial B (y, r_k)},~H_x > H_{\partial B (x, r_{k - 2})} \right] \notag \\
&\times P_{u_{\star}} \left[H_x \circ \theta_{H_{\partial B (x, r_0)}} > 
D_{m_n^- - \ell}^{x, 0} \right].
\end{align}
We will estimate the first probability of (\ref{eq:lem-2nd-moment-case4-pf-43}).

By the strong Markov property 
at $H_y$ and at $H_{\partial B (y, r_k)} \circ \theta_{H_y}$,
we have
\begin{align} \label{eq:lem-2nd-moment-case4-pf-44}
&P_z \left[
H_y < H_{\partial B (y, r_k)},~
H_x > H_{\partial B (x, r_{k - 2})} \right] \notag \\
&= P_z \left[H_y < H_{\partial B (y, r_k)} \right]
E_y \left[
P_{S_{H_{\partial B (y, r_k)}}} \left[H_x > H_{\partial B (x, r_{k - 2})} \right] \right] \notag \\
&\geq \left(1 - \frac{c_1}{L_n - k} \right)
\frac{1}{L_n - k},
\end{align}
where we have used Lemma \ref{lem:circle-to-circle-probab}
and the fact that $d (u, x) \geq r_k$ for any $u \in \partial B (y, r_k)$
in the last inequality.
Thus, 
the first probability in (\ref{eq:lem-2nd-moment-case4-pf-43}) 
is equal to
\begin{align} \label{eq:lem-2nd-moment-case4-pf-45}
&P_z \left[H_x > H_{\partial B (x, r_{k -2})} \right]
- P_z \left[
H_y < H_{\partial B (y, r_k)},~
H_x > H_{\partial B (x, r_{k - 2})} \right] \notag \\
&\leq P_z \left[H_x > H_{\partial B (x, r_{k - 2})} \right]
\left\{1 - \left(1 - \frac{c_1}{L_n - k} \right) \frac{1}{L_n - k} \right \}.
\end{align}
By (\ref{eq:lem-2nd-moment-case4-pf-36}),
(\ref{eq:lem-2nd-moment-case4-pf-43}),
(\ref{eq:lem-2nd-moment-case4-pf-45}),
and
applying the strong Markov property at $R_m^{y, k}$
and Lemma \ref{lem:circle-to-circle-probab}
to the last term of (\ref{eq:lem-2nd-moment-case4-pf-35}),
(\ref{eq:lem-2nd-moment-case4-pf-35})
is bounded from above by
\begin{align} \label{eq:lem-2nd-moment-case4-pf-46}
&\left\{1 + c_2 (\log n)^{- 2} \right \}
\left\{1 - \left(1 - \frac{c_1}{L_n - k} \right) \frac{1}{L_n - k} \right \}
 \notag \\
&\times \left\{\sum_{\ell = 0}^{m_n^- - 1}
\frac{L_n - k + 2}{L_n} P_{u_{\star}}^{(\ell)} E_{\ell} + Q \right \},
\end{align}
where
$P_{u_{\star}}^{(\ell)}$ is the probability in the second line of (\ref{eq:lem-2nd-moment-case4-pf-43})
and
\begin{equation*}
E_{\ell} := E_o \left[
1_{\left\{
D_{\ell}^{x, 0} < R_m^{y, k} < D_{\ell + 1}^{x, 0},~
H_x > R_m^{y, k},~H_y > D_{m - 1}^{y, k}
\right \}}
P_{S_{R_m^{y, k}}} \left[H_x > H_{\partial B (x, r_{k - 2})} \right]
\right],
\end{equation*}
\begin{equation*}
Q := P_o \left[
D_{m_n^-}^{x, 0} < R_m^{y, k},~
H_x > D_{m_n^-}^{x, 0},~
H_y > D_{m - 1}^{y, k}
\right].
\end{equation*}
A similar argument implies that
$P_o [B_{m - 1}]$
is bounded from below by
the factor in the second line of (\ref{eq:lem-2nd-moment-case4-pf-46})
multiplied by
$1 - c_3 (\log n)^{- 2}$.
Thus, we have
\begin{equation} \label{eq:lem-2nd-moment-case4-pf-48}
P_o [B_m]
\leq \left\{1 + c_4 (\log n)^{- 2} \right \}
\left\{1 - \left(1 - \frac{c_4}{L_n - k} \right) \frac{1}{L_n - k} \right\}
P_o [B_{m - 1}].
\end{equation}
In particular, 
$P_o \left[B_{b_n^- (k)} \right]$
is bounded from above by
\begin{equation} \label{eq:lem-2nd-moment-case4-pf-49}
\left\{1 + \frac{c_4}{(\log n)^2} \right \}^{b_n^- (k)}
\left\{1 - \left(1 - \frac{c_4}{L_n - k} \right) \frac{1}{L_n - k} \right\}^{b_n^- (k)}
P_o \left[H_x > D_{m_n^-}^{x, 0} \right].
\end{equation}
The first factor of (\ref{eq:lem-2nd-moment-case4-pf-49})
is $1 + o(1)$.
Since $k \leq \lceil (1 - \varepsilon) L_n \rceil$,
the second factor of (\ref{eq:lem-2nd-moment-case4-pf-49})
is bounded from above by
$(\ell_n)^{c_5} e^{- \frac{b_n^- (k)}{L_n - k}}$
and this is bounded from above by
\begin{align} \label{eq:lem-2nd-moment-case4-pf-50}
&(\ell_n)^{- 2 (L_n - k) + c_6} (\log n) e^{2 s_n - 4 c_{\star} w_n \log (\ell_n)}
e^{- c_7 f_n (k) \sqrt{\log (\ell_n)}} \notag \\
&\times \exp \left\{- \left(\log \log n + 2 s_n - 4 c_{\star} w_n \log (\ell_n) \right)
\frac{\log (\ell_n)}{\log n} k \right \}.
\end{align}
Since
$\alpha + \beta < \gamma < 1$,
the last factor of (\ref{eq:lem-2nd-moment-case4-pf-50})
is bounded from above by $1$.
By the strong Markov property and Lemma \ref{lem:circle-to-circle-probab},
the probability in (\ref{eq:lem-2nd-moment-case4-pf-49})
is bounded from above by
\begin{equation} \label{eq:lem-2nd-moment-case4-pf-51}
\left(\frac{\log (r_1) + \frac{c_8}{\log (r_1)}}{\log (r_0)}
\right)^{m_n^-} 
\leq (1 + o(1)) e^{- \frac{m_n^-}{L_n}} 
\leq \frac{1 + o(1)}{n^2} (\log n) e^{2 s_n} (\ell_n)^{- 2 c_{\star} w_n}.
\end{equation}
By Lemma \ref{lem:2nd-moment-case4-decoupling},
(\ref{eq:lem-2nd-moment-case4-pf-49})-(\ref{eq:lem-2nd-moment-case4-pf-51}),
$P_o \left[A_n (x) \cap A_n (y) \right]$
is bounded from above by
\begin{equation} \label{eq:lem-2nd-moment-case4-pf-52}
c_9 (\ell_n)^{- 2 (L_n - k) + c_{10}} n^{- 2}
e^{4 s_n - 6 c_{\star} w_n \log (\ell_n)} 
g_n(d_{n, k} (\eta))^2
e^{- c_7 f_n (k) \sqrt{\log (\ell_n)}}.
\end{equation}
In particular, the sum of 
$P_o [A_n (x) \cap A_n (y)]$
over $x, y \in \mathbb{Z}_n^2 \backslash B (o, r_0)$
with $\ell (x, y) = k$
is bounded from above by
(\ref{eq:lem-2nd-moment-case4-pf-52})
multiplied by
$n^2 (2 r_{k - 1})^2$. 
By this and Lemma \ref{lem:lower-bound-first-moment},
the left of (\ref{eq:2nd-moment-case4})
is bounded from above by
\begin{equation} \label{eq:lem-2nd-moment-case4-pf-53}
c_{11} (\ell_n)^{c_{12}} e^{- 2 c_{\star} w_n \log (\ell_n)} 
\left \{E_o [Z_n] \right\}^2.
\end{equation}
Thus, we have the desired result.~~~$\qed$ 
\begin{rem}
The reader may wonder why we did not apply the argument in the proof of
Lemma \ref{lem:2nd-moment-case4} to the proof of Lemma \ref{lem:2nd-moment-case3}.
This is because, in the case $k < d_n (\xi)$,
the first factor in (\ref{eq:lem-2nd-moment-case4-pf-49}) is replaced by
$(1+c(\ell_n)^{-k+2})^{b_n^- (k)}$
and this is not negligible. 
In the Brownian motion case,
thanks to the rotationally invariance,
the error factor does not appear and thus
one can apply the same method even in the case $k < d_n (\xi)$.
See \cite[Section 6.2]{BeKi} for the details.
\end{rem}

Next, we deal with the case 
$\lceil (1 - \varepsilon) L_n \rceil
< \ell (x, y) 
\leq L_n - w_n - 1$.
\begin{lem} \label{lem:2nd-moment-case5}
Let $\varepsilon$ be the constant in Lemma \ref{lem:2nd-moment-case4}.
Then,
\begin{equation} \label{eq:2nd-moment-case5}
\sum_{k = \lceil (1 - \varepsilon) L_n \rceil + 1}
^{L_n - w_n - 1} 
\sum_{\begin{subarray}{c} x, y \in \mathbb{Z}_n^2 \backslash B (o, r_0) \\ \ell (x, y) = k 
\end{subarray}}
P_o [A_n (x) \cap A_n (y)]
= o(1) \left(E_o [Z_n] \right)^2.
\end{equation}
\end{lem}
{\it Proof.}
We will follow the argument in \cite[Section 6.1]{BeKi}.
Take any $\lceil (1 - \varepsilon) L_n \rceil < k \leq L_n - w_n - 1$,
any $x, y \in \mathbb{Z}_n^2 \backslash B(o, r_0)$
with $\ell (x, y) = k$,
and sufficiently large constant $\theta > 0$.
Recall the notation 
$R_{\ell}^{z, i}, D_{\ell}^{z, i}, T_i^{z, m}$,
$d_n (\cdot)$ 
from (\ref{eq:abbreviation-stopping-time}), (\ref{eq:traversal-process}), (\ref{eq:d-function}).
We have
\begin{equation} \label{eq:2nd-moment-case5-pf-1}
P_o [A_n (x) \cap A_n (y)] \leq P_o \left[C_{b_n^- (k)} \right],
\end{equation}
where for each $m \in \mathbb{N}$, we set
\begin{equation} \label{eq:2nd-moment-case5-pf-2}
C_m := 
\bigcap_{i=w_n}^{k - d_n (\theta) - 1}
\left \{
b_n^{-} (i) \leq T_i^{x, m_n^-} \leq b_n^{+} (i) \right\}
\cap
\left\{H_x > D_{m_n^{-}}^{x, 0},~
H_y > D_m^{y, k} \right\}.
\end{equation}
To estimate $P_o [C_{b_n^- (k)}]$,
we will compare $P_o [C_m]$ with $P_o [C_{m - 1}]$.

Fix any $m \in \mathbb{N}$.
Since $R_m^{y, k}$ is contained in one of the intervals
$(D_{\ell}^{x, 0}, D_{\ell + 1}^{x, 0})$, $0 \leq \ell \leq m_n^- - 1$ or
in $(D_{m_n^-}^{x, 0}, \infty)$, we have
$P_o [C_m]$ is equal to
\begin{equation} \label{eq:2nd-moment-case5-pf-3}
\sum_{\ell = 0}^{m_n^- - 1}
P_o \left[
\left\{D_{\ell}^{x, 0} < R_m^{y, k} < D_{\ell + 1}^{x, 0} \right \}
\cap C_m \right]
+ P_o \left[\left\{D_{m_n^-}^{x, 0} < R_m^{y, k} \right \} \cap C_m \right].
\end{equation}
We define the number of traversals before time $R_m^{y, k}$ by
\begin{equation} \label{eq:2nd-moment-case5-pf-4}
T_{\star}^i := \max \left\{j \geq 1 :
D_j^{x, i} < R_m^{y, k} \right \},
~~w_n \leq i \leq k - d_n (\theta) - 1.
\end{equation}
Fix $\ell \in \left\{0, \cdots, m_n^- - 1 \right \}$.
By the strong Markov property at
$R_m^{y, k}$,
the $\ell$-th term of the sum in (\ref{eq:2nd-moment-case5-pf-3})
is equal to
\begin{equation}
\label{eq:2nd-moment-case5-pf-5}
\begin{aligned}
E_o &\Biggl[
1_{\left\{
D_{\ell}^{x, 0} < R_m^{y, k} < D_{\ell + 1}^{x, 0},~
H_x > R_m^{y, k},~H_y > D_{m - 1}^{y, k}
\right \}} \\
&\,\,\,\,\,\times
P_{S_{R_m^{y, k}}} \Biggl[
\bigcap_{i=w_n}^{k - d_n (\theta) - 1}
\left\{b_n^- (i) \leq T_{\star}^i + T_i^{x, m_n^- - \ell} \leq b_n^+ (i) \right\} \\
&\,\,\,\,\,\,\,\,\,\,\,\,\,\,\,\,\,\,\,\,\,\,\,\,\,\,\,\,\,\,\,\,\,\,\,\,\,
\,\,\,\,\,\,\,\,\,\,\,\,\,\,\,\,\,\,\,\,\,\,\,\,\,\,\,
\cap
\left\{H_x > D_{m_n^- - \ell}^{x, 0},~H_y > H_{\partial B (y, r_k)} \right\}
\Biggr]
\Biggr].
\end{aligned}
\end{equation}  

Fix any $z \in \partial B (y, r_{k + 1})$
and $t_i \geq 0$, $w_n \leq i \leq k - d_n (\theta) - 1$.
By the strong Markov property
at $H_{\partial B (x, r_{k - 2})}$ 
and at $H_{\partial B (x, r_{k - d_n (\theta)})}$,
the probability 
in (\ref{eq:2nd-moment-case5-pf-5})
with $S_{R_m^{y, k}} = z$
and $T_{\star}^i = t_i$, $w_n \leq i \leq k - d_n (\theta) - 1$
is equal to
\begin{align} \label{eq:2nd-moment-case5-pf-6}
&\sum_{u \in \partial B (x, r_{k - 2})}
\sum_{v \in \partial B (x, r_{k - d_n (\theta)})}
P_z \left[
\begin{minipage}{160pt}
$H_y > H_{\partial B (y, r_k)},~H_x > H_{\partial B (x, r_{k - 2})}$, \\
$S_{H_{\partial B (x, r_{k - 2})}} = u$
\end{minipage}
\right] \notag \\
&\times P_u \left[
H_x > H_{\partial B (x, r_{k - d_n (\theta)})},~
S_{H_{\partial B (x, r_{k - d_n (\theta)})}} = v \right] \notag \\
&\times 
P_v \left[\bigcap_{i=w_n}^{k-d_n(\theta)-1}
\left\{b_n^- (i) \leq t_i + T_i^{x, m_n^- - \ell} \leq b_n^+ (i) \right\} \cap
\left\{H_x > D_{m_n^- - \ell}^{x, 0} \right\}
\right].
\end{align}
Fix $u \in \partial B (x, r_{k - 2})$
and $v \in \partial B (x, r_{k - d_n (\theta)})$.
We will estimate 
the probability in the second line of (\ref{eq:2nd-moment-case5-pf-6}).
By the strong Markov property
at $H_x$
and at $H_{\partial B (x, r_{k - 2})} \circ \theta_{H_x}$,
one can see that
$$P_u \left[H_x < H_{\partial B (x, r_{k - d_n (\theta)})},~
S_{H_{\partial B (x, r_{k - d_n (\theta)})}} = v \right]$$
is equal to
\begin{equation} \label{eq:2nd-moment-case5-pf-7}
P_u \left[
H_x < H_{\partial B (x, r_{k - d_n (\theta)})} 
\right]
E_x \left[
P_{S_{H_{\partial B (x, r_{k - 2})}}}
\left[S_{H_{\partial B (x, r_{k - d_n (\theta)})}}^{\prime}
 = v \right]
\right].
\end{equation}
By the Harnack inequality (\ref{eq:uniform-estimate-harmonic-measure}),
(\ref{eq:2nd-moment-case5-pf-7})
is equal to
$\left\{1 + O \left((\ell_n)^2 (\log n)^{- \theta} \right) \right \}$ times
\begin{equation} \label{eq:2nd-moment-case5-pf-8}
P_u \left[H_x < H_{\partial B (x, r_{k - d_n (\theta)})}
 \right]
P_u \left[S_{H_{\partial B (x, r_{k - d_n (\theta)})}} = v \right].
\end{equation}
By (\ref{eq:2nd-moment-case5-pf-8}) and Lemma \ref{lem:circle-to-circle-probab},
the probability in the second line of
(\ref{eq:2nd-moment-case5-pf-6})
is equal to
$\{1 + O ((\ell_n)^3 (\log n)^{- \theta}) \}$ times
\begin{equation} \label{eq:2nd-moment-case5-pf-9}
P_u \left[H_x > H_{\partial B (x, r_{k - d_n (\theta)})} \right]
P_u \left[S_{H_{\partial B (x, r_{k - d_n (\theta)})}} = v \right].
\end{equation}
By the Harnack inequality (\ref{eq:uniform-estimate-harmonic-measure}),
the second probability in (\ref{eq:2nd-moment-case5-pf-9}) is equal to
\begin{equation} \label{eq:2nd-moment-case5-pf-10}
\left\{1 + O \left((\ell_n)^2 (\log n)^{- \theta} \right) \right \}
P_{u_{\star}}
\left[S_{H_{\partial B (x, r_{k - d_n (\theta)})}} = v \right],
\end{equation}
where $u_{\star}$ is a fixed vertex on $\partial B (x, r_{k - 2})$.
By Lemma \ref{lem:circle-to-circle-probab},
the first probability in (\ref{eq:2nd-moment-case5-pf-9})
is equal to
\begin{equation} \label{eq:2nd-moment-case5-pf-11}
\left\{1 + O \left(\frac{1}{(L_n - k)^2 (\log (\ell_n))^2} \right) \right \}
\frac{L_n - k + 2}{L_n - k + d_n (\theta)}.
\end{equation} 
By these and the condition that $\theta$ is large enough,
(\ref{eq:2nd-moment-case5-pf-6}) is bounded from above by
\begin{equation} 
\label{eq:2nd-moment-case5-pf-12}
\begin{aligned}
&\left\{1 + O \left(\frac{1}{(L_n - k)^2 (\log (\ell_n))^2} \right) \right \}
\frac{L_n - k + 2}{L_n - k + d_n (\theta)}
\\
&\times P_z \left[H_y > H_{\partial B (y, r_k)},~H_x > H_{\partial B (x, r_{k - 2})} \right] \\
&\times P_{u_{\star}}
\Biggl[
\bigcap_{i=w_n}^{k - d_n (\theta) - 1}
\left\{b_n^- (i) \leq t_i + T_i^{x, m_n^- - \ell} \leq b_n^+ (i) \right\} \\
&\,\,\,\,\,\,\,\,\,\,\,\,\,\,\,\,\,\,\,\,\,\,\,\,\,\,\,\,\,\,\,\,\,\,\,\,\,\,\,\,
\,\,\,\,\,\,\,\,\,\,\,\,\,\,\,
\cap
\left\{H_x \circ \theta_{H_{\partial B (x, r_{k - d_n (\theta)})}}
> D_{m_n^- - \ell}^{x, 0} \right\}
\Biggr].
\end{aligned}
\end{equation}
We will estimate the first probability in (\ref{eq:2nd-moment-case5-pf-12}).
By the strong Markov property
at $H_y$ and at $H_{\partial B (y, r_k)} \circ \theta_{H_y}$,
we have
\begin{align} \label{eq:2nd-moment-case5-pf-13}
&P_z \left[H_y < H_{\partial B (y, r_k)},~H_x > H_{\partial B (x, r_{k - 2})} \right] \notag \\
&= P_z \left[H_y < H_{\partial B (y, r_k)} \right]
E_y \left[P_{S_{H_{\partial B (y, r_k)}}} \left[H_x > H_{\partial B (x, r_{k - 2})} \right] \right].
\end{align}
By Lemma \ref{lem:circle-to-circle-probab}
and the fact that $d (u, x) \geq r_k$
for any $u \in \partial B (y, r_k)$,
the right of (\ref{eq:2nd-moment-case5-pf-13})
is bounded from below by
$\left(1 - \frac{c_1}{L_n - k} \right) \frac{1}{L_n - k}$.
Thus, we have
\begin{align} \label{eq:2nd-moment-case5-pf-14}
&P_z \left[H_y > H_{\partial B (y, r_k)},~H_x > H_{\partial B (x, r_{k - 2})} \right] \notag \\
&= P_z \left[H_x > H_{\partial B (x, r_{k - 2})} \right]
- P_z \left[H_y < H_{\partial B (y, r_k)},~H_x > H_{\partial B (x, r_{k - 2})} \right] \notag \\
&\leq \left\{1 - \left(1 - \frac{c_1}{L_n - k} \right) \frac{1}{L_n - k} \right \}
P_z \left[H_x > H_{\partial B (x, r_{k - 2})} \right].
\end{align}
By these and applying the strong Markov property
at $R_m^{y, k}$
and Lemma \ref{lem:circle-to-circle-probab}
to the last term of (\ref{eq:2nd-moment-case5-pf-3}),
$P_o [C_m]$ is bounded from above by
\begin{align} \label{eq:2nd-moment-case5-pf-15}
&\left\{1 + \frac{c_2}{(L_n - k)^2 (\log (\ell_n))^2} \right \}
\left\{1 - \left(1 - \frac{c_1}{L_n - k} \right) \frac{1}{L_n - k} \right \}
\notag \\
&\times \left\{
\sum_{\ell = 0}^{m_n^- - 1}
\frac{L_n - k + 2}{L_n - k + d_n (\theta)} 
E_{\ell} + Q \right\},
\end{align}
where
\begin{equation*}
E_{\ell} := E_o \left[
\begin{minipage}{180pt}
$1_{\left\{
D_{\ell}^{x, 0} < R_m^{y, k} < D_{\ell + 1}^{x, 0},~
H_x > R_m^{y, k},~
H_y > D_{m - 1}^{y, k}
\right \}}$ \\
$\times P_{u_{\star}}^{(\ell)} \left((T_{\star}^i)_i \right)
P_{S_{R_m^{y, k}}} \left[H_x > H_{\partial B (x, r_{k - 2})} \right]$
\end{minipage}
\right],
\end{equation*}
(we let $P_{u^{\star}}^{(\ell)} ((t_i)_i)$
be the probability in the last line of (\ref{eq:2nd-moment-case5-pf-12}))
\begin{equation*}
Q := P_o \left[
\left\{
D_{m_n^-}^{x, 0} < R_m^{y, k} \right\} \cap
C_{m-1}
\right].
\end{equation*}
A similar argument implies that
$P_o [C_{m - 1}]$
is bounded from below by 
the factor in the second line of (\ref{eq:2nd-moment-case5-pf-15})
multiplied by
$1 - c_3 (L_n - k)^{- 2} (\log (\ell_n))^{- 2}$.
Thus, 
$P_o [C_m]$ is bounded from above by
\begin{equation} \label{eq:2nd-moment-case5-pf-16}
\left(1 + \frac{c_4}{(L_n - k)^2 (\log (\ell_n))^2} \right) 
\left \{1 - \left(1 - \frac{c_1}{L_n - k} \right) \frac{1}{L_n - k} \right \}
P_o [C_{m - 1}].
\end{equation}
Thus, the right of (\ref{eq:2nd-moment-case5-pf-1})
is bounded from above by
\begin{align} \label{eq:2nd-moment-case5-pf-17}
&\left\{1 + \frac{c_4}{(L_n - k)^2 (\log (\ell_n))^2} \right \}^{b_n^- (k)}
\left\{1 - \left(1 - \frac{c_1}{L_n - k} \right) \frac{1}{L_n - k} \right \}^{b_n^- (k)}
 \notag \\
& \times
P_o \left[
b_n^{-} (i) \leq T_i^{x, m_n^{-}} \leq b_n^{+} (i),
~w_n \leq \forall i \leq k - d_n (\theta) - 1,~
H_x > D_{m_n^{-}}^{x, 0}
\right].
\end{align}
The first factor of (\ref{eq:2nd-moment-case5-pf-17}) is $1 + o(1)$.
The second factor of (\ref{eq:2nd-moment-case5-pf-17})
is bounded from above by
$(\ell_n)^{c_5} e^{- \frac{b_n^- (k)}{L_n - k}}$
and this is bounded from above by
\begin{equation} \label{eq:2nd-moment-case5-pf-17-1}
c_6 (\ell_n)^{- 2 (L_n - k) + c_5}
e^{(L_n - k) \frac{\log \log n}{\log n} \log (\ell_n)
+ (L_n - k) \frac{2 s_n}{\log n} \log (\ell_n)
- c_7 f_n (k) \sqrt{\log (\ell_n)}}.
\end{equation}
Since $k > \lceil (1 - \varepsilon) L_n \rceil > \frac{L_n - 1}{2}$,
we have
$$\frac{c_7}{2} f_n (k) \sqrt{\log (\ell_n)} - (L_n - k) \frac{\log \log n}{\log n} \log (\ell_n) > 0.$$
By this and the condition
$k > \lceil (1 - \varepsilon) L_n \rceil$,
(\ref{eq:2nd-moment-case5-pf-17-1})
is bounded from above by
\begin{equation} \label{eq:2nd-moment-case5-pf-17-2}
c_6 (\ell_n)^{- 2 (L_n - k) + c_5} 
e^{2 s_n - (1 - \varepsilon) s_n - \frac{c_7}{2} f_n (k) \sqrt{\log (\ell_n)}}.
\end{equation}
By the transfer lemma (Lemma \ref{lem:transfer-lemma}),
the probability in (\ref{eq:2nd-moment-case5-pf-17}) is bounded from above by
\begin{equation} \label{eq:2nd-moment-case5-pf-22}
(1 + o(1))
P_{m_n^-}^{\text{GW}} \left[\bigcap_{i=w_n}^{k-d_n(\theta)-1}
\left\{b_n^- (i) \leq T_i \leq b_n^+ (i) \right\} \cap
\left\{T_{L_n - 1} = 0 \right\}
\right].
\end{equation}
By the Markov property,
the probability in (\ref{eq:2nd-moment-case5-pf-22})
is equal to
\begin{equation} 
\label{eq:2nd-moment-case5-pf-22-1} 
\begin{aligned}
\sum_{m, m^{\prime}}
&P_{m_n^-}^{\text{GW}} [T_{w_n} = m] \\
&\times P_m^{\text{GW}} \Biggl[
\bigcap_{i=1}^{k - d_n (\theta) - w_n - 2}
\left\{b_n^- (w_n + i) \leq T_i \leq b_n^+ (w_n + i) \right\} \\
&\,\,\,\,\,\,\,\,\,\,\,\,\,\,\,\,\,\,\,\,\,\,\,\,\,\,\,\,\,\,\,\,\,\,\,\,\,\,\,\,\,\,\,\,\,\,\,\,\,\,
\,\,\,\,\,\,\,\,\,\,\,\,\,\,\,\,\,\,\,\,\,\,\,\,\,\,\,\,\,\,\,\,\,\,\,\,\,\,\,\,\,\,\,\,\,\,\,\,\,\,
\cap \left\{T_{k - d_n (\theta) - w_n - 1} = m^{\prime} \right\}
\Biggr] \\
&\times P_{m^{\prime}}^{\text{GW}} [T_{L_n - k + d_n (\theta)} = 0],
\end{aligned}
\end{equation}
where the sum is taken over
$m \in \{b_n^- (w_n), \dotsc, b_n^+ (w_n) \}$
and $m^{\prime} \in \{b_n^- (k - d_n (\theta) - 1), \dotsc, b_n^+ (k - d_n (\theta) - 1) \}$.
By the barrier estimate (Lemma \ref{lem:GW-process-barrier-estimate})
with $a = (1 - \frac{w_n}{L_n}) \sqrt{2 m_n^-}$,
$b = \frac{L_n - k + d_n (\theta) + 1}{L_n} \sqrt{2 m_n^-}$,
$x = \sqrt{2m}$, $y = \sqrt{2 m^{\prime}}$
(we can take $\eta = c \sqrt{\log (\ell_n)}$ for some positive constant),
the second probability in (\ref{eq:2nd-moment-case5-pf-22-1})
is bounded from above by
\begin{equation} \label{eq:2nd-moment-case5-pf-22-2}
c_8 (\ell_n)^{c_9} \frac{(L_n - k + d_n (\theta) + 1)^{\frac{1}{2} + \delta}}{\log n}
e^{- \frac{(\sqrt{m} - \sqrt{m^{\prime}})^2}{k - d_n (\theta) - w_n - 1}}.
\end{equation}
The third probability in (\ref{eq:2nd-moment-case5-pf-22-1})
is equal to 
$(1 - \frac{1}{L_n - k + d_n (\theta) + 1})^{m^{\prime}}$
and this is bounded from above by
$e^{- \frac{m^{\prime}}{L_n - k + d_n (\theta) + 1}}$.
The product of this and the exponential factor in 
(\ref{eq:2nd-moment-case5-pf-22-2}) is bounded from above by
$e^{- \frac{m}{L_n - w_n}}$.
By this, 
the contribution of the sum over $m$ in (\ref{eq:2nd-moment-case5-pf-22-1})
comes from $\sum_{m = b_n^- (w_n)}^{b_n^+ (w_n)} 
e^{- \frac{m}{L_n - w_n}} P_{m_n^-}^{\text{GW}} [T_{w_n} = m]$ and 
this is bounded from above by 
\begin{equation} \label{eq:2nd-moment-case5-pf-22-3}
E_{m_n^-}^{\text{GW}} \left[1_{\{T_{w_n} \geq b_n^- (w_n) \}} e^{- \frac{T_{w_n}}{L_n - w_n}} \right]
\leq e^{- \frac{b_n^- (w_n)}{L_n - w_n}}.
\end{equation}
Therefore, (\ref{eq:2nd-moment-case5-pf-22}) is bounded from above by
\begin{equation} \label{eq:2nd-moment-case5-pf-22-4}
c_{10} n^{- 2} e^{2 s_n} (\ell_n)^{- 2 c_{\star} w_n + c_{11} w_n}
(L_n - k + d_n (\theta) + 1)^{\frac{5}{2} + \delta}.
\end{equation}

Note that the number of pairs
of points $x, y \in \mathbb{Z}_n^2 \backslash B (o, r_0)$
with $\ell (x, y) = k$
is bounded from above by $n^2 (2 r_{k - 1})^2$.
Thus, 
by (\ref{eq:2nd-moment-case5-pf-17}), (\ref{eq:2nd-moment-case5-pf-17-2}),
(\ref{eq:2nd-moment-case5-pf-22-4}),
and Lemma \ref{lem:lower-bound-first-moment},
the left of (\ref{eq:2nd-moment-case5})
is bounded from above by
\begin{align} \label{eq:2nd-moment-case5-pf-27}
c_{12} (\ell_n)^{ c_{13} w_n} 
e^{- (1 - \varepsilon) s_n}
\left\{E_o [Z_n] \right \}^2.
\end{align}
Since $\gamma > \alpha + \beta$
by the assumption (\ref{eq:assumption-parameter}),
we have the desired result.
~~~$\qed$ \\

Next, we deal with the case $\ell (x, y) = 0$.
\begin{lem} \label{lem:2nd-moment-case1}
As $n \to \infty$,
\begin{equation} \label{eq:2nd-moment-case1}
\sum_{\begin{subarray}{c} x, y \in \mathbb{Z}_n^2 \backslash B (o, r_0) \\ \ell (x, y) = 0 \end{subarray}}
P_o [A_n (x) \cap A_n (y)] \leq (1 + o (1)) \left(E_o [Z_n] \right)^2.
\end{equation}
\end{lem}
{\it Proof.}
Fix any $x, y \in \mathbb{Z}_n^2 \backslash B (o, r_0)$
with $\ell (x, y) = 0$.
Since $A_n (z)$ ($z \in \{x, y\}$) can be written as an event
for $m_n^-$ excursions from $\partial B (z, r_1)$ to $\partial B (z, r_0)$,
we can apply Lemma \ref{lem:decomposition-trajectories}
with $R = r_0$, $r = r_1$, $k = 0$, $\ell = m = m_n^-$.
By this and the transfer lemma
(Lemma \ref{lem:transfer-lemma}(i)),
one can show that $P_o [A_n (x) \cap A_n (y)]$ is bounded from above by
\begin{equation} \label{eq:lem-2nd-moment-case1-1}
(1 + o (1)) P_{m_n^-}^{\text{GW}} \left[\bigcap_{i=w_n}^{L_n-w_n-1} 
\left\{b_n^{-} (i) \leq T_i \leq b_n^{+} (i) \right\}
\cap \left\{T_{L_n - 1} = 0 \right\} \right]^2,
\end{equation}
where $o (1) \to 0$
as $n \to \infty$ uniformly in $x, y$.
By (\ref{eq:lem-lower-bound-first-moment-3}),
the right-hand side of (\ref{eq:lem-2nd-moment-case1-1})
is bounded from above by
\begin{equation}
\label{eq:lem-2nd-moment-case1-1-almost-sure-independence}
(1 + o(1)) P_o [A_n (x)] P_o [A_n (y)].
\end{equation}
Therefore, we have (\ref{eq:2nd-moment-case1}). $\qed$ \\\\

{\it Proof of (\ref{eq:lower-main-prop}).}
By Lemma \ref{lem:lower-bound-first-moment},
we have
\begin{align} \label{eq:main-lower-bound-pf-step1}
\sum_{\begin{subarray}{c} x, y \in \mathbb{Z}_n^2 \backslash B(o, r_0) : 
\\ \ell (x, y) \geq L_n - w_n \end{subarray}}
P_o [A_n (x) \cap A_n (y)]
&\leq c_1 (\ell_n)^{2 (w_n + 1)} E_o [Z_n]  \notag \\
&\leq c_1 e^{- 2 s_n} (\ell_n)^{c_2 w_n}
\{E_o [Z_n] \}^2.
\end{align}
Since $\gamma > \alpha + \beta$ by the assumption (\ref{eq:assumption-parameter}),
the right of (\ref{eq:main-lower-bound-pf-step1})
is equal to $o(1) E_o [Z_n]^2$.
By this and
and Lemmas \ref{lem:lower-bound-first-moment},
\ref{lem:2nd-moment-case3},
\ref{lem:2nd-moment-case2},
\ref{lem:2nd-moment-case4},
\ref{lem:2nd-moment-case5},
\ref{lem:2nd-moment-case1},
we have
\begin{align} \label{eq:main-lower-bound-pf-step2}
&P_o \left[
\exists x \in \mathbb{Z}_n^2 \backslash B (o, r_0),
H_x > D_{m_n^{-}}^{x, 0} \right] \notag \\
&\geq P_o [Z_n \geq 1]
\geq \frac{\left \{E_o [Z_n] \right \}^2}{E_o [Z_n^2]} = 1 + o (1),
~~\text{as}~n \to \infty.~~~\qed
\end{align}

{\it Proof of the lower bound of Theorem \ref{thm:main} via (\ref{eq:lower-time}).}
(\ref{eq:lower-main-prop}) and (\ref{eq:lower-time})
immediately yield the lower bound. $\qed$

\appendix \section{Excursion length} \label{sec:appendix}
In this section, we give a proof of Proposition \ref{prop:time-number}.
We follow arguments in \cite[Section 8]{BeKi}.
Recall definitions of 
$R_m (x, R, r)$, $D_m (x, R, r)$
from (\ref{eq:stopping-times}).
\begin{lem} (\cite[Lemma 2.1]{Ue}) \label{lem:appendix-equilibrium-law}
Fix $n \in \mathbb{N}$ and $1 < r < R/4 < n/8$.
For any $y \in \mathbb{Z}_n^2$,
there exists a pair of probability measures
$\mu_r^{y, R}$
on $\partial B (y, R)$
and $\mu_R^{y, r}$
on $\partial B (y, r)$
such that
\begin{equation} \label{eq:appendix-equilibrium-law-statement1}
\mu_r^{y, R} (z) = P_{\mu_R^{y, r}} [S_{H_{\partial B (y, R)}} = z],~z \in \partial B (y, R),
\end{equation}
\begin{equation} \label{eq:appendix-equilibrium-law-statement2}
\mu_R^{y, r} (z) = P_{\mu_r^{y, R}} [S_{H_{\partial B (y, r)}} = z],~z \in \partial B (y, r).
\end{equation}
\end{lem}
\begin{lem} \label{lem:appendix-expectation-d1}
Fix $n \in \mathbb{N}$ and $1 < r < R/4 < n/8$.
For any $y \in \mathbb{Z}_n^2$,
\begin{equation} \label{eq:lem-appendix-expectation-d1-statement}
E_{\mu_r^{y, R}} [D_1 (y, R, r)] = \frac{2}{\pi} n^2 
\left\{\log \left(\frac{R}{r} \right) + O \left(\frac{1}{r} \right) \right \}.
\end{equation}
\end{lem}
{\it Proof.}
We define the measure $m$ on $\mathbb{Z}_n^2$ by
\begin{equation} \label{eq:lem-appendix-expectation-d1-step1}
m (\cdot) := \sum_{v \in \partial B (y, r)} \mu_R^{y, r} (v) G_{B (y, R)} (v, \cdot)
+ \sum_{v \in \partial B (y, R)} \mu_r^{y, R} (v) G_{B (y, r)} (v, \cdot),
\end{equation}
where 
$$G_{B (y, R^{\prime})} (u, v)
:= E_u \left[\sum_{i = 0}^{H_{\partial B (y, R^{\prime})} - 1}
1_{\{S_i = v \}} \right],~~~0 < R^{\prime} < n/2,~u, v \in \mathbb{Z}_n^2.$$ 
One can easily check that
$m$ is the stationary measure on $\mathbb{Z}_n^2$.
Thus, there exists $c_0 > 0$ such that
$m = c_0 \nu$
where $\nu$ is the uniform measure on $\mathbb{Z}_n^2$.
By Green's function estimate (see, for example, \cite[Proposition 1.6.7]{La}),
we have
\begin{equation} \label{eq:lem-appendix-expectation-d1-step2}
m (y) = \sum_{v \in \partial B (y, r)} \mu_R^{y, r} (v) G_{B (y, R)} (v, y)
= \frac{2}{\pi} \log \left(\frac{R}{r} \right) + O \left(\frac{1}{r} \right).
\end{equation}
By (\ref{eq:appendix-equilibrium-law-statement2}), we have
\begin{align} \label{eq:lem-appendix-expectation-d1-step3}
&E_{\mu_r^{y, R}} [D_1 (y, R, r)] \notag \\
&= E_{\mu_r^{y, R}} [H_{\partial B (y, r)}] + E_{\mu_r^{y, R}} 
\left[E_{S_{H_{\partial B (y, r)}}} [H_{\partial B (y, R)}] \right] \notag \\
&= E_{\mu_r^{y, R}} [H_{\partial B (y, r)}]
+ E_{\mu_R^{y, r}} \left[H_{\partial B (y, R)} \right] \notag \\
&= \sum_{u \in \mathbb{Z}_n^2} m (u) 
= c_0.
\end{align}
By (\ref{eq:lem-appendix-expectation-d1-step2}),
(\ref{eq:lem-appendix-expectation-d1-step3}),
and the fact that $m (y) = c_0 \nu (y) = c_0 n^{- 2}$,
we have (\ref{eq:lem-appendix-expectation-d1-statement}). $\qed$

\begin{lem} \label{lem:appendix-tail-d1}
There exist $c_1, c_2, c_3 \in (0, \infty)$
such that
for any $n \in \mathbb{N}$,
$\lambda \geq 0$,
$x, y \in \mathbb{Z}_n^2$,
and $c_1 \leq r < R/4 < n/8$,
\begin{equation} \label{eq:lem-appendix-tail-d1-statement}
P_x [D_1 (y, R, r) \geq \lambda]
\leq c_2 \exp \left\{- c_3 \frac{\lambda}{n^2 \log (n/r)} \right \}.
\end{equation}
\end{lem}
{\it Proof.}
By the exponential Chebyshev inequality,
for any $\theta > 0$,
the left of (\ref{eq:lem-appendix-tail-d1-statement})
is bounded from above by
\begin{equation} \label{eq:lem-appendix-tail-d1-pf-1}
e^{- \theta \lambda} E_x \left[\exp \left\{\theta D_1 (y, R, r) \right\} \right].
\end{equation}
By the strong Markov property at $H_{\partial B (y, r)}$,
(\ref{eq:lem-appendix-tail-d1-pf-1}) 
is bounded from above by
\begin{equation} \label{eq:lem-appendix-tail-d1-pf-2}
e^{- \theta \lambda}
E_x \left[\exp \left \{\theta H_{\partial B (y, r)} \right \} \right]
\cdot \max_{z \in \partial B (y, r)} 
E_z \left[\exp \left \{\theta H_{\partial B (y, R)} \right \} \right].
\end{equation} 
To estimate the exponential moments,
we use the following:
\begin{itemize}
\item \cite[Lemma 3.1]{DPRZ2}
There exist $c, c^{\prime} \in (0, \infty)$
such that for any $c \leq r < n/6$,
\begin{equation} \label{eq:lem-appendix-tail-d1-pf-3}
\max_{y \in \mathbb{Z}_n^2}
\max_{x \in \mathbb{Z}_n^2}
E_x \left[H_{\partial B (y, r)} \right]
\leq c^{\prime} n^2 \log (n/r).
\end{equation}
\item \cite[(1.21)]{La}
For any $z \in B (y, R)$,
\begin{equation} \label{eq:lem-appendix-tail-d1-pf-4}
R^2 - (d (z, y))^2
\leq E_z \left[H_{\partial B (y, R)} \right]
\leq (R + 1)^2 - (d (z, y))^2.
\end{equation}
\item Kac's moment formula
\cite[(6)]{FP}
For any first hitting time $T$,
$m \in \mathbb{N}$,
and $y \in \mathbb{Z}_n^2$,
\begin{equation} \label{eq:lem-appendix-tail-d1-pf-5}
E_y [T^m]
\leq m! E_y [T]
\left(\max_{z \in \mathbb{Z}_n^2} E_z [T] \right)^{m - 1}.
\end{equation}
\end{itemize}
By (\ref{eq:lem-appendix-tail-d1-pf-3})
and (\ref{eq:lem-appendix-tail-d1-pf-5}),
we have
\begin{align} \label{eq:lem-appendix-tail-d1-pf-6}
E_x \left[\exp \left\{\theta H_{\partial B (y, r)} \right\} \right] 
&= 1 +
\sum_{k = 1}^{\infty}
\frac{\theta^k}{k !}
E_x \left[\left(H_{\partial B (y, r)} \right)^k \right] \notag \\
&\leq 1 + \sum_{k = 1}^{\infty}
\theta^k
\left(\max_{z \in \mathbb{Z}_n^2} E_z \left[H_{\partial B (y, r)} \right] \right)^k \notag \\
&\leq 1 + \sum_{k = 1}^{\infty}
\left\{\theta c_2 n^2 \log (n/r) \right \}^k.
\end{align}
Taking $\theta \leq \frac{1}{2} (c_2 n^2 \log (n/r))^{- 1}$,
we have
\begin{equation} \label{eq:lem-appendix-tail-d1-pf-7}
E_x \left[\exp \left\{\theta H_{\partial B (y, r)} \right \} \right] \leq 2.
\end{equation}
Fix any $z \in \partial B (y, r)$.
Similarly, by (\ref{eq:lem-appendix-tail-d1-pf-4})
and (\ref{eq:lem-appendix-tail-d1-pf-5}),
we have
\begin{align} \label{eq:lem-appendix-tail-d1-pf-8}
E_z \left[\exp \left\{\theta H_{\partial B (y, R)} \right\} \right]
&= E_z \left[\exp \left\{\theta H_{B (y, R)^c} \right\} \right] \notag \\
&\leq 1 + \sum_{k = 1}^{\infty}
\left\{\theta \max_{v \in B (y, R)}
E_v \left[H_{B (y, R)^c} \right] \right\}^k \notag \\
& \leq 1 + \sum_{k = 1}^{\infty}
\left\{\theta (R + 1)^2 \right \}^k.
\end{align} 
Taking $\theta \leq \frac{1}{2} (R + 1)^{- 2}$,
we have
\begin{equation} \label{eq:lem-appendix-tail-d1-pf-9}
E_z \left[\exp \left\{\theta H_{\partial B (y, R)} \right\} \right]
\leq 2.
\end{equation}
Note that $R < n/2$ and $r < n/8$
imply
$(R + 1)^2 \leq n^2 \log (n/r)$.
Thus, taking
$\theta = \frac{1}{2} \left((c_2 \vee 1) n^2 \log (n/r) \right)^{- 1}$,
we have
\begin{equation} \label{eq:lem-appendix-tail-d1-pf-10}
P_x \left[D_1 (y, R, r) \geq \lambda \right]
\leq 4 \exp \left\{- \frac{\lambda}{2 (c_2 \vee 1) n^2 \log (n/r)} \right\}.~~~\qed
\end{equation}

\begin{lem} \label{lem:appendix-sum-exit-time}
There exist $c_1, c_2, c_3 \in (0, \infty)$ such that
for any $n \in \mathbb{N}$, $1 < r < R/4 < n/8$,
$x, y \in \mathbb{Z}_n^2$,
$m \geq 2$, and $\delta \in (0, 1)$,
\begin{align} \label{eq:lem-appendix-sum-exit-time-statement}
&P_x \left[
\Biggl|
\frac{\sum_{i = 2}^m \left(D_i (y, R, r) - R_i (y, R, r) \right)}
{E_{\mu_R^{y, r}} \left[H_{\partial B (y, R)} \right] (m - 1)}
- 1 \Biggr| > \delta \right] \notag \\
&\leq c_1 \exp \left \{
- c_2 \delta^2 (m - 1) + c_3 \delta \frac{1}{R} (m - 1)
\right \}.
\end{align}
\end{lem}
{\it Proof.}
By the exponential Chebyshev inequality,
for any $\theta > 0$,
we have
\begin{align} \label{eq:lem-appendix-sum-exit-time-pf-1}
&P_x \left[
\sum_{i = 2}^m \left(D_i (y, R, r) - R_i (y, R, r) \right)
\geq (1 + \delta) E_{\mu_R^{y, r}} \left[H_{\partial B (y, R)} \right] (m - 1) 
\right] \notag \\
&\leq
e^{- \theta (1 + \delta) E_{\mu_R^{y, r}} \left[H_{\partial B (y, R)} \right] (m - 1)}
E_x \left[
\exp \left \{
\theta \sum_{i = 2}^m \left(D_i (y, R, r) - R_i (y, R, r) \right)
\right \}
\right].
\end{align}
By the strong Markov property at $R_i (y, R, r)$, $2 \leq i \leq m$,
the expectation on the right of (\ref{eq:lem-appendix-sum-exit-time-pf-1})
is bounded from above by
\begin{equation} \label{eq:lem-appendix-sum-exit-time-pf-2}
\left(
\max_{z \in \partial B (y, r)}
E_z \left[\exp \left\{\theta H_{\partial B (y, R)} \right \} \right] 
\right)^{m - 1}.
\end{equation}
By (\ref{eq:lem-appendix-tail-d1-pf-4}) and
Kac's moment formula (\ref{eq:lem-appendix-tail-d1-pf-5}), 
for any $z \in \partial B (y, r)$ and $\theta \leq 1/(2 (R + 1)^2)$,
we have
\begin{align} \label{eq:lem-appendix-sum-exit-time-pf-3}
E_z \left[\exp \{\theta H_{\partial B (y, R)} \} \right]
& = E_z \left[\exp \{\theta H_{B (y, R)^c} \} \right] \notag \\
&\leq 1 + \theta E_z \left[H_{B (y, R)^c} \right]
+ \sum_{k = 2}^{\infty}
\left(\theta \max_{w \in B (y, R)}
E_w \left[H_{B (y, R)^c} \right]
\right)^k \notag \\
&\leq 1 + \theta \left \{(R + 1)^2 - r^2 \right \}
+ \sum_{k = 2}^{\infty} \left \{\theta (R + 1)^2 \right \}^k \notag \\
&\leq \exp \left \{\theta \left((R + 1)^2 - r^2 \right)
+ 2 \theta^2 (R + 1)^4
\right \}.
\end{align}
Note that (\ref{eq:lem-appendix-tail-d1-pf-4})
implies that 
$E_{\mu_R^{y, r}} \left[H_{\partial B (y, R)} \right] \geq R^2 - (r + 1)^2$.
By this, (\ref{eq:lem-appendix-sum-exit-time-pf-2}),
and(\ref{eq:lem-appendix-sum-exit-time-pf-3}),
for $\theta = \frac{\delta}{8} \frac{R^2 - (r + 1)^2}{(R + 1)^4}$,
the right of (\ref{eq:lem-appendix-sum-exit-time-pf-1})
is bounded from above by a value
of the form of the right of (\ref{eq:lem-appendix-sum-exit-time-statement}).
By almost the same argument,
one can obtain a similar estimate
of the left tail.
We omit the detail. $\qed$

For $y \in \mathbb{Z}_n^2$ and $0 < r < R < n/2$,
set
\begin{equation} \label{eq:appendix-def-q}
q (y, R, r) := 
\min_{\begin{subarray}{c} u \in \partial B (y, R), \\ v \in \partial B (y, r) \end{subarray}}
\frac{P_v [S_{H_{\partial B (y, R)}} = u]}{\mu_r^{y, R} (u)}.
\end{equation}
\begin{lem} \label{lem:appendix-estimate-q}
There exists $c_1 > 0$ such that
for any $n \in \mathbb{N}$, $x \in \mathbb{Z}_n^2$,
and $1 < r < R/4 < n/8$,
\begin{equation} \label{eq:lem-appendix-estimate-q-statement}
q (y, R, r) \geq 1 - c_1 \frac{r}{R}.
\end{equation}
\end{lem}
{\it Proof.}
Fix any $v \in \partial B (y, r)$ and $u \in \partial B (y, R)$.
By Lemma \ref{lem:appendix-equilibrium-law},
we have
\begin{equation} \label{eq:lem-appendix-estimate-q-pf-1}
\mu_r^{y, R} (u) = \sum_{z \in \partial B (y, r)} \mu_R^{y, r} (z)
P_z \left[S_{H_{\partial B (y, R)}} = u \right]
\leq \max_{z \in \partial B (y, r)} P_z \left[S_{H_{\partial B (y, R)}} = u \right].
\end{equation}
By the Harnack inequality (\ref{eq:uniform-estimate-harmonic-measure}),
for any $z \in \partial B (y, r)$,
\begin{equation} \label{eq:lem-appendix-estimate-q-pf-2}
P_z \left[S_{H_{\partial B (y, R)}} = u \right]
= \left(1 + O\left(\frac{r}{R} \right) \right)
P_v \left[S_{H_{\partial B (y, R)}} = u \right].
\end{equation}
By (\ref{eq:lem-appendix-estimate-q-pf-1}) and (\ref{eq:lem-appendix-estimate-q-pf-2}),
we have (\ref{eq:lem-appendix-estimate-q-statement}). $\qed$ \\\\

For $y \in \mathbb{Z}_n^2$, $1 < r < R/4 < n/8$,
let $W(y, R, r)$ be the space of paths from $\partial B (y, R)$ to $\partial B (y, r)$.
\begin{lem} \label{lem:appendix-excursion}
Fix $n \in \mathbb{N}$, $x, y \in \mathbb{Z}_n^2$, and $1 < r < R/4 < n/8$.
Then, for the SRW on $\mathbb{Z}_n^2$ starting at $x$,
the sequence of excursions
$\left(S_{\cdot \wedge H_{\partial B (y, r)}} 
\circ \theta_{D_{i + 1} (y, R, r)} \right)_{i \geq 0}$
is a $W (y, R, r)$-valued Markov chain
with initial distribution
$$P_x \left[S_{\cdot \wedge H_{\partial B (y, r)}} 
\circ \theta_{D_1 (y, R, r)} \in ~~\cdot ~~\right],$$
and with transition probabilities
\begin{equation} \label{eq:lem-appendix-excursion-transition-prob}
K (\omega, \omega^{\prime})
:= 
P_{\omega (I_{\omega})}
\left[S_{\cdot \wedge H_{\partial B (y, r)}}
\circ \theta_{H_{\partial B (y, R)}} = \omega^{\prime} \right],
\end{equation}
where $I_{\omega}$ is the terminal time of $\omega$
defined by
\begin{equation} \label{eq:terminal-time-appendix}
I_{\omega} = \min \{i \geq 0 : \omega (i) \in \partial B (y, r) \}.
\end{equation}
\end{lem}
The proof of Lemma \ref{lem:appendix-excursion}
is straightforward,
 so we omit the proof.
Fix $n \in \mathbb{N}$, $x, y \in \mathbb{Z}_n^2$, and
$1 < r < R/4 < n/8$.
Assume that the SRW $S$ starts at $x$.
Recall the definition of $q (y, R, r)$ from (\ref{eq:appendix-def-q}).
For each $z \in \partial B (y, r)$,
we define the probability measure $\nu_z$ on $\partial B (y, R)$ by
\begin{equation} \label{eq:appendix-def-measure-nu}
\nu_z (u) := \frac{P_z \left[S_{H_{\partial B (y, R)}} = u \right] - 
q (y, R, r) \mu_r^{y, R} (u)}{1 - q (y, R, r)},
~~~~u \in \partial B (y, R).
\end{equation}
We take $W (y, R, r)$-valued
random variables $X_{\ell}, \ell \geq 0$
and $\{0, 1 \}$-valued random variables $I_{\ell}, \ell \geq 0$ as follows:
\begin{itemize}
\item $X_0$ has the law
$$\mathbb{P} [X_0 = \omega] 
= P_x \left[S_{\cdot \wedge H_{\partial B (y, r)}}
\circ \theta_{D_1 (y, R, r)} = \omega \right],~
\omega \in W (y, R, r).$$
\item $I_0$ is a Bernoulli random variable with success probability $q (y, R, r)$
which is independent of $X_0$. 
\item Suppose that we have constructed 
$X_0, \cdots, X_{\ell}$ and $I_0, \cdots, I_{\ell}$.
Then, take $X_{\ell + 1}$ as follows: 
for each $\omega \in W (y, R, r)$,
\begin{align} \label{eq:appendix-construction-excursion} 
&\mathbb{P} \left[X_{\ell + 1} = \omega \biggm | \sigma (X_j, I_j : 0 \leq j \leq \ell) \right]
\notag \\
&= \begin{cases}
P_{\mu_r^{y, R}} \left[S_{\cdot \wedge H_{\partial B (y, r)}} = \omega \right] 
&\text{if}~I_{\ell} = 1, \\
P_{\nu_{X_{\ell} (I_{X_{\ell}})}} 
\left[S_{\cdot \wedge H_{\partial B (y, r)}} = \omega \right] 
&\text{if}~I_{\ell} = 0.
\end{cases}
\end{align}
Take $I_{\ell + 1}$ as a Bernoulli random variable
with success probability $q (y, R, r)$
which is independent of $X_j, 1 \leq j \leq \ell + 1$ and of $I_j, 1 \leq j \leq \ell$.
\end{itemize}
One can check that
$(X_{\ell})_{\ell \geq 0}$
is a $W (y, R, r)$-valued Markov chain.
By a direct calculation of transition probabilities of $(X_{\ell})_{\ell \geq 0}$
and Lemma \ref{lem:appendix-excursion}, we have
\begin{equation} \label{eq:appendix-excursion-equivalence}
(X_{\ell})_{\ell \geq 0} \stackrel{\text{law}}{=}
\left(S_{\cdot \wedge H_{\partial B (y, r)}}
\circ \theta_{D_{\ell + 1} (y, R, r)} \right)_{\ell \geq 0}.
\end{equation}
We define $J_{\ell}, \ell \geq 0$ inductively as follows:
\begin{equation} \label{eq:appendix-def-j}
J_0 = \min \{i \geq 0 : I_i = 1 \},~
J_{\ell} := \inf \{i > J_{\ell - 1} : I_i = 1 \},~\ell \geq 1.
\end{equation}
Recall the definition of the terminal time
$I_{\omega}$ of $\omega$ from (\ref{eq:terminal-time-appendix}).
Set
\begin{align} \label{eq:appendix-def-sum-hitting-times}
&G_0 := \sum_{0 \leq \ell \leq J_0} I_{X_{\ell}}, \notag \\
&G_m := \sum_{J_{m - 1} < \ell \leq J_m} I_{X_{\ell}},~~m \geq 1.
\end{align} 
It is straightforward to show that $G_m, m \geq 0$ are independent and
$G_m, m \geq 1$ are identically distributed.
We need moment estimates of $G_m$:

\begin{lem} \label{lem:appendix-moment-g}
There exists $c_1 > 0$ such that
for any $n \in \mathbb{N}$,
$x, y \in \mathbb{Z}_n^2$,
$1 < r < R/4 < n/8$,
the following hold: \\
(i) 
\begin{equation} \label{eq:lem-appendix-moment-g-statement1}
\mathbb{E} [G_1] = \frac{E_{\mu_r^{y, R}} \left[H_{\partial B (y, r)} \right]}{q (y, R, r)}.
\end{equation}
(ii) For any $m \in \{0, 1 \}$ and $k \geq 1$,
\begin{equation} \label{eq:lem-appendix-moment-g-statement2}
\mathbb{E} [(G_m)^k] \leq \frac{k !}{\left(q (y, R, r) \right)^k} 
\left\{c_1 n^2 \log \left(\frac{n}{r} \right) \right \}^k.
\end{equation}
\end{lem}
{\it Proof.}
(i) By the definition of $G_1$, we have
\begin{equation} \label{eq:lem-appendix-moment-g-1-pf-1}
\mathbb{E} [G_1]
= \mathbb{E} \left[I_{X_{J_0 + 1}} \right]
+ \sum_{\ell = 2}^{\infty}
\mathbb{E}
\left[1_{\{J_0 + \ell \leq J_1 \}}  I_{X_{J_0 + \ell}} \right].
\end{equation}
By the definition of $(X_{\ell})_{\ell \geq 0}$ and $J_0$,
we have
\begin{equation} \label{eq:lem-appendix-moment-g-1-pf-2}
\mathbb{E} \left[I_{X_{J_0 + 1}} \right]
= E_{\mu_r^{y, R}} \left[H_{\partial B (y, r)} \right].
\end{equation}
By the definitions of $J_0$, $J_1$, and $(X_{\ell})_{\ell \geq 0}$,
for each $\ell \geq 2$,
we have
\begin{equation} \label{eq:lem-appendix-moment-g-1-pf-3}
\mathbb{E} \left[1_{\{J_0 + \ell \leq J_1 \}} I_{X_{J_0 + \ell}} \right]
= E_{\mu_r^{y, R}} \left[
E_{\nu_{X_{\ell - 1} (I_{X_{\ell - 1}})}} \left[H_{\partial B (y, r)} \right]
 \right]
\left(1 - q (y, R, r) \right)^{\ell - 1}.
\end{equation}
By Lemma \ref{lem:appendix-equilibrium-law}
and the definition of $(X_{\ell})_{\ell \geq 0}$,
for any $\ell \geq 2$,
we have
\begin{equation} \label{eq:lem-appendix-moment-g-1-pf-3-2}
E_{\mu_r^{y, R}} \left[
E_{\nu_{X_{\ell - 1} (I_{X_{\ell - 1}})}} \left[H_{\partial B (y, r)} \right]
 \right]
= E_{\mu_r^{y, R}}
\left[H_{\partial B (y, r)} \right].
\end{equation}
By (\ref{eq:lem-appendix-moment-g-1-pf-1})-(\ref{eq:lem-appendix-moment-g-1-pf-3-2}),
we have (\ref{eq:lem-appendix-moment-g-statement1}).

(ii) 
We only prove the $m = 0$ case.
By Kac's moment formula (\ref{eq:lem-appendix-tail-d1-pf-5}),
\begin{align} \label{eq:lem-appendix-moment-g-1-pf-4}
&\mathbb{E} [(G_0)^k]
= \sum_{j = 0}^{\infty}
\mathbb{E}
\left[
1_{\{J_0 = j \}} \left(\sum_{\ell = 0}^j I_{X_{\ell}} \right)^k 
\right] \notag \\
&\leq \sum_{j = 0}^{\infty}
\left(1 - q (y, R, r) \right)^j q (y, R, r) 
\sum_{\begin{subarray}{c} k_0, \cdots, k_j \geq 0, \\ 
k_0 + \cdots + k_j = k \end{subarray}}
\frac{k!}{k_0! \dotsc k_j!}
\prod_{\ell = 0}^j
\max_{v \in \mathbb{Z}_n^2}
E_v \left[\left(H_{\partial B (y, r)} \right)^{k_{\ell}} \right] \notag \\
&\leq \left \{k! 
\sum_{j \geq 0} 
\begin{pmatrix}
k + j \\ k
\end{pmatrix}
\left(1 - q (y, R, r) \right)^j q (y, R, r) \right \}
\left(\max_{v \in \mathbb{Z}_n^2} E_v \left[H_{\partial B (y, r)} \right] \right)^k.
\end{align}
Let $(E_i)_{i \geq 0}$
be independent standard exponential random variables
which are independent of $J_0$.
Then, one can check that
$\mathbb{E} \left[\left(\sum_{\ell = 0}^{J_0} E_{\ell} \right)^k \right]$
is equal to the big bracket on
the right-hand side of (\ref{eq:lem-appendix-moment-g-1-pf-4}).
By this and the fact that
$\sum_{\ell = 0}^{J_0} E_{\ell}$
has the same law as an exponential random variable with mean $q (y, R, r)^{- 1}$,
$\mathbb{E} [(\sum_{\ell = 0}^{J_0} E_{\ell})^k]$
is equal to $\frac{k !}{\left(q (y, R, r) \right)^k}$.
By this and (\ref{eq:lem-appendix-tail-d1-pf-3}),
we have the statement (\ref{eq:lem-appendix-moment-g-statement2})
with $m = 0$.
The proof for the $m = 1$ case is almost the same, so we omit the detail.
~~~$\qed$

\begin{lem} \label{lem:appendix-estimate-sum-g}
There exist $c_1, c_2 \in (0, \infty)$ such that
for any 
$n \in \mathbb{N}$, 
$x, y \in \mathbb{Z}_n^2$,
$1 < r < R/4 < n/8$,
$\delta \in (0, 1)$ and $m \geq 1$,
\begin{equation} \label{eq:lem-appendix-estimate-sum-g-statement}
\mathbb{P} \left[
\left| \frac{\sum_{i = 0}^{m - 1} G_i}{\mathbb{E} [G_1] m} - 1 \right| > \delta \right]
\leq c_1 \exp \left \{- c_2 (m - 1) \delta^2
\left(\frac{\mathbb{E}[G_1]}{\frac{n^2 \log (n/r)}{q (y, R, r)}} \right)^2 \right \}.
\end{equation}
\end{lem}
{\it Proof.}
By the exponential Chebyshev inequality,
for any $\theta > 0$,
we have
\begin{align} \label{eq:lem-appendix-estimate-sum-g-pf-1}
\mathbb{P} \left[
\sum_{i = 0}^{m - 1} G_i
> (1 + \delta) \mathbb{E} [G_1] m
\right] 
&\leq
e^{- \theta (1 + \delta) \mathbb{E} [G_1] m}
\mathbb{E} \left[
\exp \left\{\theta \sum_{i = 0}^{m - 1} G_i \right \} \right] \notag \\
&= e^{- \theta (1 + \delta) \mathbb{E} [G_1] m}
\mathbb{E} \left[e^{\theta G_0} \right]
\left(\mathbb{E} [e^{\theta G_1}] \right)^{m - 1},
\end{align}
where we have used the fact that $G_m,~m \geq 0$ are independent
and that $G_m, m \geq 1$ are identically distributed.
By Lemma \ref{lem:appendix-moment-g},
for any $0 < \theta \leq \frac{1}{2} \frac{q (y, R, r)}{c_1 n^2 \log (n/r)}$,
we have
\begin{align} \label{eq:lem-appendix-estimate-sum-g-pf-2}
\mathbb{E} [e^{\theta G_0}]
&=
1 + \sum_{k = 1}^{\infty} \frac{\theta^k}{k !}
\mathbb{E} [(G_0)^k] \notag \\
&\leq 1 + \sum_{k = 1}^{\infty}
\left(\theta \frac{c_1 n^2 \log (n/r)}{q (y, R, r)} \right)^k
\leq 2.
\end{align}
Similarly,
by Lemma \ref{lem:appendix-moment-g},
for any $0 < \theta \leq \frac{1}{2} \frac{q (y, R, r)}{c_1 n^2 \log (n/r)}$,
we have
\begin{align} \label{eq:lem-appendix-estimate-sum-g-pf-3}
\mathbb{E} [e^{\theta G_1}]
&=
1 + 
\theta \mathbb{E} [G_1]
+ \sum_{k = 2}^{\infty} \frac{\theta^k}{k !}
\mathbb{E} [(G_1)^k] \notag \\
&\leq 1 + 
\theta \mathbb{E} [G_1]
+ \sum_{k = 2}^{\infty}
\left(\theta \frac{c_1 n^2 \log (n/r)}{q (y, R, r)} \right)^k \notag \\
&\leq \exp \left\{\theta \mathbb{E} [G_1]
+ 2 \theta^2
\left(\frac{c_1 n^2 \log (n/r)}{q (y, R, r)} \right)^2 \right \}.
\end{align}
By (\ref{eq:lem-appendix-estimate-sum-g-pf-2})
and (\ref{eq:lem-appendix-estimate-sum-g-pf-3}),
for the optimal value $\theta = \frac{\delta}{4} 
\mathbb{E} [G_1] \left(\frac{q (y, R, r)}{c_1 n^2 \log (n/r)} \right)^2$,
the right of (\ref{eq:lem-appendix-estimate-sum-g-pf-1})
is bounded from above by a value of the form of
the right of (\ref{eq:lem-appendix-estimate-sum-g-statement}).
Repeating a similar argument, we can obtain 
a similar estimate of the left tail.
We omit the detail.~~~$\qed$

\begin{lem} \label{lem:appendix-sum-explore-times}
There exist $c_1, c_2 \in (0, \infty)$ such that
for any $n \in \mathbb{N}$,
$x, y \in \mathbb{Z}_n^2$, 
$1 < r < R/4 < n/8$,
$m > 1 + 4 \left(q (y, R, r) \right)^{- 1}$, 
$\delta \in (\frac{4}{(m - 1) q (y, R, r)}, 1)$,
\begin{align} \label{eq:lem-appendix-sum-explore-times-statement}
&P_x \left[
\left| \frac{\sum_{i = 1}^{m - 1} \left(R_{i + 1} (y, R, r) - D_i (y, R, r) \right)}
{E_{\mu_r^{y, R}} [H_{\partial B (y, r)}] (m - 1)} - 1 \right| > \delta \right] \notag \\
&\leq c_1 \exp \left \{- c_2 (m - 1) q (y, R, r) \delta^2
\left(\frac{E_{\mu_r^{y, R}} [H_{\partial B (y, r)}]}{n^2 \log (n/r)} \right)^2 \right \}.
\end{align}
\end{lem}
{\it Proof.}
Set
$m^+ := \left \lceil \frac{(m - 1) q (y, R, r)}{1 - \delta/8} \right \rceil$.
By (\ref{eq:appendix-excursion-equivalence}),
we have
\begin{align} \label{eq:lem-appendix-sum-explore-times-pf-1}
&P_x \left[
\sum_{i = 1}^{m - 1}
\left(R_{i + 1} (y, R, r) - D_i (y, R, r) \right)
> (1 + \delta) E_{\mu_r^{y, R}} \left[H_{\partial B (y, r)} \right] (m - 1) \right] \notag \\
&= \mathbb{P}
\left[
\sum_{\ell = 0}^{m - 2}
I_{X_{\ell}}
> (1 + \delta) E_{\mu_r^{y, R}}
\left[H_{\partial B (y, r)} \right] (m - 1) \right] \notag \\
&\leq \mathbb{P} \left[
\sum_{\ell = 0}^{m^+ - 1} G_{\ell}
> (1 + \delta) E_{\mu_r^{y, R}}
\left[H_{\partial B (y, r)} \right] (m - 1) \right]
+ \mathbb{P} \left[J_{m^+ - 1} < m - 2 \right].
\end{align}
By Lemma \ref{lem:appendix-moment-g} (i),
the definition of $m^+$,
and the condition
$\delta \in (\frac{4}{(m - 1) q (y, R, r)}, 1)$,
we have
$$(1 + \delta) E_{\mu_r^{y, R}} \left[H_{\partial B (y, r)} \right] (m - 1)
\geq \left(1 + \frac{\delta}{4} \right) \mathbb{E} [G_1] m^+.$$
By this and Lemmas \ref{lem:appendix-moment-g}, \ref{lem:appendix-estimate-sum-g},
the first term on the right of (\ref{eq:lem-appendix-sum-explore-times-pf-1})
is bounded from above by
a value of the form of the right of (\ref{eq:lem-appendix-sum-explore-times-statement}).
By the definition of $J_{m^+ - 1}$
and the exponential Chebyshev inequality,
for any $\theta \in (0, 1)$,
the second term on the right of (\ref{eq:lem-appendix-sum-explore-times-pf-1})
is bounded from above by
\begin{align} \label{eq:lem-appendix-sum-explore-times-pf-2}
&\mathbb{P} \left[I_0 + \cdots + I_{m - 2} \geq m^+ \right]
\leq e^{- \theta m^+}
\mathbb{E} \left[e^{\theta \sum_{\ell = 0}^{m - 2} I_{\ell}} \right].
\end{align}
Since $I_{\ell},~\ell \geq 0$
are i.i.d. Bernoulli random variables with parameter $q (y, R, r)$,
the expectation on the right of (\ref{eq:lem-appendix-sum-explore-times-pf-2})
is bounded from above by
\begin{equation} \label{eq:lem-appendix-sum-explore-times-pf-3}
\left\{e^{\theta} q (y, R, r) + 1 - q (y, R, r) \right\}^{m - 1} 
\leq \exp \left\{(\theta + 2 \theta^2) q (y, R, r) (m - 1) \right \},
\end{equation}
where we have used the inequality
$e^{\theta} \leq 1 + \theta + 2 \theta^2$
for any $\theta \in (0, 1)$
in the last inequality.
Optimizing $\theta$
(take $\theta = \frac{1}{4} \left(\frac{m^+}{(m - 1) q (y, R, r)} - 1 \right)$),
the right of (\ref{eq:lem-appendix-sum-explore-times-pf-2})
is bounded from above by
\begin{equation} \label{eq:lem-appendix-sum-explore-times-pf-4}
\exp \left\{- c_3 (m - 1) q (y, R, r) \delta^2 \right \}.
\end{equation}
This is bounded from above by a value 
of the form of the right of (\ref{eq:lem-appendix-sum-explore-times-statement})
since $E_{\mu_r^{y, R}} [H_{\partial B (y, r)}]$ $\leq c n^2 \log (n/r)$
by (\ref{eq:lem-appendix-tail-d1-pf-3}).
Therefore, the right of (\ref{eq:lem-appendix-sum-explore-times-pf-1})
is bounded from above by a value of the form of 
the right of (\ref{eq:lem-appendix-sum-explore-times-statement}).
Repeating a similar argument,
we can obtain a similar estimate of the left tail.
We omit the detail.~~~$\qed$

\begin{prop} \label{prop:appendix-estimate-time-d}
There exist $c_i \in (0, \infty),~1 \leq i \leq 8$ such that
for any $n \in \mathbb{N}$, 
$x, y \in \mathbb{Z}_n^2$,
$c_1 \leq r < R/4 < n/8$,
$m \geq 1 + 4 \left(q (y, R, r) \right)^{- 1}$, 
$\delta \in (\frac{4}{(m - 1) q (y, R, r)}, 1)$,
\begin{align} \label{eq:prop-appendix-estimate-time-d-statement}
&P_x \left[
\left| \frac{D_m (y, R, r)}{E_{\mu_r^{y, R}} [D_1 (y, R, r)] (m - 1)} - 1 \right| \geq \delta \right] 
\notag \\
&\leq c_2 \exp \left \{- c_3 \delta 
\frac{E_{\mu_r^{y, R}} \left[D_1 (y, R, r) \right]}{n^2 \log (n/r)} (m - 1) \right \} \notag \\
&+ c_4 \exp \left \{- c_5 \delta^2 (m - 1)
+ c_6 \delta \frac{1}{R} (m - 1) \right \} \notag \\
&+ c_7 \exp \left\{- c_8 q (y, R, r) \delta^2 
\left(\frac{E_{\mu_r^{y, R} \left[H_{\partial B (y, r)} \right]}}{n^2 \log (n/r)} \right)^2
(m - 1) \right \}.
\end{align}
\end{prop}
{\it Proof.}
By Lemma \ref{lem:appendix-equilibrium-law},
we have
\begin{equation} \label{eq:prop-appendix-estimate-time-d-pf-1}
E_{\mu_r^{y, R}} \left[D_1 (y, R, r) \right]
= E_{\mu_r^{y, R}} \left[H_{\partial B (y, r)} \right]
+ E_{\mu_R^{y, r}} \left[H_{\partial B (y, R)} \right].
\end{equation}
Note that
$D_m (y, R, r)$ is equal to
\begin{equation} \label{eq:prop-appendix-estimate-time-d-pf-2}
D_1 (y, R, r)
+ \sum_{i = 2}^m
\left(D_i (y, R, r) - R_i (y, R, r) \right)
+ \sum_{i = 1}^{m - 1}
\left(R_{i + 1} (y, R, r) - D_i (y, R, r) \right).
\end{equation}
By (\ref{eq:prop-appendix-estimate-time-d-pf-1})
and (\ref{eq:prop-appendix-estimate-time-d-pf-2}),
the left of (\ref{eq:prop-appendix-estimate-time-d-statement})
is bounded from above by
\begin{align} \label{eq:prop-appendix-estimate-time-d-pf-3}
&P_x \left[
D_1 (y, R, r) \geq \frac{\delta}{2}
E_{\mu_r^{y, R}} \left[D_1 (y, R, r) \right] (m - 1) \right] \notag \\
&+ P_x \left[
\left|\frac{\sum_{i = 2}^m \left(D_i (y, R, r) - R_i (y, R, r) \right)}
{E_{\mu_R^{y, r}} \left[H_{\partial B (y, R)} \right] (m - 1)}
- 1 \right| \geq \frac{\delta}{2} \right] \notag \\
&+ P_x \left[
\left|\frac{\sum_{i = 1}^{m - 1} \left(R_{i + 1} (y, R, r) - D_i (y, R, r) \right)}
{E_{\mu_r^{y, R}} \left[H_{\partial B (y, r)} \right] (m - 1)}
- 1 \right| \geq \frac{\delta}{2} \right].
\end{align}
By (\ref{eq:prop-appendix-estimate-time-d-pf-3})
and Lemmas \ref{lem:appendix-tail-d1},
\ref{lem:appendix-sum-exit-time},
\ref{lem:appendix-sum-explore-times},
we have (\ref{eq:prop-appendix-estimate-time-d-statement}).~~~$\qed$
\\\\

{\it Proof of Proposition \ref{prop:time-number}.}
We only prove (\ref{eq:upper-time}).
Set $r_1^{-} := \left(1 - \frac{\sqrt{2}}{(\log n)^2} \right) r_1$,
$r_0^{+} := \left(1 + \frac{\sqrt{2}}{(\log n)^2} \right) r_0$.
We define the set $F$ by
\begin{equation*}
F := \left \{\left(i \left \lfloor \frac{r_0}{\ell_n (\log n)^2} \right \rfloor, 
~j \left \lfloor \frac{r_0}{\ell_n (\log n)^2} \right \rfloor \right)
: 0 \leq i, j \leq \left \lfloor \frac{n}{\left \lfloor \frac{r_0}{\ell_n (\log n)^2} \right \rfloor} 
\right \rfloor
\right \}.
\end{equation*}
By a simple observation, one can show that 
for any $x \in \mathbb{Z}_n^2$,
there exists $y \in F$ such that 
\begin{equation*}
x \in \left(y + \left[0,~\left \lfloor \frac{r_0}{\ell_n (\log n)^2} \right \rfloor \right]^2 \right)
\cap \mathbb{Z}^2~~~\text{mod}~n \mathbb{Z}^2
\end{equation*}
and that
\begin{equation*}
B (y, r_1^{-}) \subset B (x, r_1) \subset B (x, r_0) \subset B (y, r_0^{+}).
\end{equation*}
In particular, we have
\begin{equation*}
D_{m_n^{+}} (x, r_0, r_1) \leq D_{m_n^{+}} (y, r_0^{+}, r_1^{-}).
\end{equation*}
By this observation, we have
\begin{align} \label{eq:upper-time-pf-step1}
&P_o \left[
\begin{minipage}{220pt}
$\exists x \in \mathbb{Z}_n^2$, \\
$D_{m_n^{+}} (x, r_0, r_1) > \frac{4}{\pi} n^2 (\log n)^2
\left(1 - \frac{\log \log n}{2 \log n} + \frac{2 s_n}{\log n} \right)$
\end{minipage}
\right] \notag \\
&\leq \sum_{y \in F}
P_o \left[
D_{m_n^{+}} (y, r_0^{+}, r_1^{-}) > \frac{4}{\pi} n^2 (\log n)^2
\left(1 - \frac{\log \log n}{2 \log n} + \frac{2 s_n}{\log n} \right)
\right].
\end{align}
By Lemma \ref{lem:appendix-expectation-d1},
for sufficiently large $n$,
we have
\begin{align} \label{eq:upper-time-pf-step2}
&E_{\mu_{r_1^-}^{y, r_0^+}} \left[D_1 (y, r_0^+, r_1^-) \right]
(m_n^+ - 1) \left(1 + \frac{s_n}{4 \log n} \right) \notag \\
&\leq \frac{4}{\pi} n^2 (\log n)^2
\left(1 - \frac{\log \log n}{2 \log n} + \frac{2 s_n}{\log n} \right).
\end{align}
Note that
by Lemma \ref{lem:appendix-expectation-d1},
(\ref{eq:lem-appendix-tail-d1-pf-4}),
and (\ref{eq:prop-appendix-estimate-time-d-pf-1}),
we have
$$E_{\mu_{r_1^-}^{y, r_0^+}} \left[H_{\partial B (y, r_1^-)} \right]
\geq c_1 n^2 \log (\ell_n).$$
By this, (\ref{eq:upper-time-pf-step2}),
Lemmas \ref{lem:appendix-expectation-d1},
\ref{lem:appendix-estimate-q},
and Proposition \ref{prop:appendix-estimate-time-d},
each term on the right of (\ref{eq:upper-time-pf-step1})
is bounded from above by
\begin{equation} \label{eq:upper-time-pf-step3}
c_2 e^{- c_3 (\log \log n)^{2 \gamma - \alpha - 2 \beta}}.
\end{equation}
By the definition of the set $F$, we have
\begin{equation} \label{eq:upper-time-pf-step4}
|F| \leq c_4 (\log n)^4 (\ell_n)^4 e^{2 c_{\star} (\log \log n)^{\alpha + \beta}}.
\end{equation}
Note that by the assumption (\ref{eq:assumption-parameter}),
we have
$2 \gamma - \alpha - 2 \beta > 1 > \alpha + \beta$.
By this, (\ref{eq:upper-time-pf-step1}), (\ref{eq:upper-time-pf-step3}),
and (\ref{eq:upper-time-pf-step4}),
we have (\ref{eq:upper-time}).
The proof of (\ref{eq:lower-time}) is almost the same as that of (\ref{eq:upper-time}).
So, we omit the detail.
$\qed$

\section{Proofs of Lemma \ref{lem:transfer-lemma} 
and (\ref{eq:lem-2nd-moment-case3-pf-pc-last})} \label{sec:pf-transfer-lem}
In this section, we prove Lemma \ref{lem:transfer-lemma}. 
For simplicity, we drop $\overline{R}$ from the notation.
\\\\
{\it Proof of Lemma \ref{lem:transfer-lemma}(i).}
We will only prove the upper bound in (\ref{eq:transfer-lemma-statement-1}).
We will keep track of jumps among
$\partial B (x, R_{\ell})$,~$\ell \in I_{k, \widetilde{k}}^L$,
where
$I_{k, \widetilde{k}}^L := \{0, 1, k, k + 1, \dotsc, L - \widetilde{k}, L \}$.
Note that 
there are two types of
excursions from $\partial B (x, R_1)$ to $\partial B (x, R_0)$:
\begin{itemize}
\item[(a)] excursions which hit $\partial B (x, R_{k + 1})$ ~(we will call them of type (a)),
\item[(b)] excursions which do not hit $\partial B (x, R_{k + 1})$ ~(we will call them of type (b)).
\end{itemize}
Let $\mathbf{e}$ be an excursion 
from $\partial B (x, R_1)$ to $\partial B (x, R_0)$.
When $\mathbf{e}$ is of type (a),
we will define the process $\text{tr}^{\mathbf{e}} (i)$, $i \geq 0$,
which records subscripts of radii of the circles SRW visits
as follows:
Set
\begin{equation*}
\sigma_0 := 0,~~\text{tr}^{\mathbf{e}} (0) := 1,
\end{equation*}
\begin{equation*}
\sigma_1 := \min \{i : \mathbf{e} (i) \in \partial B (x, R_{k + 1}) \},
~~\text{tr}^{\mathbf{e}} (1) := k + 1.
\end{equation*}
Suppose that we have constructed 
$\sigma_{\ell}$, $\text{tr}^{\mathbf{e}} (\ell)$, $0 \leq \ell \leq i - 1$ ($i \geq 2$)
and that $\text{tr}^{\mathbf{e}} (\ell) \neq 0$ for all $0 \leq \ell \leq i - 1$.
Then, we define $\sigma_i$ and $\text{tr}^{\mathbf{e}} (i)$ by
\begin{equation*}
\sigma_i := \min \left\{j > \sigma_{i - 1} : \mathbf{e} (j)
\in \bigcup_{\ell \in I_{k, \widetilde{k}}^L \backslash \{1 \},~
\ell \neq \text{tr}^{\mathbf{e}} (i - 1)}
\partial B (x, R_{\ell}) \right \},
\end{equation*}
\begin{equation*}
\text{tr}^{\mathbf{e}} (i) := \ell~\text{such that}~\mathbf{e} (\sigma_i) \in \partial B (x, R_{\ell}).
\end{equation*}
Stop the construction at the first time when $\text{tr}^{\mathbf{e}}$ reaches $0$.
When $\mathbf{e}$ is of type (b),
we define $\text{tr}^{\mathbf{e}} = (\text{tr}^{\mathbf{e}} (i))_{i \in \{0, 1 \}}$ 
just by $\text{tr}^{\mathbf{e}} (0) := 1$, $\text{tr}^{\mathbf{e}} (1) = 0$.

Let $\mathbf{e}_i$ ($1 \leq i \leq m$) be the $i$-th excursion
from $\partial B (x, R_1)$ to $\partial B (x, R_0)$.
We define the process $\text{tr}^S$ by concatenating
$\text{tr}^{\mathbf{e}_1}, \dotsc, \text{tr}^{\mathbf{e}_m}$
in this order.
Let $W$ be the state space of $\text{tr}^S$.
For each $w \in W$ and
$i, j \in I_{k, \widetilde{k}}^L$ with $i < j$, let
$T_{i \to j}^w$
\footnote{More precisely, $T_{i \to j}^w$ is defined as follows:
set 
$R_1^w := \min \{\ell : w (\ell) = j \}$,
$D_p^w := \min \{\ell > R_p^w : w (\ell) = i \}$,
$R_{p+1}^w := \min \{\ell > D_p^w : w(\ell) = j \}$, $p \geq 1$.
Then, we define $T_{i \to j}^w$ by
$T_{i \to j}^w := \max \{\ell : R_{\ell}^w < t_w \}$,
where $t_w$ is the terminal time of $w$.}
be the number of traversals from $i$ to $j$ by $w$.
Let $\widetilde{W}$ be the set of all $w \in W$ which satisfies
$T_{i \to i + 1}^w = m_i$, $k \leq \forall i \leq L - \widetilde{k} - 1$
and $T_{L - \widetilde{k} \to L}^w = 0$.
(Recall that $(m_i)_{i = k}^{L - \widetilde{k} - 1}$ are integers on the left-hand side
of (\ref{eq:transfer-lemma-statement-1}).)

Fix any $w \in \widetilde{W}$.
Recall the definition of $p_{i_1, i_2}^{i_3, +}$ from (\ref{eq:prob-circle-to-circle-1})
and (\ref{eq:prob-circle-to-circle-2}).
By the strong Markov property and Lemma \ref{lem:circle-to-circle-probab}, we have
\begin{align} \label{eq:pf-transfer-lemma-1-1}
&P_y [\text{tr}^S = w] \notag \\
&\leq \left(p_{1, 0}^{k + 1, +} \right)^{m - T_{1 \to k + 1}^w}
\left(p_{1, k + 1}^{0, +} \right)^{T_{1 \to k + 1}^w} \notag \\
&\times \left(p_{k, 0}^{k + 1, +} \right)^{T_{1 \to k + 1}^w}
\left(p_{k, k + 1}^{0, +} \right)^{m_k - T_{1 \to k + 1}^w} \notag \\
&\times \prod_{i = k + 1}^{L - \widetilde{k} - 1}
\left\{\left(p_{i, i - 1}^{i + 1, +} \right)^{m_{i - 1}} \left(p_{i, i + 1}^{i - 1, +} \right)^{m_i} \right \} \notag \\
&\times \left(p_{L - \widetilde{k}, L - \widetilde{k} - 1}^{L, +} \right)^{m_{L - \widetilde{k} - 1}}.
\end{align}
Recall the definitions of $\Delta_{i_2, i_3}^{i_1, +}$, $\Delta_1^+$,
and $\Delta_{\star}^+$
from (\ref{eq:error}) and (\ref{eq:transfer-lemma-main-error-1}).
By these definitions and the strong Markov property of the SRW 
$Z = (Z_i, i \geq 0, P_k^{1D}, k \in \{0, \dotsc, L\})$ 
on $\{0, \dotsc, L \}$,
the right-hand side of (\ref{eq:pf-transfer-lemma-1-1}) is bounded from above by
\begin{equation} \label{eq:pf-transfer-lemma-1-2}
\Delta_1^+ \cdot \Delta_{\star}^+ P_1^{1D} [\text{tr}^Z = w],
\end{equation}
where $\text{tr}^Z$ is 
the process which records vertices SRW visits
and is defined in the same manner as $\text{tr}^S$.  
By (\ref{eq:pf-transfer-lemma-1-1}) and (\ref{eq:pf-transfer-lemma-1-2}),
summing over $w \in \widetilde{W}$,
we can bound the numerator on the left-hand side of (\ref{eq:transfer-lemma-statement-1})
from above by
\begin{equation} \label{eq:pf-transfer-lemma-1-3}
\Delta_1^+ \Delta_{\star}^+
P_1^{1D} \left[T_i^m = m_i,~k \leq \forall i \leq L - \widetilde{k} - 1,~
T_{L - 1}^m = 0 \right],
\end{equation}
where $T_i^m$
is the number of traversals from $i$ to $i + 1$
by $Z$ starting at $1$ up to the $m$-th return to $0$
\footnote{\label{fn:1d-traversal-process} More precisely, $T_i^m$ is defined as follows:
Set $R_1^i := \min \{j : Z_j = i + 1 \}$,
$D_p^i := \min \{j > R_p^i : Z_j = i \}$,
$R_{p+1}^i := \min \{j > D_p^i : Z_j = i + 1 \}$, $p \geq 1$.
Then, we define $T_i^m$ by $T_i^m := \max \{j : R_j^i < D_m^0 \}$.}.
Since $P_1^{1D}$-law of $(T_i^m)_{i \geq 0}$ 
is the same as $P_m^{\text{GW}}$-law of $(T_i)_{i \geq 0}$,
we have obtained the upper bound.
Other statements of (i) can be proved similarly.
We omit the details. \\\\
{\it Proof of Lemma \ref{lem:transfer-lemma} (ii).}
We have
\begin{equation} \label{eq:pf-transfer-lemma-2-1}
P_y \left[T_i^{x, m} \geq \ell, T_{L - 1}^{x, m} = 0 \right]
= \sum_{k = 0}^{m - 1} P_y \left[T_i^{x, k} < \ell \leq T_i^{x, k + 1}, T_{L - 1}^{x, m} = 0 \right].
\end{equation}
Fix $0 \leq k \leq m - 1$.
By the strong Markov property, 
the $k$-th term on the right-hand side of (\ref{eq:pf-transfer-lemma-2-1})
is bounded from above by
\begin{align} \label{eq:pf-transfer-lemma-2-2}
&\sum_{j = 0}^{\ell - 1}
P_y \left[T_i^{x, k} = j, T_i^{x, k + 1} \geq \ell, ~T_{L - 1}^{x, m} = 0 \right] \notag \\
&\leq \sum_{j = 0}^{\ell - 1}
P_y \left[T_i^{x, k} = j, ~T_{L - 1}^{x, k} = 0 \right]
\cdot \max_{z \in \partial B (x, R_0)} P_z
\left[T_i^{x, 1} \geq \ell - j, ~T_{L - 1}^{x, m - k} = 0 \right].
\end{align}
By Lemma \ref{lem:transfer-lemma} (i), for each $0 \leq j \leq \ell - 1$, 
the first factor of the $j$-th term on the right of (\ref{eq:pf-transfer-lemma-2-2})
is bounded from above by 
\begin{equation} \label{eq:pf-transfer-lemma-2-3}
\left(\Delta_{1, 0}^{i + 1, +} \vee \Delta_{1, i + 1}^{0, +} \right)^k
\cdot \left(\Delta_{i, 0}^{i + 1, +} \vee \Delta_{i, i + 1}^{0, +} \right)^j
\cdot \left(\Delta_{i + 1, i}^{L, +} \right)^j
\cdot P_1^{1D} \left[T_i^k = j,~T_{L - 1}^k = 0 \right].
\end{equation}
By Lemma \ref{lem:circle-to-circle-probab},
the second factor of the $j$-th term on the right of (\ref{eq:pf-transfer-lemma-2-2})
is bounded from above by
\begin{equation} \label{eq:pf-transfer-lemma-2-4}
p_{1, i + 1}^{0, +} \left(p_{i, i + 1}^{0, +} \right)^{\ell - j - 1}
\left(p_{i + 1, i}^{L, +} \right)^{\ell - j - 1} p_{i + 1, 0}^{L, +}
\left(p_{1, 0}^{L, +} \right)^{m - k - 1}.
\end{equation}
By the definition of $\Delta_{i_2, i_3}^{i_1, +}$ (see (\ref{eq:error}))
and the strong Markov property of SRW on $\{0, \dotsc, L \}$,
(\ref{eq:pf-transfer-lemma-2-4}) is bounded from above by
\begin{align} \label{eq:pf-transfer-lemma-2-5}
&\Delta_{1, i + 1}^{0, +}
\left(\Delta_{i, i + 1}^{0, +} \right)^{\ell - j - 1}
\left(\Delta_{i + 1, i}^{L, +} \right)^{\ell - j - 1}
\Delta_{i + 1, 0}^{L, +}
\left(\Delta_{1, 0}^{L, +} \right)^{m - k - 1} \notag \\
&\times P_1^{1D} \left[T_i^1 \geq \ell - j,~T_{L - 1}^{m - k} = 0 \right].
\end{align}
By (\ref{eq:pf-transfer-lemma-2-1})-(\ref{eq:pf-transfer-lemma-2-5})
and the strong Markov property,
(\ref{eq:pf-transfer-lemma-2-1})
is bounded from above by
\begin{align} \label{eq:pf-transfer-lemma-2-6}
\Delta_2^+ P_1^{1D} \left[T_i^m \geq \ell,~T_{L - 1}^m = 0 \right].
\end{align}
Since $P_1^{1D}$-law of $(T_j^m)_{j \geq 0}$
is the same as $P_m^{\text{GW}}$-law of $(T_j)_{j \geq 0}$,
we have obtained the desired result.
$\qed$ \\\\
{\it Proof of Lemma \ref{lem:transfer-lemma} (iii).}
By the strong Markov property,
the probability on the left-hand side of (\ref{eq:transfer-lemma-statement-3})
is bounded from above by
\begin{align} \label{eq:pf-transfer-lemma-3-1}
&\sum \prod_{\ell = 1}^m \max_{z \in \partial B (x, R_{k+1})}
P_z \left[
T_i^{k, x, 1} = m_i^{(\ell)},~
k + 1 \leq \forall i \leq L - \widetilde{k} - 1,~
T_{L - 1}^{k, x, 1} = 0 \right] \notag \\
&\times \left(\max_{z \in \partial B (x, R_k)} 
P_z \left[H_{\partial B (x, R_{k + 1})} < H_{\partial B (x, R_0)} \right] \right)^m \notag \\
&\times \max_{z \in \partial B (x, R_k)} 
P_z \left[H_{\partial B (x, R_0)} < H_{\partial B (x, R_{k+1})} \right],
\end{align}
where the sum is taken over all nonnegative integers 
$m_i^{(\ell)}$, $k + 1 \leq i \leq L - \widetilde{k} - 1$, $1 \leq \ell \leq m$
satisfying $\sum_{\ell = 1}^m m_i^{(\ell)} = m_i$ for each $k + 1 \leq i \leq L - \widetilde{k} - 1$.
By Lemmas \ref{lem:transfer-lemma} (i) and \ref{lem:circle-to-circle-probab},
(\ref{eq:pf-transfer-lemma-3-1}) is bounded from above by
\begin{align} \label{eq:pf-transfer-lemma-3-2}
&\sum \prod_{i = k + 1}^{L - \widetilde{k} - 1}
\left\{\left(\Delta_{i, i - 1}^{i + 1, +} \vee \Delta_{i, i + 1}^{i - 1, +} \right)^{m_{i - 1} + m_i} \right\}
\left(\Delta_{L - \widetilde{k}, L - \widetilde{k} - 1}^{L, +} \right)^{m_{L - \widetilde{k} - 1}} \notag \\
&\times \prod_{\ell = 1}^m P_1^{\text{GW}}
\left[T_i = m_{k+i}^{(\ell)},~1 \leq \forall i \leq L - k - \widetilde{k} - 1,~T_{L - k - 1} = 0 \right] \notag \\
&\times \left(\Delta_{k, k + 1}^{0, +} \right)^m \Delta_{k, 0}^{k + 1, +}
\left(\frac{k}{k+1} \right)^m \frac{1}{k+1}.
\end{align}
(\ref{eq:pf-transfer-lemma-3-2}) immediately yields the desired result.
$\qed$ \\\\
{\it Proof of (\ref{eq:lem-2nd-moment-case3-pf-pc-last}).}
Recall the definitions $P_c$, $P_1$, $C_{m^{\prime}}^1$
from (\ref{eq:lem-2nd-moment-case3-pf-pc}), (\ref{eq:lem-2nd-moment-case3-pf-12}),
and above (\ref{eq:lem-2nd-moment-case3-pf-3}).
By Lemma \ref{lem:circle-to-circle-probab}, $P_c$ is bounded from above by
$(1 + o(1))$ times
\begin{align} \label{eq:pf-pc-last-1}
\left\{\prod_{\ell = 1}^p
\left(\frac{1}{k - 2} \right) \left(\frac{k - 3}{k - 2} \right)^{j_{\ell}^c - 1}
\left(\frac{1}{k - 2} \right) \right\} 
\left\{\left(\frac{1}{k - 2} \right) 
\left(\frac{k - 3}{k - 2} \right)^{m_n^- - m^{\prime} - \sum_{\ell = 1}^p j_{\ell}^c -1} \right\}.
\end{align}
We will relate (\ref{eq:pf-pc-last-1})
to a law of the $1$-dimensional SRW.
Recall the definitions of $D_{\ell}^k$, $R_{\ell}^k$ 
from the footnote \ref{fn:1d-traversal-process}.
For each $m \in \mathbb{N}$, set
$$\widetilde{T}^m := \max \left\{\ell \geq 1 : D_{\ell}^0 < R_m^{k - 3} \right\}.$$
By the strong Markov property,
the right of (\ref{eq:pf-pc-last-1}) is equal to
\begin{equation} \label{eq:pf-pc-last-2}
(1 + o(1))
P_{k - 2}^{1D}
\left[\bigcap_{\ell = 1}^p 
\left\{\widetilde{T}^1 \circ \theta_{D_{\ell}^{k - 3}} = j_{\ell}^c \right\}
\cap \left\{\widetilde{T}^1 \circ \theta_{D_{p+1}^{k-3}} \geq
m_n^- - m^{\prime} - \sum_{\ell = 1}^p j_{\ell}^c \right\} \right].
\end{equation}
By (\ref{eq:pf-pc-last-1}) and (\ref{eq:pf-pc-last-2}), we have
\begin{equation} \label{eq:pf-pc-last-3}
\sum_{(j_{\ell}^c)} P_c \leq (1 + o(1)) P_{k - 3}^{1D} 
\left[\widetilde{T}^p < m_n^- - m^{\prime} \leq \widetilde{T}^{p+1} \right],
\end{equation}
where the sum is taken over $(j_{\ell}^c)_{\ell = 1}^p$
satisfying the condition in (\ref{eq:condition-traversal-c}).
Similarly, by Lemma \ref{lem:circle-to-circle-probab}
and the strong Markov property, we have
\begin{equation} \label{eq:pf-pc-last-4}
P_o [C_{m^{\prime}}^1]
\leq (1 + o(1)) P_1^{1D} \left[D_{m^{\prime}}^0 < H_{k-2} < D_{m^{\prime}+1}^0 \right].
\end{equation}
By (\ref{eq:pf-pc-last-3}), (\ref{eq:pf-pc-last-4}), and the strong Markov property, 
$P_o [C_{m^{\prime}}^1] \cdot \sum_{(j_{\ell}^c)} P_c$
is bounded from above by $(1 + o(1))
P_1^{1D} [D_{m^{\prime}}^0 < H_{k-2} < D_{m^{\prime} + 1}^0,~
\widetilde{T}^p \circ \theta_{D_1^{k-3}}
 < m_n^- - m^{\prime} \leq \widetilde{T}^{p+1} \circ \theta_{D_1^{k-3}}]$.
This probability is equal to
$P_1^{1D} [
D_{m^{\prime}}^0 < H_{k-2} < D_{m^{\prime} + 1}^0,~
T_{k - 3}^{m_n^-} = p + 1]$.
Therefore, summing over 
$p = b_n^- (k-3) - 1, \dotsc, b_n^+ (k-3) - 1$ and
$m^{\prime} = 0, \dotsc, m_n^- - 1$
and using the fact that $P_1^{1D}$-law of $T_{k-3}^{m_n^-}$
is the same as $P_{m_n^-}^{GW}$-law of $T_{k-3}$,
we have the desired result. $\qed$

\section{Barrier estimates} \label{sec:barrier-estimates}
In this section, we will prove Lemma \ref{lem:GW-process-barrier-estimate}.
The proof is heavily based on arguments in \cite{BRZ1}
and we use the same notation as in \cite{BRZ1}.

As mentioned in Remark \ref{rem:barrier-estimate},
we need to know how the constant $c$ in \cite[Theorem 1.1]{BRZ1}
depends on $\eta$. 
The proof of Lemma \ref{lem:GW-process-barrier-estimate} (i)
is almost the same as that of \cite[Theorem 1.1]{BRZ1},
so we only show how we should slightly modify the argument.
The constant $c$ in \cite[Lemma 2.5]{BRZ1} depends on $\eta$ as follows:
\begin{lem} \label{lem:barrier-estimate-0-dim-bessel-process}
Under the same assumption as in \cite[Lemma 2.5]{BRZ1},
\begin{align} \label{eq:barrier-estimate-0-dim-bessel-process}
&P_x^Y \left[f_{a, b} (\ell; L) - C \ell_L^{\frac{1}{2} - \varepsilon}
\leq Y_{\ell},~\ell = 1, \dotsc, L - 1,~Y_L \in H_{y, \delta} \right] \notag \\
&\leq c e^{\delta \eta} 
\frac{(1 + x - a)(1 + y - b)}{L}
\sqrt{\frac{x}{yL}}
e^{- \frac{(y - x)^2}{2L}},
\end{align}
where $c > 0$ depends only on $\delta, \varepsilon, C$ (not on $\eta$).
\end{lem}
{\it Proof.}
As mentioned in \cite[Remark 2.2]{BRZ1},
the second line in \cite[(2.3)]{BRZ1}
is bounded from above by \cite[(2.4)]{BRZ1}
multiplied by a constant depending only on $\delta, \varepsilon, C$.
By the assumption $|x - y| \leq \eta L$, we have
$(x - z)^2 \geq (x - y)^2 - 2 \eta L \delta$ for any $z \in H_{y, \delta}$.
Thus, the exponential factor in \cite[(2.4)]{BRZ1}
is bounded from above by
$e^{\delta \eta} e^{- \frac{(x - y)^2}{2L}}$.
By this and the proofs of \cite[Lemmas 2.3, 2.5]{BRZ1},
we have the desired result.
$\qed$ \\\\
{\it Proof of Lemma \ref{lem:GW-process-barrier-estimate} (i).}
As mentioned before,
the proof is almost the same as that of \cite[Theorem 1.1]{BRZ1}, but
we need some minor modifications as follows:
\begin{itemize}
\item Below \cite[(4.4)]{BRZ1}:
Using our additional condition $a \geq b$,
we have
\begin{equation*}
\sqrt{T_0 + T_1} \geq \frac{\sqrt{T_0} + \sqrt{T_1}}{\sqrt{2}}
\geq a^{\prime} + \frac{x - a}{2} - C.
\end{equation*}

\item \cite[(4.7)]{BRZ1}:
When $\frac{\ell}{L} - \frac{\ell^{\prime}}{L^{\prime}} \geq 0$,
using the condition $a \geq b$, we have
$f_{a, b} (\ell; L) \leq f_{a, b} (\ell^{\prime}; L^{\prime})$.
We use this in place of \cite[(4.7)]{BRZ1} in our proof.

\item \cite[(4.8)]{BRZ1}:
By a simple calculation, for each $2 \leq \ell \leq L - 1$, we have
\begin{equation*}
(\ell - 1)_L^{\frac{1}{2} - \varepsilon} \leq 2^{\frac{1}{2} - \varepsilon}
\ell_L^{\frac{1}{2} - \varepsilon},~~~
\ell_L^{\frac{1}{2} - \varepsilon} \leq 2^{\frac{1}{2} - \varepsilon}
(\ell - 1)_{L - 1}^{\frac{1}{2} - \varepsilon}.
\end{equation*}
These inequalities imply the following:
\begin{equation*}
\sqrt{\frac{T_{\ell - 1} + T_{\ell}}{2}}
\geq \frac{1}{\sqrt{2}} \left\{f_{a^{\prime}, b} (\ell - 1; L - 1) - 2^{1 - 2\varepsilon} C 
(\ell - 1)_{L - 1}^{\frac{1}{2} - \varepsilon} \right\},
\end{equation*}
Thus, we can replace $\frac{\overline{C}}{2}$ above \cite[(4.9)]{BRZ1} with $2^{1 - 2\varepsilon} C$. 
(In particular, this does not depend on $\eta$.)
\item \cite[(4.16)]{BRZ1}:
By Lemma \ref{lem:barrier-estimate-0-dim-bessel-process},
we can replace the constant $C$ in \cite[(4.16)]{BRZ1}
with $c e^{\eta}$,
where $c > 0$ is independent of $\eta$.
\item \cite[(4.17)]{BRZ1}:
For any $j, k \geq 1$ with $|j - k| \leq 2 \eta L$ and $z \in H_j$,
we have
\begin{equation*}
(k - z)^2 = \left \{(k - j) + (j - z) \right \}^2
\geq (k - j)^2 + 2 (k - j)(j - z)
\geq (k - j)^2 - 4 \eta L.
\end{equation*}
Thus, we can replace $c$
(the first factor in the first line of \cite[(4.17)]{BRZ1})
with $c^{\prime} e^{5 \eta}$,
where $c^{\prime} > 0$ is independent of $\eta$.

\item \cite[(4.18)]{BRZ1}:
Let us look at the exponential factors in the first sum of \cite[(4.17)]{BRZ1}.
Since
\begin{equation*}
(k - j)^2 = \left\{(k - y) + (y - x) + (x - j) \right\}^2
\geq (y - x)^2 + 2 (y - x)(k - y) + 2 (y - x)(x - j),
\end{equation*}
we have
\begin{equation*}
e^{- \frac{(j - k)^2}{2 (L - 1)}}
\leq e^{- \frac{(y - x)^2}{2L}} e^{- \frac{y - x}{L - 1} (k - y)} e^{- \frac{y - x}{L - 1} (x - j)}.
\end{equation*}
Since $\sqrt{2} \leq x, y \leq \eta L$, we have
\begin{equation*}
c (x - j)^2 + \frac{y - x}{L - 1} (x - j)
= c \left\{(x - j) + \frac{y - x}{2c(L - 1)} \right\}^2
- \frac{(y - x)^2}{4c(L - 1)^2}
\geq c \alpha^2 - c_1 \eta^2,
\end{equation*}
where we set
$\alpha := (x - j) + \frac{y - x}{2c (L - 1)}$
and $c_1 > 0$ is independent of $\eta$.
($c_i$, $i = 2, \dotsc, 10$ below are also positive constants independent of $\eta$.)
Similarly, we have
\begin{equation*}
c (k - y)^2 + \frac{y - x}{L - 1} (k - y) \geq c \beta^2 - c_2 \eta^2,
\end{equation*}
where we set $\beta := (k - y) + \frac{y - x}{2c (L - 1)}$.
Thus, the product of three exponential factors in the first sum of \cite[(4.17)]{BRZ1}
is bounded from above by
\begin{equation*}
e^{(c_1 + c_2) \eta^2} e^{- c (\alpha^2 + \beta^2)} e^{- \frac{(y - x)^2}{2L}}.
\end{equation*} 
By the assumptions
$|a - b| \leq \eta L$ and $|x - y| \leq \eta L$,
we have
\begin{equation*}
1 + j - a^{\prime}
= 1 - \alpha + (x - a) + \frac{y - x}{2c (L - 1)} + \frac{a - b}{L}
\leq 1 + |\alpha| + (x - a) + c_3 \eta.
\end{equation*}
Similarly, we have
\begin{equation*}
1 + k - b \leq 1 + |\beta| + (y - b) + c_4 \eta.
\end{equation*}
Since $k \geq \frac{y + b}{2} \geq \frac{y}{2}$
and $\sqrt{2} \leq x, y \leq \eta L$,
we have
\begin{equation*}
\sqrt{\frac{j}{k}} = \sqrt{\frac{-\alpha + x + \frac{y - x}{2c (L - 1)}}{k}}
\leq \sqrt{\frac{2x}{y}} \sqrt{\frac{|\alpha|}{\sqrt{2}} + 1 + c_5 \eta}
\leq c_6 \sqrt{\frac{x}{y}} \left(\sqrt{|\alpha|} + \sqrt{\eta} \right).
\end{equation*}
By the above estimates,
the first sum in \cite[(4.17)]{BRZ1}
is bounded from above by
\begin{equation*}
c_7 e^{c_8 \eta^2} \frac{(\eta + x - a)(\eta + y - b)}{L^{3/2}}
\sqrt{\frac{x}{y}} e^{- \frac{(y - x)^2}{2L}}.
\end{equation*}
By the argument below \cite[(4.18)]{BRZ1},
the second sum in \cite[(4.17)]{BRZ1}
is bounded from above by
$c_9 e^{- c_{10} \eta^2 L^2}$.
\end{itemize}
By the above modifications,
we have the desired result.
$\qed$ \\\\
{\it Proof of Lemma \ref{lem:GW-process-barrier-estimate} (ii).}
We will use notation in \cite[Section 3]{BRZ1}.
Set
\begin{equation} \label{eq:constants-barrier-estimate-lower-bound}
C^{\prime} := 3 \cdot 2^{\frac{1}{2} - \varepsilon} C,~~~
C^{\prime \prime} := \left(\frac{1}{2} \right)^{\frac{3}{2} + \varepsilon} \widetilde{C}.
\end{equation}
Let $A$ be the event on the left hand side of (\ref{eq:GW-process-barrier-lower-bound}).
Set
\begin{equation*}
B:= 
\bigcap_{i=r}^{L-r}
\left\{
f_{a, 0} (i; L) + C^{\prime} i_L^{\frac{1}{2} - \varepsilon} \leq \sqrt{2 \mathcal{L}_i}
\leq f_{x, 0} (i; L) + C^{\prime \prime} i_L^{\frac{1}{2} + \varepsilon} \right\}
\cap
\left\{\mathcal{L}_L = 0 \right\}.
\end{equation*}
By \cite[Lemma 3.1 c)]{BRZ1}, 
$\mathbb{Q}_1^{\frac{x^2}{2}} [A \cap B]$
is equal to
\begin{equation} \label{eq:pf-GW-process-barrier-lower-bound-1}
\mathbb{Q}_1^{\frac{x^2}{2}} \left[1_B
\prod_{i = r}^{L - 1 - r}
\mathbb{Q}_1^{\frac{x^2}{2}} \left[
f_{a, 0} (i; L) + C i_L^{\frac{1}{2} - \varepsilon} \leq \sqrt{2 T_i}
\leq f_{x, 0} (i; L) + \widetilde{C} i_L^{\frac{1}{2} + \varepsilon}
| \mathcal{L}_i, \mathcal{L}_{i + 1} \right] \right].
\end{equation}
Fix $i \in \{r, \dotsc, L - 1 - r \}$.
Since $a \leq \eta L$,
we have
\begin{equation} \label{eq:pf-GW-process-barrier-lower-bound-2}
f_{a, 0} (i + 1; L) \geq f_{a, 0} (i; L) - \eta,~~~f_{x, 0} (i + 1; L) \leq f_{x, 0} (i; L).
\end{equation}
By a simple calculation, we have
\begin{equation} \label{eq:pf-GW-process-barrier-lower-bound-3}
(i + 1)_L^{\frac{1}{2} + \varepsilon} \leq 2^{\frac{1}{2} + \varepsilon} i_L^{\frac{1}{2} + \varepsilon},
~~~(i + 1)_L^{\frac{1}{2} - \varepsilon} \geq 2^{- \frac{1}{2} + \varepsilon} i_L^{\frac{1}{2} - \varepsilon}.
\end{equation}
Recall the definitions of $C^{\prime}$ and $C^{\prime \prime}$ from
(\ref{eq:constants-barrier-estimate-lower-bound}).
By (\ref{eq:pf-GW-process-barrier-lower-bound-2}),
(\ref{eq:pf-GW-process-barrier-lower-bound-3}),
and the assumption $C r^{\frac{1}{2} - \varepsilon} > \eta$,
under the event $B$, we have
\begin{equation} \label{eq:pf-GW-process-barrier-lower-bound-4}
f_{a, 0} (i; L) + 2C i_L^{\frac{1}{2} - \varepsilon}
\leq (4 \mathcal{L}_i \mathcal{L}_{i + 1})^{\frac{1}{4}}
\leq f_{x, 0} (i; L) + \frac{\widetilde{C}}{2} i_L^{\frac{1}{2} + \varepsilon}.
\end{equation}
By (\ref{eq:pf-GW-process-barrier-lower-bound-4}) and
\cite[Lemma 3.4, c)]{BRZ1}, under the event $B$,
the big product in (\ref{eq:pf-GW-process-barrier-lower-bound-1})
is bounded from below by a positive constant (not depending on $\eta$).
Thus, we have
\begin{equation} \label{eq:pf-GW-process-barrier-lower-bound-5}
\mathbb{Q}_1^{\frac{x^2}{2}} [A]
\geq \mathbb{Q}_1^{\frac{x^2}{2}} [A \cap B]
\geq c \cdot \mathbb{Q}_1^{\frac{x^2}{2}} [B]
= c \cdot \mathbb{Q}_1^{\frac{x^2}{2}} [B | \mathcal{L}_L = 0] \cdot
\mathbb{Q}_1^{\frac{x^2}{2}} [\mathcal{L}_L = 0].
\end{equation}
By \cite[Lemma 3.1, e)]{BRZ1}, 
$\mathbb{Q}_1^{\frac{x^2}{2}} [B | \mathcal{L}_L = 0]$ is bounded from below by
\begin{equation} 
\label{eq:pf-GW-process-barrier-lower-bound-6}
\begin{aligned}
P_{x \to 0}^Y \Biggl[
\left\{f_{a, 0} (s; L) + C^{\prime} (s_L)^{\frac{1}{2} - \varepsilon} 
\leq Y_s
\leq f_{x, 0} (s; L) + C^{\prime \prime} (s_L)^{\frac{1}{2} + \varepsilon},~
\forall s \in [r, ~L - r] \right\} \\
\cap \left\{Y_{s^{\prime}} \geq f_{a, 0} (s^{\prime}; L) - r^{\frac{1}{2} + 2 \varepsilon},~
\forall s^{\prime} \in [0, r] \right\}
\Biggr],
\end{aligned}
\end{equation}
where $P_{x \to 0}^{Y} [~\cdot~] := P_x^Y [~\cdot~ | Y_L = 0]$.
By an argument similar to the proof of \cite[Lemma 7.6]{BeKi},
(\ref{eq:pf-GW-process-barrier-lower-bound-6})
is bounded from below by 
\begin{equation} \label{eq:pf-GW-process-barrier-lower-bound-8}
c^{\prime} \frac{r}{L - 2r},
\end{equation}
where $c^{\prime} > 0$ is independent of $\eta$.
By (\ref{eq:pf-GW-process-barrier-lower-bound-5}),
(\ref{eq:pf-GW-process-barrier-lower-bound-8}),
and
$\mathbb{Q}_1^{\frac{x^2}{2}} [\mathcal{L}_L = 0] = (1 - \frac{1}{L})^{\frac{x^2}{2}}$,
we have the desired result. 
$\qed$
\\\\

{\bf Acknowledgements.} \\ 
This work was done almost when the author was
a postdoctoral researcher at Kobe University.
The author would like to thank Dr. Naotaka Kajino for valuable comments and encouragement.
The author is grateful to the referee for constructive suggestions.
This work was partially supported by JSPS KAKENHI Grant Numbers 16J00347, 18K13429.

\end{document}